\documentclass[openany, amssymb, psamsfonts]{amsart}
\usepackage{amssymb}
\usepackage{mathrsfs,comment}
\usepackage{stmaryrd}
\usepackage{tikz-cd}
\usepackage[normalem]{ulem}
\usepackage{url}
\usepackage{mathtools}
\usepackage{tikz}
\usetikzlibrary{calc,arrows}
\usepackage[all,arc,2cell]{xy}
\UseAllTwocells
\usepackage{enumitem}
\usepackage{hyperref} 
\usepackage[utf8]{inputenc}
\usepackage[T1]{fontenc}
\hypersetup{%
  bookmarksnumbered=true,%
  bookmarks=true,%
  colorlinks=true,%
  linkcolor=blue,%
  citecolor=blue,%
  filecolor=blue,%
  menucolor=blue,%
  pagecolor=blue,%
  urlcolor=blue,%
  pdfnewwindow=true,%
  pdfstartview=FitBH}

\def\makeautorefname#1#2{\expandafter\def\csname#1autorefname\endcsname{#2}}
\def\equationautorefname~#1\null{(#1)\null}
\makeautorefname{footnote}{footnote}%
\makeautorefname{item}{item}%
\makeautorefname{figure}{Figure}%
\makeautorefname{table}{Table}%
\makeautorefname{part}{Part}%
\makeautorefname{appendix}{Appendix}%
\makeautorefname{chapter}{Chapter}%
\makeautorefname{section}{Section}%
\makeautorefname{subsection}{Section}%
\makeautorefname{subsubsection}{Section}%
\makeautorefname{theorem}{Theorem}%
\makeautorefname{thm}{Theorem}%
\makeautorefname{cor}{Corollary}%
\makeautorefname{lem}{Lemma}%
\makeautorefname{prop}{Proposition}%
\makeautorefname{pro}{Property}
\makeautorefname{conj}{Conjecture}%
\makeautorefname{defn}{Definition}%
\makeautorefname{notn}{Notation}
\makeautorefname{notns}{Notations}
\makeautorefname{rem}{Remark}%
\makeautorefname{quest}{Question}%
\makeautorefname{exmp}{Example}%
\makeautorefname{ax}{Axiom}%
\makeautorefname{claim}{Claim}%
\makeautorefname{ass}{Assumption}%
\makeautorefname{asss}{Assumptions}%
\makeautorefname{con}{Construction}%
\makeautorefname{prob}{Problem}%
\makeautorefname{warn}{Warning}%
\makeautorefname{obs}{Observation}%
\makeautorefname{conv}{Convention}%

\newtheorem{thm}{Theorem}[section]
\newtheorem{cor}{Corollary}[section]
\newtheorem{prop}{Proposition}[section]
\newtheorem{lem}{Lemma}[section]

\newtheorem*{thm*}{Theorem}

\theoremstyle{definition}
\newtheorem{defn}{Definition}[section]

\newtheorem{rem}{Remark}[section]

\makeatletter
\let\c@obs=\c@thm
\let\c@cor=\c@thm
\let\c@prop=\c@thm
\let\c@lem=\c@thm
\let\c@prob=\c@thm
\let\c@con=\c@thm
\let\c@conj=\c@thm
\let\c@defn=\c@thm
\let\c@notn=\c@thm
\let\c@notns=\c@thm
\let\c@exmp=\c@thm
\let\c@ax=\c@thm
\let\c@pro=\c@thm
\let\c@ass=\c@thm
\let\c@warn=\c@thm
\let\c@rem=\c@thm
\let\c@sch=\c@thm
\let\c@equation\c@thm
\numberwithin{equation}{section}
\makeatother

\DeclareMathOperator{\SL}{SL}
\DeclareMathOperator{\GL}{GL}
\DeclareMathOperator{\Sp}{Sp}
\DeclareMathOperator{\SO}{SO}
\DeclareMathOperator{\SU}{SU}
\DeclareMathOperator{\Isom}{Isom}
\DeclareMathOperator{\injrad}{injrad}

\DeclareMathOperator{\Hess}{Hess}
\DeclareMathOperator{\Hom}{Hom}

\DeclareMathOperator{\diam}{diam}

\DeclareMathOperator{\ad}{ad}
\DeclareMathOperator{\Ad}{Ad}
\DeclareMathOperator{\tr}{tr}
\DeclareMathOperator{\rk}{rk}
\DeclareMathOperator{\grad}{grad}
\DeclareMathOperator{\vol}{vol}
\DeclareMathOperator{\absys}{absys}
\DeclareMathOperator{\supp}{supp}
\DeclareMathOperator{\FF}{FF}
\DeclareMathOperator{\FA}{FA}
\DeclareMathOperator{\EO}{E}

\DeclareMathOperator{\ff}{\mathfrak{f}}

\DeclareMathOperator{\Sym}{Sym}

\newcommand{\bs}{\backslash}

\newcommand{\fg}{\mathfrak{g}}      
\newcommand{\fk}{\mathfrak{k}}      
\newcommand{\fa}{\mathfrak{a}}  

\newcommand{\fn}{\mathfrak{n}}

\renewcommand{\setminus}{-} 
\makeatletter
\newtheorem*{rep@theorem}{\rep@title}
\newcommand{\newreptheorem}[2]{%
\newenvironment{rep#1}[1]{%
 \def\rep@title{#2 \ref{##1}}%
 \begin{rep@theorem}}%
 {\end{rep@theorem}}}
\makeatother

\newreptheorem{theorem}{Theorem}
\newreptheorem{cor}{Corollary}

\title{Minimal submanifolds and waists of locally symmetric spaces}

\author{Mikołaj Frączyk}
\address{Faculty of Mathematics and Computer Science, Jagiellonian University, ul. Łojasiewicza 6, 30-348 Krak{\'o}w, Poland}
\email{mikolaj.fraczyk@uj.edu.pl}

\author{Ben Lowe}
\address{Department of Mathematics\\University of Chicago\\Chicago, IL 60637}
\email{loweb24@uchicago.edu}

\date{\today}

\begin{document}

\maketitle

\begin{abstract} 
We study the higher expansion properties of locally symmetric spaces, with a particular focus on octonionic hyperbolic manifolds. We show that codimension two minimal submanifolds of compact octonionic locally symmetric spaces must have large volume, at least linear in the volume of the ambient space.  As a corollary we prove linear waist inequalities for octonionic hyperbolic manifolds in codimension two and construct the first locally symmetric examples of power-law systolic freedom.  We also show that any codimension two submanifold of small volume can be homotoped to a lower dimensional set. We use this to prove that branched covers of octonionic hyperbolic manifolds are stable in the sense of Dinur-Meshulam and to establish a uniform lower bound on the non-abelian Cheeger constants of octonionic hyperbolic manifolds. 

In a more general setting, we prove that maps from locally symmetric spaces to low dimensional euclidean spaces admit fibers whose fundamental group has large exponent of growth. We show as a consequence that cocompact lattices in $\SL_n(\mathbb R)$  have property ${\rm FA_{\lfloor n/8\rfloor-1}}$: any action on a contractible ${\rm CAT}(0)$ simplicial complex of dimension at most $ \lfloor n/8\rfloor -1$ has a global fixed point. 
\end{abstract} 

\tableofcontents
\section{Introduction}

Expander graphs have been widely useful in combinatorics,  number theory, group theory, and many other branches of mathematics and computer science \cite{lubotzkyEXP}. They have a natural Riemannian analogue in the form of manifolds with a uniform spectral gap for the Laplace-Beltrami operator. A family of closed Riemannian manifolds $M_n$ with a uniform spectral gap is called an \textit{expander family}, and has the following property: there is a uniform $\lambda>0$ so that for every $M_n$ and every continuous map $f:M_n \rightarrow \mathbb{R}$ with rectifiable fibers, there is some $p \in \mathbb{R}$ so that 
\begin{equation}\label{eq-topexpansion}
\vol (f^{-1}(p)) \geq \lambda \vol(M_n).  
\end{equation}
These types of estimates are called \emph{waist inequalities}. In a series of works \cite{gr03, gr09,gromov2010pt2} Gromov initiated a program to prove statements of the following kind.  Suppose we are given two simplicial complexes, $X$ and $Y$, where $X$ is ``complicated" and $Y$ is ``simple". Then, any map $f\colon X\to Y$ must have at least one ``complicated" fiber $f^{-1}(y), y\in Y$. Historically volume was the first measure of complexity considered, but it is just one of many possible choices. 

\begin{defn} Let $X,Y$ be  Riemannian manifolds (or simplicial complexes) of dimensions $n,d$ respectively. Let $f\colon X\to Y$ be a continuous map. The \emph{waist} of $f$ is 
 $w(f):=\sup_{y\in Y} \mathcal H^{n-d}(f^{-1}(y))$, where $\mathcal H^{n-d}$ stands for the $(n-d)$-dimensional Hausdorff measure (or the number of simplices intersecting the set). The \emph{$d$-waist} of $X$ is then defined as $w_d(X):=\inf_{f\colon X\to \mathbb R^d} w(f)$, where the infimum is taken over a suitable class of maps with $(n-d)$-rectifiable fibers.  
\end{defn}

The families of spaces $\{ X_i\}_{i\in I}$ satisfying $w_d(X_i)\geq c \vol(X_i)$ are called \emph{ topological $d$-expanders}. Gromov asked for examples of natural families of bounded geometry Riemannian manifolds or CW-complexes that are ``higher expanders", that is, $d$-expanders for $d\geq 2$. The topological expansion property subsequently motivated other definitions of higher expansion. In this paper we study these kinds of questions for locally symmetric spaces. While some of our results hold for general locally symmetric spaces, a large part of this work focuses on octonionic hyperbolic manifolds. This case seemed most approachable, being rank one and having good quantitative spectral gap estimates.  
 
We believe our work highlights a range of new possibilities for applications of minimal surface techniques to Gromov's program for non-positively curved locally symmetric spaces more generally, as well as to low dimensional actions of lattices in semisimple Lie groups. The interaction between representation theory and minimal surfaces is at the heart of our approach\footnote{Recent work by Song on  minimal surfaces in high dimensional spheres \cite{song2024random} also used representation theory, although in a completely different way.}. Our work draws on the min-max theory of minimal surfaces (via a connection to waists), which has seen a revival of activity (see \cite{marques2013applications} for a survey.) We hope there can be more applications of min-max theory to the study of the geometry of locally symmetric spaces and higher expansion.  

\subsection{Waist inequalities.} Our first theorem provides a natural class of $2$-expander manifolds, answering  the Riemannian version of Gromov's question in the $d=2$ case. 
\begin{thm}\label{mt-waist} There exists a constant $c>0$ such that any compact octonionic hyperbolic manifold $M$ satisfies $w_2(M)\geq c \vol(M).$ That is, compact octonionic hyperbolic manifolds form a topological $2$-expander family. 
\end{thm}
\begin{rem}
It is well known that if the fundamental group of a Riemannian manifold $M$ has Kazhdan's property $(T)$, then the family of finite covers of $M$ is an expander family. In a very recent work, Bader and Sauer \cite{bader2024} show that this family is also a topological $2$-expander thus giving another answer to Gromov's question in the Riemannian setting.
\end{rem}

\begin{rem}
The first examples of higher expanders came with the breakthrough work of Kaufman, Kazhdan and Lubotzky \cite{kaufman2014ramanujan}, who constructed simplicial $2$-complexes of bounded degrees which were topological $2$-expanders. This was later extended by Evra and Kaufman \cite{EK24} to a construction of bounded geometry $d$-expanders for any $d\geq 2$,  using quite refined combinatorial criteria for higher expansion. In both cases the spaces were derived from quotients of Bruhat-Tits building of simple $p$-adic groups. This settled Gromov's question in the setting of bounded geometry simplicial complexes.  Fox-Naor-Gromov-Lafforgue-Pach also studied higher expansion in the simplicial complex setting, in terms of overlap properties of maps of bounded degree hypergraphs to lower dimensional Euclidean spaces \cite{fox2012overlap}.    We also mention work by Zung \cite{zung2024expansion} and Abdurrahman-Adve-Giri-Lowe-Zung \cite{abdurrahman2024hyperbolic}  that established some higher expansion properties for families of respectively negatively curved and hyperbolic 3-manifolds.  
\end{rem}

The proof of Theorem \ref{mt-waist} is based on  new monotonicity formulas for minimal submanifolds of low codimension in rank one locally symmetric spaces.  We use them to show that codimension-two minimal submanifolds of octonionic hyperbolic manifolds have to be large, of volume comparable to the volume of the ambient space. By Almgren-Pitts theory, this implies linear lower bounds on waists and shows that octonionic hyperbolic manifolds form a topological $2$-expander family. Using similar methods, we prove that octonionic hyperbolic manifolds satisfy an even stronger condition: we show that they are homotopy $2$-expanders (see Definition \ref{def-hhdx}, Theorem \ref{mt-collapse}), which means that small codimension $2$ rectifiable sets can be efficiently homotoped into a lower dimensional set. This answers a question of Belolipetsky and Weinberger \cite{BelWein23} regarding the absolute systoles of locally symmetric spaces. We believe that our methods will also apply in the higher rank setting, see \ref{sec-IntroHR} for a discussion of some of the difficulties. 

Our methods differ significantly from the more algebraic approach of \cite{bader2024}. We are restricted to the setting of very specific locally symmetric spaces, while \cite{bader2024} applies to covers of any Riemannian manifold with fundamental group having Kazhdan's property (T). This in particular includes covers of any fixed compact hyperbolic octonionic manifold. On the flip side, we can deal with cycles which are not assumed to be boundaries and work homotopically (as opposed to homologically) which can be used to establish non-abelian expansion results like those described later in the introduction. 

\subsection{Topological Waist Inequalities} \label{subsec-topwaistineq}
In the previous section we used volume as the measure of complexity of a fiber. An alternative fruitful perspective is to look at the size of some measure of the complexity of the topology of the fibers, for example the size of the image of $\pi_1(f^{-1}(y))$ in $\pi_1(M)$ as $y\in Y$ varies. For rank one locally symmetric spaces we can prove the following.

\begin{thm} \label{intro:critical exponent}
Let $\Gamma$ be a uniform lattice in a rank one symmetric space $X$, and let $M = \Gamma \backslash X$. Then if $f: M \rightarrow \mathbb{R}^d$ is a generic piecewise-linear map, there is a constant $C_X(d)>0$ depending only on $X$ and $d$ so that for some $p\in \mathbb{R}^d$ the fundamental group of the inclusion of some connected component of $f^{-1}(p)$ in $\pi_1(M)$ has critical exponent at least $C_X(d)$.  
 \end{thm}

\begin{rem}
The constant $C_X(d)$ can be computed in a straightforward way from the symmetric space, see Definition \ref{defn-CXk} and Table \ref{tab-cxd}. For example, for real hyperbolic space $X=\mathbb{H}^n$, $C_X(d)=n-d-1$. 
\end{rem}

The critical exponent of a discrete subgroup of a rank one group measures the asymptotic growth rate of orbits of that group.  The analogue in higher rank is a function from the positive Weyl chamber $\mathfrak{a}^+$ to $\mathbb{R}$ called the growth indicator function.  We prove a version of Theorem \ref{intro:critical exponent} for general nonpositively curved locally symmetric spaces that includes it as a special case and that, for a generic PL map, gives a lower bound on the growth indicator function of some fiber (Theorem \ref{largefiber}.)  Combined with the upper bounds on the growth indicator functions of infinite covolume subgroups \cite{LeeOh24}, \cite{oh02},  this gives the following corollary. We show that a generic map from $M$ to a low-dimensional $\mathbb{R}^d$, $d \leq r(X)$, has a fiber that includes to a finite index subgroup of $\pi_1(M)$. The constant $r(X)$ can in principle be computed from the root system of $G$, which we do in several cases. This is a topological analogue of a linear waist inequality.

\begin{thm}\label{mt-topwaist} Let $X$ be a symmetric space of non-compact type. For any compact quotient $\Gamma\backslash X$ and generic piecewise-linear map $f\colon \Gamma\bs X\to \mathbb{R}^d$, $d\leq r(X)$, we can find a fiber $f^{-1}(y)$ such that $\pi_1(f^{-1}(y))$ includes to a finite index subgroup of $\Gamma$.
\end{thm}

\begin{rem}
Lower bounds for the constants $r(X)$ are given in Propositions \ref{prop-octonionicr(x)}, \ref{prop-novereight}.  In particular,  $r(\mathbb H^2_\mathbb O)\geq 2$ and $r(\SL_n(\mathbb R)/\SO(n))\geq \lfloor \frac{n}{8}\rfloor-1.$  In contrast to Theorem \ref{mt-waist} we do not need to restrict to rank one spaces.  
\end{rem}

\begin{rem}
In formulating waist inequalities, it is necessary to restrict to a class of maps for which the relevant measure of complexity makes sense for the fibers.  One can accomplish this, and this is the approach we take here, by considering maps that are generic in an appropriate sense. This is also the approach taken in the paper by Freedman that proves various topological waist inequalities for topologically stable maps from a smooth manifold $M \rightarrow \mathbb{R}^d$  \cite{freedmanwidths}, and the papers by Gromov \cite{gr09},\cite{gromov2010pt2} mentioned above. One alternative approach is to work with real analytic maps as Bader-Sauer in their recent work do \cite{bader2024}. We comment that the statement of Theorem \ref{mt-topwaist} also holds for generic smooth maps provided that generic smooth maps are stratified fibrations (see \cite{gr09}[Section 2.4,3] for statements in this direction), but we choose to work in the PL category because we found it most convenient for the applications.  
\end{rem}

The proof of Theorem \ref{mt-topwaist} depends on recent work by Connell-McReynolds-Wang, which we describe in more detail further down \cite{cmw23}.  For higher rank and certain hyperbolic infinite volume locally symmetric spaces, they construct vector fields whose flows contract volume in low codimension.  One application is to prove vanishing theorems for cohomology in low codimension.    

A result similar to Theorem \ref{mt-topwaist} was obtained by Alagalingam for tori \cite{alagalingam2019algebraic}.  Alagalingam proved that for a generic map from $\mathbb{T}^n$ to $\mathbb{R}^d$ there is some $p\in \mathbb{R}^d$ so that the image of  some connected component of $\pi_1(f^{-1}(\{p\}))$ under the map induced by inclusion has rank at least $n-d$ in $\pi_1(\mathbb{T}^n) = \mathbb{Z}^n$. Using the same strategy as the proof of Theorem \ref{mt-topwaist}, it is possible to give an alternative proof of Alagalingam's result (see Section \ref{subsec-ideasoftheproofs} for a more detailed explanation.) 

Using Theorem \ref{mt-topwaist} we are able to prove the following theorem, which leads into the next section.   

\begin{thm}\label{mt-finiteorbit} Let $\Gamma$ be a lattice in $\Isom(X)$ and suppose that $\Gamma$ acts by simplicial automorphisms on a contractible $d$-dimensional simplical complex $\Delta$, for $d\leq r(X).$ Then the action of $\Gamma$ has a finite orbit. 
\end{thm}

\subsection{Property $\FA_r$ and low dimensional actions on $\rm{CAT}(0)$ spaces.}
Theorem \ref{mt-finiteorbit} is closely related to property $\FA_r$, which we now describe.  A group has Serre's Property $\FA$ if every action of that group on a tree has a global fixed point.  Property $\FA$ is a central notion in Bass-Serre theory, and it puts strong constraints on the structure of a group that possesses it.  For example, a group with property $\FA$ cannot split as an amalgamated free product or HNN extension, and its irreducible representations in $\GL_2(\mathbb{C})$ must have algebraic integer traces. It is well known that property FA is implied by Kazhdan's property (T), so with the exception of $G=\SO(n,1),\SU(n,1)$ (and isogeneous groups), lattices in any simple real Lie group $G$ have property FA. In particular, all lattices in higher rank simple groups enjoy property $\FA$. This was first proved by Serre for $\SL(3,\mathbb{Z})$ \cite{serre02}, and for higher rank lattices in general by Margulis \cite{margulis1981decomposition}.  

On the other hand, higher rank lattices \textit{do} admit natural isometric actions on simply connected higher dimensional simplicial complexes without global fixed points, for example the action on the symmetric space of the ambient semi-simple Lie group. Isometric actions of higher rank lattices on non-positively curved contractible manifolds are super-rigid by \cite{MokSiuYeung} so they are quite well understood. However, the actions on simplicial complexes (either contractible or non-positively curved and simply connected) are still mysterious. The lattices in $p$-adic simple Lie groups like $\SL_n(\mathbb Q_p)$ act on their Bruhat-Tits complex, which is a non-positively curved simplicial complex of dimension equal to the rank of the group. We currently do not know of an example of a lattice in a simple (real or $p$-adic) group $G$ that acts on a non-positively curved simplicial complex of dimension less than the rank of $G$ without a global fixed point.

\begin{defn}
A group $\Lambda$ has property $\text{FA}_r, r>0$ if every isometric action by simplicial automorphisms on a $d$-dimensional contractible $\text{CAT}(0)$-simplicial complex by isometries, $d\leq r$, has a fixed point. 
\end{defn}

The action of a group on a simply connected $d$-dimensional simplicial complex by simplicial automorphisms without a global fixed point gives a presentation of the group as a pushout of a $d$-dimensional complex of groups \cite{soule1977chevalley}. Complexes of groups are a higher dimensional generalization of graphs of groups, see \cite{corson92}, \cite{haefliger1991complexes}.For $d=1$, this is the fact that an action of a group on a tree without a fixed point gives a splitting of the group as an amalgamated free product or HNN extension. Stallings realized that it was necessary to add a non-positive curvature condition on the simplicial complexes under consideration to extend the theory of group actions on trees to actions on higher dimensional complexes \cite{stallings1991non}. Asking whether a group has property $\FA_r$ is thus equivalent to asking whether a given group can be written nontrivially as a nonpositively curved $r$-dimensional complex of groups.

Property $\FA_r$ was introduced by Farb in \cite{fa09}, where he proved $\SL_n(\mathbb Z)$ has property $\text{FA}_{n-2}.$  Besides $\SL_n(\mathbb Z)$, Farb was able to prove property $\FA_r$ for a range of $r>1$ for non-uniform higher rank lattices, $S$-arithmetic groups, Chevalley groups, and certain discrete subgroups of $\SO(n,1)$ generated by reflections.  Subsequently much work has been done to prove similar lower bounds over broad classes of discrete groups (see e.g. \cite{bridson2010mappingclass}, \cite{varghese2011fixed}, \cite{barnhill2006fan}, \cite{ye2014low}.) 

Farb's proof of  property $\text{FA}_{n-2}$ for nonuniform lattices in $\SL_n(\mathbb{R})$ relied on the presence of unipotent subgroups, however, and to the authors' knowledge no progress has been made towards establishing property $\FA_r$ for $r>1$ in the uniform case.  Our Theorem \ref{mt-finiteorbit} implies:

\begin{cor} \label{corfanover8}
\begin{enumerate}\item Cocompact lattices in $\SL_n(\mathbb R)$ have property $\FA_{\lfloor \frac{n}{8}\rfloor -1}.$
\item Cocompact lattices in $\rm{F}_4^{(-20)}$ have property $\FA_{2}.$
\end{enumerate}
\end{cor}

Fixed point theorems for groups acting on non-positively curved spaces have also been obtained in a variety of different contexts using harmonic maps and ergodic theory (\cite{gromov1992harmonic}, \cite{margulis1991discrete} \cite{gao1997superrigidity}, \cite{burger1996cat}.)  Farb's proofs are based on facts about the combinatorics of intersections of convex sets.  By contrast, Corollary \ref{corfanover8} follows from a topological waist inequality (Theorem \ref{largefiber}.)   

By a recent work of Bader and Sauer \cite{bader2023higher}, the groups $\SL_n(\mathbb R)$ enjoy property $(T_{n-2})$, called higher property $(T)$ for $n\geq 4$. Currently the relation between higher property $(T)$ and other higher expansion phenomena is not clear. Motivated by the relation between property (T) and $\FA$ one could conjecture that property $(T_n)$ should imply $(\FA_n)$. If that were true, cocompact lattices in $\SL_n(\mathbb R)$ should have $\FA_{n-2}$. Our method doesn't go that far but we still establish $\FA_r$ for $r$ linear in the conjectural range.

There are some parallels between our work and work restricting the actions of lattices by diffeomorphisms on manifolds. According to Zimmer's conjectures, which after recent breakthroughs have been established in a number of important cases \cite{brown2020zimmer,brown2022zimmer},  every action of a lattice in a Lie group on a manifold of sufficently low dimension should factor through the action of a finite group.  Property $\FA_r$ is similar, in that it restricts the actions of lattices on low-dimensional spaces, which are allowed to belong to a very general class (they are not required to be locally compact) but with the added conditions that the action be by simplicial automorphisms and that the space is contractible.  



\subsection{Minimal submanifolds of locally symmetric spaces}
Let $M$ be a Riemannian manifold. A submanifold $S\subset M$ of dimension $m$ is \emph{minimal} if it is a critical point of the $m$-dimensional volume functional. In other words, for every differentiable deformation $S_t, t\in [0,\varepsilon)$  with $S_0=S$ we have 
$$\frac{\partial}{\partial t}\mathcal H^m(S_t)|_{t=0}=0,$$ where $\mathcal H^m$ stands for the $m$-dimensional Hausdorff measure \cite[\S 1.2]{simongmt}.
This definition can be extended to stationary varifolds \cite[Chap. 4. \S 2.4]{simongmt}, which are a well behaved and useful generalization of minimal submanifolds. Minimal surfaces in hyperbolic $3$-manifolds have been studied and successfully used as a tool to study the ambient manifold \cite{calegari2006shrinkwrapping,gabai1997geometric,  Lackenby2002HeegaardST}. In contrast, the topic of low codimension minimal submanifolds of other locally symmetric spaces remains mostly unexplored.
\begin{thm}\label{mt-stationary}
Let $M$ be a compact  octonionic hyperbolic $16$-manifold. There is a constant $c>0$ such that any stationary integral rectifiable $(16-d)$-varifold $S\subset M$ satisfies $\mathcal{H}^{16-d}(S)\geq c \vol(M)$ for $d=1,2$. 
\end{thm}

This phenomenon is new and does not occur in real or complex hyperbolic spaces\footnote{In real hyperbolic spaces one can find small codimension one closed totally geodesic hypersurfaces. For complex hyperbolic case consider a complex hyperbolic manifold $M$ admitting a non-constant holomorphic map to a closed surface $S$ and take a sequence of covers of $M_n$ corresponding to some $\mathbb Z/n\mathbb Z$- cyclic covers $S_n$ of $S$. Then the lift of a short closed geodesic on $S_n$ to $M_n$ can be deformed to a small minimal stationary varifold of codimension $1$ (it cannot collapse to smaller dimension because it is homologically non-trivial).}. It is also special to low codimensions. 
While we believe the analogous results should hold for higher rank locally symmetric spaces, extending our methods to the higher rank setting will require several new ideas. Theorem \ref{mt-stationary} is a crucial ingredient in the proof of Theorem \ref{mt-waist}. 
As part of the proof of Theorem \ref{mt-stationary}, we improve Anderson's monotonicity estimates \cite{anderson1982complete} for negatively curved symmetric spaces. Define $\kappa(k)$ as in the Table \ref{tab-kappaintro}.
\begin{table}[ht]
\centering
\begin{tabular}{ l l |c |c }
 & & $X$ & $\dim X$\\ \hline
 $\kappa(k)=k-1$ & & $\mathbb H_\mathbb R^n$ & $n$\\
 $\kappa(k)=\begin{cases}k-1 \\ 2n \end{cases}$ & $\begin{array}{l} k=1,\ldots, 2n-1 \\  k=2n \end{array}$ & $\mathbb H_\mathbb C^n$ & $2n$\\
 $\kappa(k)=\begin{cases}k-1\\ 4n-4+2i \end{cases}$ & $\begin{array}{l} k=1,\ldots, 4n-3 \\ k=4n-3+i  \end{array}$ &$\mathbb H_\mathbb H^n$ & $4n$\\
$\kappa(k)=\begin{cases}k-1\\ 8+2i\end{cases}$ & $\begin{array}{l} k=1,\ldots,9 \\ k=9+i  \end{array}$ & $\mathbb H_\mathbb O^2$ & $16$\\
\end{tabular}
\caption{$\kappa$ in rank one groups.}
\label{tab-kappaintro}
\end{table}

\begin{thm}\label{mt-monotonicity}
Let $S\subset X, X=H^n_{\mathbb R}, \mathbb H^n_{\mathbb C}, \mathbb H^n_{\mathbb H}, \mathbb H^2_{\mathbb O},$ be a stationary integral $k$-varifold. Normalize the metric so that the highest sectional curvature of $X$ is $-1$. Then 
\begin{enumerate}
    \item $\mathcal H^k(S\cap B(x,r))\geq e^{\kappa(k)(r-s)}\mathcal H^k(S\cap B(x,s))$ for $r>s$,
    \item $\mathcal H^k(S\cap B(x,r))\gtrsim  e^{\kappa(k)r}$  for any regular point $x\in S,$
\end{enumerate} where $B(x,r)$ is the $r$-ball around $x$.
\end{thm}
Anderson's estimate \cite{anderson1982complete} yields  $\mathcal H^k(S\cap B(x,r))\gtrsim e^{(k-1)r}$. We gain better exponents for $X=\mathbb H^n_{\mathbb C}, \mathbb H^n_{\mathbb H}, \mathbb H^2_{\mathbb O}$ and $k$ close to the top dimension. This gain is most significant for the octonionic hyperbolic plane, which is part of the reason why Theorems  \ref{mt-waist} and \ref{mt-stationary} are restricted to octonionic hyperbolic manifolds.

\subsection{Absolute systoles and homotopy higher expanders}\label{sec-hdx}
Let $M$ be a Riemannian manifold. The $k$-th absolute systole of $M$, denoted $\absys_k(M)$ is the infimum of $k$-dimensional Hausdorff measures of rectifiable subsets of $M$ which cannot be homotoped to a $(k-1)$-dimensional subset \cite{BelWein23}. This quantity was studied by Belolipetsky and Weinberger for arithmetic locally symmetric spaces in \cite{BelWein23}. For each symmetric space $X$ they show that the absolute $k$-systole $\Gamma\bs X$ grows at least as fast as a power of volume when $k$ is strictly above the rank. We strengthen this result in octonionic hyperbolic spaces for dimensions near the top.

\begin{thm}\label{mt-collapse}
Let $\varepsilon>0$ and let $M$ be a compact octonionic hyperbolic manifold with injectivity radius at least $\varepsilon.$ There exists an $\alpha=\alpha(\varepsilon)>0$ with the following property. Any rectifiable $k$-varifold $S$  with $\mathcal H^{k}(S)\leq \alpha \vol(M)$ and $k=14,15$ admits a homotopy $h\colon S\times [0,1]\to M$ such that
\begin{enumerate}
\item $h(s,0)=s$ for all $s\in S$,
\item $\mathcal H^{k-1}(h(S,1)))<\infty,$
\item $\mathcal H^{k+1}(h(S,[0,1]))\leq c \mathcal H^k(S),$ with $c=c(\varepsilon).$
\end{enumerate}
\end{thm}

\begin{rem}
A consequence of theorem \ref{mt-collapse} is that every integral cycle $F$ in codimension one or two with area smaller than some constant  times the volume of the ambient space can be homotoped to a lower-dimensional set.  This simpler statement admits a short proof using Theorem \ref{mt-stationary},  work of White \cite{white1984mappings}, and the natural flows from \cite{cmw23}.  See Section \ref{simplerproofofmt-collapse} for a sketch of the argument. 
\end{rem}

\noindent Theorem \ref{mt-collapse} implies that any small codimension $2$ submanifold of an octonionic hyperbolic manifold can be deformed to a codimension-$3$ subset. Theorem \ref{mt-collapse} should be compared with a similar result for infinite volume quotients \cite{cmw23}. In infinite volume spaces $\Lambda\bs X$ Connell, McReynolds and Wang construct flows which squeeze the volumes of $k$-submanifolds for $k$-sufficiently close to $\dim X$. Such flows cannot exists on compact quotients so our argument follows a different approach, but we do share some common ingredients. 

The conclusions of Theorem \ref{mt-collapse} are reminiscent of the definition of cosystolic expanders, introduced independently by Gromov \cite{gromov2010pt2} and Linial-Meshulam \cite{linialHDX}. We recall the definition below, following \cite{EK24}.
Let $Y$ be a finite $n$-dimensional simplicial complex. Let $\delta\colon C^{\bullet}(Y,\mathbb F_2)\to C^{\bullet+1}(Y,\mathbb F_2)$ be the standard coboundary map on the simplicial cochain complex of $Y$ with coefficients in $\mathbb F_2$. We endow $C^i(Y,\mathbb F_2)$ with a norm $\|a\|$ which, up to a universal constant, measures the size of the support of the cocycle $a\in C^i(Y,\mathbb F_2)$ divided by the total number of $i$-cells. Denote by $Z^i(Y,\mathbb F_2)$ and $B^i(Y,\mathbb F_2)$ the cocycles and coboundaries respectively. The simplicial complex $Y$ is a \emph{$(\varepsilon,\mu)$-cosystolic $k$-expander} if 
\begin{align*}
\|\delta f\| \geq& \varepsilon \|f+z\| &\text{for all }f\in C^k(Y,\mathbb F_2)\setminus Z^k(Y,\mathbb F_2) \text{ and } z\in  Z^k(Y,\mathbb F_2),\\
\|z\|\geq& \mu &\text{for all }z\in Z^k(Y,\mathbb F_2)\setminus B^k(Y,\mathbb F_2).
\end{align*}
Being a cosystolic expander means that small coboundaries must have small fillings and small cocycles must be coboundaries. Theorem \ref{mt-collapse} implies a dual statement for cycles. Small closed cycles should be boundaries and small boundaries should admit small fillings. Theorem \ref{mt-collapse} motivates the following definition of a homotopy expander. We use the language of integral rectifiable varifolds, explained in Section \ref{sec-varifolds}. The $\mathcal H^{m}$ stands for the $m$-dimensional Hausdorff measure \cite[\S1.2]{simongmt}. Readers unfamiliar with varifolds can replace them by $C^1$-submanifolds and $\mathcal H^m$ by the $m$-dimensional Riemannian volume.  

\begin{defn}\label{def-hhdx} Let $\{M_i\}$ be a family of bounded geometry Riemannian manifolds. We say that $\{M_i\}$ is a homotopy $d$-expander family if there are constants $\alpha,c>0$ for which the following holds. Any integral rectifiable $(n-d)$-varifold $S\subset M_i$, with ${\mathcal H}^{n-d}(S)\leq \alpha \vol(M_i)$ admits a homotopy $h\colon S\times [0,1]\to M_i$ with $h(s,0)=s$ for all $s\in S$, $\mathcal H^{n-d-1}(h(S,1))<\infty$ and $\mathcal H^{n-d+1}(h(S,[0,1]))\leq c \mathcal H^{n-d}(S).$
\end{defn}

Theorem \ref{mt-collapse} says that compact octonionic hyperbolic manifold with injectivity radius at least $\varepsilon$ form a homotopy $2$-expander family. The deformations responsible for Theorem \ref{mt-collapse} are explicit enough to derive consequences on branched cover stability, explained in section \ref{sec-introbranched}. 

\subsection{Power-Law Systolic Freedom}

Loewner\footnote{The result was first mentioned in print in the paper by Pu \cite{pu1952some} proving a systolic inequality for $\mathbb{RP}^2$, where it is credited to Loewner.} proved that any unit area Riemannian torus has a non-contractible loop of length at most $\sqrt{\frac{2}{\sqrt{3}}}$.  More generally, Gromov proved that for any unit volume closed aspherical Riemannian manifold $M$, there is an upper bound on the length of the smallest non-contractible loop in $M$ depending only on the dimension of $M$ \cite{gromov1983filling}.   This result is an example of a \textit{systolic inequality}, or inequality comparing the volume and the \textit{systole}, or infimal k-volume of a homologically or homotopically non-trival k-cycle. As we have seen in the previous section, in certain geometric situations it makes sense to ask for lower bounds on systoles in terms of volume. 
The $(p,q)$-systolic ratio $\text{(p,q)-SR}(M)$ of a closed $(p+q)$-dimensional Riemannian manifold $M$ is defined to be 
\[
\text{(p,q)-SR}(M):= \inf_{P,Q  } \frac{\text{Vol}_p (P) \text{Vol}_q(Q)}{\text{Vol}(M)}, 
\]
\noindent where the infimum is taken over all rectifiable cycles $P$ and $Q$ that are non-trivial in respectively $H_p(M;R)$, $H_q(M;R)$ for a coefficients ring $R$.
We say that a sequence of bounded geometry Riemannian manifolds $M_n$ exhibits \textit{systolic freedom} if $\text{(p,q)-SR}(M_n) \rightarrow  \infty $ for some $p$ and $q$ and that $\{M_n\}$ exhibits \textit{power-law systolic freedom} if $\text{(p,q)-SR}(M_n) \gtrapprox \text{Vol}(M_n)^\beta$ for some $p$ and $q$, and $\beta>0$.  The study of systolic freedom was pioneered by Gromov (see \cite{gromovintersystolic} for a survey.)  

It follows from the results of the previous paragraph that there cannot be systolic freedom in two dimensional Riemannian manifolds.  In contrast, there are examples in three or greater dimensions of systolic freedom. Gromov gave the first examples of systolic freedom, which were homeomorphic to $S^3 \times S^1$ \cite{gromovintersystolic}. Freedman later gave examples in three dimensions with mod-2 coefficients \cite{freedman1999z2}.  Both of these came as a surprise at the time.   Freedman-Hastings gave the first examples of Riemannian manifolds exhibiting power-law systolic freedom for $R=\mathbb{F}_2$ \cite{freedmanhastings}. Their construction has a combinatorial flavor; the Riemannian manifolds they construct are built from $\mathbb{F}_2$-chain complexes associated to quantum codes \footnote{Freedman-Hastings make the following comment regarding their examples: ``In retrospect it is not surprising that the most efficient constructions will necessarily have a combinatorial complexity beyond a geometer’s intuition, perhaps even requiring computer search, and are
better discovered within coding theory and then translated into geometry."}.  We are able to give a very different class of examples exhibiting power-law systolic freedom.

\begin{thm}\label{mt-sysfreedom}
Any sequence of congruence covers $M_n$ of a fixed closed octonionic hyperbolic manifold exhibits power-law systolic freedom, for any ring $R$.    
\end{thm}
\begin{proof}
 \noindent  We claim that for $M_n$ as in the statement of the theorem
\[
\text{(2,14)-SR}(M_n) \gtrsim\text{Vol}(M_n)^\beta
\]

for some $\beta>0$.  We know by Theorem \ref{mt-collapse} that the 14-volume of any homologically non-trivial 14-cycle in $M_n$ is $ \gtrsim \text{Vol}(M_n)$.  On the other hand, it follows from \cite[Theorem A]{BelWein23} that the 2-volume of any non-trivial 2-cycle is $ \gtrsim \text{Vol}(M_n)^\beta$ for some $\beta>0$.  

\end{proof}

We comment that work by Guth and Lubotzky also used techniques from minimal surface theory, for example the Anderson monotonicity formula for minimal surfaces in hyperbolic manifolds \cite{anderson1982complete}, along with techniques from systolic geometry, to construct new quantum error correcting codes from arithmetic hyperbolic 4-manifolds \cite{guth2014quantum}.

\subsection{Branched cover stability and non-abelian expansion}\label{sec-introbranched}

We are also able to prove results to the effect that branched covers with small branching locus must actually have been genuine covers in disguise.  

\begin{thm}\label{mt-branched} Let $\varepsilon>0$. There exists $\alpha, c>0$ dependent only on $\varepsilon$ with the following properties. Let $M$ be a compact octonionic hyperbolic manifold of systole at least $\varepsilon$ and let $\iota_1\colon M_1\to M$ be a branched cover of $M$ with branching locus $N\subset M$, which is a closed rectifiable $14$-varifold. Suppose, $\mathcal H^{14}(N)\leq \alpha \vol(M).$ Then, there exists a $15$-dimensional rectifiable varifold $W\subset M$ with $\mathcal H^{15}(W)\leq c\mathcal H^{14}(N)$, a genuine cover $\iota_2\colon M_2\to M$ and an isometry $\mu\colon M_1\setminus \iota_1^{-1}(W)\to M_2\setminus \iota_2^{-1}(W)$ such that $\iota_2\circ \mu=\iota_1.$ 
\end{thm}

This theorem implies a corresponding combinatorial statement for triangulations of octonionic hyperbolic manifolds (Theorem \ref{thm-CovStab}). In the language of \cite{dinur2022near} we prove that compact octonionic hyperbolic manifolds of systole at least $\varepsilon>0$ admit triangulations such that the resulting family of simplicial complexes is a non-abelian cosystolic expander family. We refer to sections \ref{sec-branched} and \ref{sec-ccstability} for a brief introduction to non-abelian Cheeger constants and the cover stability setup of \cite{dinur2022near}.
\begin{thm}\label{mt-nacheeger}
Let $\Gamma\subset F_4^{(-20)}$ be a cocompact lattice. Let $\mathcal T$ be a triangulation of $\Gamma\bs \mathbb H^2_\mathbb O$ and let $Y$ be the resulting simplicial complex. There is a $\delta>0$ such that the non-abelian Cheeger constants satisfy
$h_1(Y_1,H)\geq \delta$, for any finite cover $Y_1$ of $Y$ and any group $H$. 
\end{thm}

We note that in a recent work \cite{Dikstein2023CoboundaryAC} Dikstein and Dinur proved that families of bounded degree simplicial complexes satisfying certain expansion conditions for the links are cover stable. This is a consequence of isoperimetric inequalities which they establish for 1-cochains with values in non-abelian groups. This approach is somewhat dual to one we took to prove Theorem \ref{mt-branched}, in the sense that we work with subvarifolds which are more like chains as opposed to cochains. Interestingly, for the proof both \cite{EK24} and \cite{Dikstein2023CoboundaryAC} consider \emph{locally minimal co-chains} which raises the question of whether there is some deeper connection between the methods of \cite{EK24,Dikstein2023CoboundaryAC} and our approach based on minimal submanifolds. 
Cover stability theorems for a simplicial complex are related to the problem of flexible stability of its fundamental group. The type of expansion we prove is not sufficient for flexible stability, the two are compared in \cite[\S 4.4]{chapman2023stability}. It is an interesting question whether our methods can be adapted to show that octonionic hyperbolic lattices are flexibly stable.
\subsection{Ideas of the proofs} \label{subsec-ideasoftheproofs}
The proof of \textbf{Theorem \ref{mt-stationary}} is inspired by Anderson's monotonicity estimates in negatively curved manifolds \cite{anderson1982complete}, the flow version of Besson-Courtois-Gallot maps \cite{bcgamalgamated} developed recently by Connell, McReynolds and Wang \cite{cmw23} and Corlette's gap theorem \cite{corlette1990}. Anderson's monotonicity estimate predicts that a minimal $k$-submanifold $S$ in a contractible space $X$ with sectional curvature $\leq -1$ will have volume growth at least that of the real hyperbolic $k$-space $\mathbb H^k_\mathbb R$. More precisely, for any point $x_0\in S$ 
$${\mathcal H^k}(S\cap B_X(x_0,r))\gtrsim \vol(B_{\mathbb H^k_\mathbb R}(r))\gtrsim e^{(k-1)r}.$$
For $X=\mathbb H^n_\mathbb R, n\geq k$ this estimate is sharp, as one can take $S$ to be a totally geodesic copy of $\mathbb H^k_\mathbb R.$ If we know that $X$ has varying negative curvature, this estimate can be improved in low codimensions. The idea is that if $k$ is near the top dimension, the $k$-submanifolds should witness the directions with more negative curvature. To make this precise, we study the Hessians of distance functions, and establish lower bounds for the sum of their bottom $k$-eigenvalues. The same type of estimates were required in \cite{cmw23}, but for Busemann functions instead of distances. In this way we get improved monotonicity estimate in low codimensions for $X=\mathbb H^n_\mathbb C, \mathbb H^n_\mathbb H, \mathbb H^2_\mathbb O$ (Theorem \ref{mt-monotonicity}). 

\begin{rem}
Yuping Ruan recently found and corrected a mistake in the work of Besson-Courtois-Gallot in the octonionic case \cite{ruan2024cayley}. He shows that the formula for the Hessian of Busemann functions they write in the octonionic case in terms of complex structures $J_1,..,J_7$ is incorrect. In fact, he shows that complex structures that behave well with respect to the curvature tensor cannot exist in the octonionic case, unlike in the complex and quaternionic hyperbolic cases.  Our computations of Hessians in Section \ref{sec-Hess} are carried out by performing computations in the Lie algebra and make no reference to complex structures, so they are unaffected by this issue.  
\end{rem}

One interesting observation, which motivated our approach, is that for $X=\mathbb H^2_{\mathbb O}$ (which is $16$-dimensional)  and $k=14,15$ the rate of volume growth of a minimal submanifold of $X$ would have to be faster than the maximal possible rate of growth for an infinite covolume subgroup of $\Isom(\mathbb H^2_\mathbb O)$ \cite{corlette1990}. This indicates a tension between monotonicity estimates and the rate of decay of matrix coefficients which stood behind Corlette's growth gap.

Suppose that $S$ is a minimal $k$-submanifold of a finite volume locally symmetric space $\Gamma\bs \mathbb H^2_\mathbb O$. Let $G=\Isom(\mathbb H^2_\mathbb O)$, and let $K$ be a maximal compact subgroup so that $\mathbb H^2_\mathbb O=G/K$. The monotonicity estimate gives us a lower bound on the mass of $S$ in $r$-balls around points $s\in S$. On the other hand, we can access the average mass spectrally, by considering integrals of certain matrix coefficients (Lemma \ref{lem-BallIntForm}). First we thicken $S$ and consider its smoothened characteristic function $\mathfrak u\in L^2(\Gamma\bs \mathbb H^2_\mathbb O)=L^2(\Gamma\bs G)^K$. The average mass of an $r$-ball around a point of $S$ can be read from the matrix coefficient $\Psi(g):= \langle  \mathfrak ug,\mathfrak u\rangle, g\in \Isom(\mathbb H^2_\mathbb O)$ (formula \ref{eq-PsiIntegralformula}). The proof of Theorem \ref{mt-stationary} is based on the tension between the decay of matrix coefficients of the group $G$ (Lemma \ref{lem-rankonedecay}), which up to isogeny is the exceptional rank one group $F_4^{(-20)}$, and the lower bound on growth of balls coming from the monotonicity estimate. 
We show that the lower bound on volume growth and the upper bounds on matrix coefficients can be reconciled only if $\mathfrak u$ has large projection onto the space of constant functions in $L^2(\Gamma\bs G).$ The latter translates to a lower bound on $\mathcal H^k(S)$ which is linear in $\vol(\Gamma\bs \mathbb H^2_\mathbb O).$

The proof of \textbf{Theorem \ref{mt-collapse}} expands upon the ideas from the argument explained above. Here $S$ is no longer assumed to be stationary, but we first try to apply monotonicity estimates nonetheless. Suppose that the $k$-volume of $S$ is much smaller than $\vol(\Gamma\bs X)$. If the monotonicity estimate works at most points of $S$, we can run the same argument as for the proof of Theorem \ref{mt-stationary}. This would give a lower bound on the volume of $S$ contradicting our assumption. Hence, there are points where the monotonicity estimate fails. For those points we unfold the proof of the monotonicity estimates and construct a vector field whose flow squeezes the $k$-volume of $S$ for a small but definite amount of time. There are cases when the sketch above fails, for example when $S$ is a thin tube around a lower dimensional set, but in all these cases we are able to show that Federer-Fleming deformation (see Section \ref{sec-FF}) decreases the $k$-volume of $S$ by a definite amount. In this way we construct a sequence of deformations, either by flows or Federer-Fleming deformations, that will eventually bring the $k$-volume of $S$ near $0$. At that point one more application of the Federer-Fleming deformation collapses $S$ to a $(k-1)$-dimensional set. 

For the proof of branched cover stability (\textbf{Theorem \ref{mt-branched}}), we introduce a topological property called $W$-domination, for triples of closed subsets $N_1,N_2\subset W$ of a Riemannian manifold $M$, denoted $N_1\xRightarrow{W}N_2$. This property holds if and only if any cover of $M$ branched at $N_1$ can be isometrically embedded into a cover branched at $N_2$ after discarding the lift of $W$. It admits a simple algebraic description (see Definition \ref{def-Wdomination}). 
Now suppose $M$ is a compact octonionic hyperbolic manifold and $M_1\to M$ is a cover branched over a small closed rectifiable $14$-varifold $N\subset M$. We apply the sequence of deformations used in the proof of Theorem \ref{mt-collapse} to the branching locus $N$. The bulk of the proof of Theorem \ref{mt-branched} is showing that these deformations (either by flows or Federer-Fleming deformations)  produce a chain of subsets $N_i,W_i$ such that  $N=N_0\xRightarrow{W_0} N_1\xRightarrow{W_1}\ldots \xRightarrow{W_{i_0-1}}N_{i_0}$, where the final set $N_{i_0}$ is of codimension $3$. This means that any cover of $M$ branched at $N$ after discarding $\bigcup W_i$ embeds isometrically into a cover branched over $N_{i_0}$. But a cover branched over a codimension $3$ set is a genuine cover. The quantitative aspects of Theorem \ref{mt-branched} are handled by controlling the $(k+1)$-dimensional volumes of the sets $W_i$, which can be extracted from the proof of Theorem \ref{mt-collapse}. Several substantial technical difficulties have been omitted for clarity.


Finally we describe the proof of \textbf{Theorem \ref{mt-topwaist}}, which is the key ingredient in establishing property $\FA_r$ for higher rank lattices (Corollary \ref{corfanover8}.). Given a generic map $f: \Gamma \backslash X \rightarrow \mathbb{R}^d$, the preimages give a sweepout $\mathcal{F}$ of $\Gamma \backslash X$ defined by $\mathcal{F}: p \mapsto f^{-1}(p)$. We assume for contradiction that the fundamental group of each connected component of each preimage includes to an infinite index subgroup of  $\Gamma \backslash X$.  The general strategy is then to construct new sweepouts $\mathcal{F}_t$ so that the largest area of a cycle $\mathcal{F}_t(p)$ as $p$ ranges over $\mathbb{R}^d$ tends to zero as $t \to \infty$.  Then by making $t$ large enough, we can get a contradiction by using a theorem of Almgren (Theorem \ref{positivewaist}.)  

We build the $\mathcal{F}_t(p)$ by piecing together the vector fields constructed in (\ref{section:naturalvectorfields}) over a triangulation of $\mathbb{R}^d$. We consider the cover $\Lambda\bs X$ of $\Gamma\bs X$, where $\Lambda$ is the image of the fundamental group of $f^{-1}(p)$ in $\Gamma=\pi_1(\Gamma\bs X)$. The space $\Lambda\bs X$ admits a vector field inducing the \emph{natural flow} which contracts volume of low codimensional subsets, in particular $f^{-1}(p)$. These flows were first constructed and studied in \cite{cmw23}. It is here that we use the higher rank assumption\footnote{Or that $X=\mathbb H^2_\mathbb O$.}, because the range of codimensions $r(X)$ where the natural flow squeezes the volume depends on the growth gap phenomenon for infinite covolume subgroups. 
The next step is to flow lifts of the cycles of $\mathcal{F}$ to infinite volume covers of $\Gamma \backslash X$  along these vector fields, and then project them back down to $\Gamma \backslash X$. To ensure that the flows of the fibers glue together to the flow of a whole sweepout, we need to carefully interpolate the flows between fibers while ensuring that their volume squeezing properties are preserved. The implementation of this piecing-together procedure depends on the concept of a stratified fibration (Definition \ref{sf}), introduced by Gromov in \cite{gr09}.  An important fact is that convex combinations of the above vector fields still have the same area-contraction properties.  We comment that a proof of the result due to Alagalingam \cite{alagalingam2019algebraic} concerning tori mentioned in Section \ref{subsec-topwaistineq} can be given by following a similar strategy.  In place of the Connell-McReynolds-Wang natural flows, one instead uses vector fields whose flows contract $\mathbb{R}^n$ to an affine subspace, namely the gradient of the distance squared function to that affine subspace.  Their key property is that if the affine subspace is $i$-dimensional, then the flows contract $(i+1)$-dimensional volume.  

\begin{rem}
Connell-McReynolds-Wang \cite{cmw23} applied their natural flows to prove vanishing theorems for the cohomology of infinite volume higher rank locally symmetric spaces, following the work by Connell-Farb-McReynolds \cite{cfm16} and Carron-Pedon \cite{carron2004differential} in the rank one setting. This work in turn built on earlier work by Kapovich \cite{kapovich2009homological}, and also work by Besson-Courtois-Gallot \cite{bcgamalgamated} that initiated this approach to proving vanishing theorems for cohomology by constructing volume contracting natural maps. 
We comment that for the $\SL_n(\mathbb R)$ case, a lower bound for the range of dimensions for which the cohomology vanishes that is linear in rank follows from the computation of $r(  \SL_n(\mathbb R) / \SO(n))$ in Section \ref{computingr(x)}. The authors were inspired to begin thinking about the questions that led to this paper by the paper by Connell-Farb-McReynolds \cite{cfm16}. 
\end{rem}

\subsection{Comments on higher rank}\label{sec-IntroHR}

This description of the difficulties in higher rank  will probably make sense only after reading the proof of Theorem \ref{mt-stationary} in Section \ref{sec-stationaryproof}. To describe the nature of the difficulties in a little more detail, fix a Weyl chamber and denote by  $A^+$ its image under the exponential map.  Let $G=\SL_n(\mathbb R)$, and take a Cartan decomposition $KA^+K$ of $G$, for $K=\SO(n)$. Then given a minimal submanifold $S$ in $\Gamma \bs X=  \Gamma\bs G/K$, we expect similar arguments to the rank one case to work in establishing a monotonicity formula unless  $g^{-1}S$ ``sticks" to the wall in $
X= G /K =  KA^+K /K,$ for most $\Gamma gK\in S$.
By sticking to the wall we mean that the $A^+$ coordinate is near the walls of the Weyl chamber (the boundary of $A^+$.)  
We weren't able to find a test function playing the role of $\mathfrak f$ in our proof, whose Hessian would have desired properties near the walls of the Weyl chamber. In the ``sticking to the walls" scenario, a non-trivial contributions to the integrals from section \ref{sec-stationaryproof} would come from the part of the space where we have little control over the Hessian of $\mathfrak f$.  
We believe we will be able to rule out this improbable scenario by a mass transport argument (see e.g. \cite{aldous2007processes}), but it remains a significant obstacle compared to the rank one case. We expect that after dealing with this difficulty we could get analogues of Theorems \ref{mt-stationary},\ref{mt-collapse} and \ref{mt-branched} in the same range of codimensions as in Corollary \ref{corfanover8}. The second natural barrier for these methods is the minimal codimension where totally geodesic subspaces start to appear. For $\SL_n(\mathbb R)$ it would be codimension $n$. 


\subsection{Outline of the paper} 

In Section \ref{sec:background}, we give background from the theory of minimal submanifolds and  representation theory that will be useful in the paper.  In Section \ref{sec-Fed}, we describe the Federer-Fleming deformation.  In Section \ref{sec-Hess}, we perform the necessary computations of Hessians of distance and Busemann functions in symmetric spaces.  We also describe how one obtains the pluri-superharmonic functions for infinite covolume discrete subgroups, whose gradients give the Connell-McReynolds-Wang natural flows.  In Section \ref{sec-mono} we prove our monotonicity formulas for minimal submanifolds of low codimension in rank one symmetric spaces.  In Section \ref{sec-VarifoldsMonoMC} we give the proof of Theorem \ref{mt-stationary} by playing the decay of matrix coefficients against the monotonicity formula for stationary varifolds in octonionic locally symmetric spaces.  In Section \ref{sec-collapsing}, we prove  Theorem \ref{mt-collapse} by strengthening the methods of Section \ref{sec-stationaryproof} so that they also apply to varifolds that are not necessarily stationary.  In Section \ref{sec-stratified} we explain the concept of a stratified fibration, and show that generic piecewise linear maps are stratified fibrations.  In Section \ref{sec-topwaist}, we prove Theorem \ref{largefiber} on topological waists of higher rank locally symmetric spaces.  In Section \ref{section:thm2} we use the results of the previous section to establish Property $\FA_r$ for lattices in octonionic hyperbolic and certain higher rank symmetric spaces.  In Section \ref{sec-branched}, we apply Theorem \ref{mt-collapse} to prove non-abelian expansion and branched cover stability for uniform lattices in octonionic hyperbolic space.   In Appendix \ref{sec-sweepouts} we collect some facts about sweepouts that are used in the main body of the paper.  

 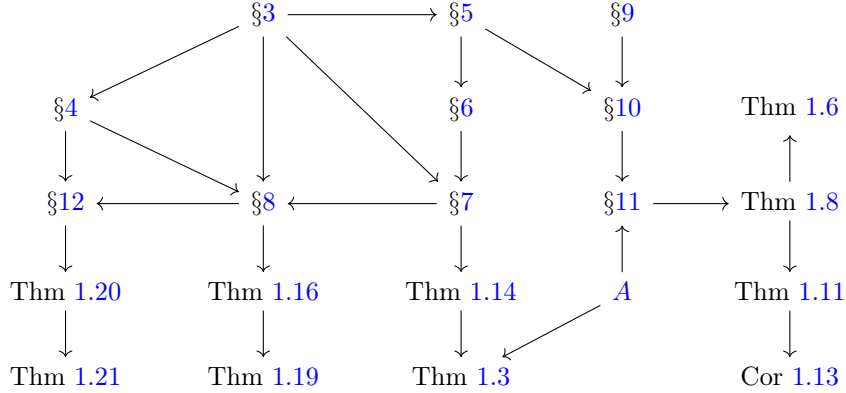
\begin{figure}[h!]
\[\begin{tikzcd}
	& \S\ref{sec:background} & \S\ref{sec-Hess} & \S\ref{sec-stratified} \\
	\S\ref{sec-Fed} & & \S\ref{sec-mono} & \S\ref{sec-topwaist} & \begin{array}{c} \text{Thm \ref{intro:critical exponent}}\\ \end{array} \\
	\S\ref{sec-branched} & \S\ref{sec-collapsing} & \S\ref{sec-VarifoldsMonoMC} & \S\ref{section:thm2} & {\text{Thm \ref{mt-topwaist}}} \\
	{\text{Thm \ref{mt-branched}}} & {\text{Thm \ref{mt-collapse}}} & {\text{Thm \ref{mt-stationary}}} & {\ref{sec-sweepouts}} & {\text{Thm \ref{mt-finiteorbit}}} \\
	{\text{Thm \ref{mt-nacheeger}}} & {\text{Thm \ref{mt-sysfreedom}}} & {\text{Thm \ref{mt-waist}}} & &{\text{Cor \ref{corfanover8}}}
	\arrow[from=1-2, to=3-3]
        \arrow[from=1-2, to=1-3]
	\arrow[from=1-2, to=2-1]
	\arrow[from=1-2, to=3-2]
        \arrow[from=1-3, to=2-3]
	\arrow[from=1-3, to=2-4]
	\arrow[from=1-4, to=2-4]
	\arrow[from=2-1, to=3-1]
	\arrow[from=2-1, to=3-2]
	\arrow[from=2-3, to=3-3]
	\arrow[from=2-4, to=3-4]
	\arrow[from=3-1, to=4-1]
	\arrow[from=3-2, to=3-1]
	\arrow[from=3-2, to=4-2]
	\arrow[from=3-3, to=3-2]
	\arrow[from=3-3, to=4-3]
	\arrow[from=3-4, to=3-5]
	\arrow[from=3-5, to=2-5]
	\arrow[from=3-5, to=4-5]
	\arrow[from=4-1, to=5-1]
	\arrow[from=4-2, to=5-2]
	\arrow[from=4-3, to=5-3]
	\arrow[from=4-4, to=3-4]
	\arrow[from=4-4, to=5-3]
        \arrow[from=4-5, to=5-5]
\end{tikzcd}\]
\caption{Pathways to the main results of the paper.}
\end{figure}
Reader interested only in one result, only needs to read its antecedents in the graph above. For example, for Theorem \ref{mt-stationary} it is enough to read \S\ref{sec:background}, \S\ref{sec-Hess}, \S\ref{sec-mono} and \S\ref{sec-VarifoldsMonoMC} while for Theorem \ref{mt-finiteorbit} it sufficient to read \S\ref{sec:background}, \S\ref{sec-Hess}, \S\ref{sec-stratified}, \S\ref{sec-topwaist}, \S\ref{section:thm2} and the Appendix \ref{sec-sweepouts}.

\section{Acknowledgements}  We thank Shmuel Weinberger for helpful discussions. We are also grateful to Benson Farb for a helpful discussion regarding property $\FA_r$.  We thank Larry Guth for bringing to our attention the work by Alagalingam. We also thank Yangyang Li. MF was supported by the Dioscuri programme initiated by the Max Planck Society, jointly managed with the National Science
Centre in Poland, and mutually funded by the Polish Ministry of Education
and Science and the German Federal Ministry of Education and Research. BL was supported by the NSF under grant DMS-2202830.

\section{Background} \label{sec:background}
\subsection{Varifolds and the first variation formula}\label{sec-varifolds}
In this section we collect basic results and definitions concerning rectifiable varifolds. We fix an ambient Riemannian manifold $M$ of dimension $n$. Write $g$ for the metric tensor and $\nabla$ for the Levi-Civita connection. Let $\mathcal H^k$ denote the $k$-dimensional Hausdorff measure on $M$ \cite[Chap.1 \S 2]{simongmt}. Recalling that Lipschitz maps are differentiable a.e., for any $k$-dimensional Lipschitz embedded submanifold $S\subset M$, the measure $\mathcal H^k|_S$ coincides with the induced Riemannian volume on $S$.
\begin{defn}[{\cite[Chap.3 \S 1]{simongmt}}]
A subset of $M$ is countably $k$-rectifiable if it is a countable union of Lipschitz images of subsets of $\mathbb R^k$ and a set of $\mathcal H^k$-measure $0$.
\end{defn}
This definition actually implies a stronger characterization \cite[Chap 3. Lemma 1.2]{simongmt}. A countably $k$-rectifiable subset $A\subset M$ admits a decomposition
\begin{equation}\label{eq-recfiabledecomp}
A=A_0\cup \bigcup_{i=1}^\infty N_i,
\end{equation}
where each $N_i$ is a $C^1$-embedded (but not necessarily closed) $k$-submanifold and $\mathcal H^k(A_0)=0.$
\begin{defn}
Let $A\subset M$ be an $\mathcal H^k$-measurable, countably $k$-rectifiable set such that $\mathcal H^k(A\cap K)$ is finite for all compact sets $K\subset M$. Then, for $x\in A$, a $k$-dimensional subspace $W\subset T_xM$ is said to be tangent to $A$ if for any compactly supported continuous function $f\colon T_x M \to \mathbb R_{\geq 0}$ we have 
$$\lim_{\lambda\to 0}\lambda^{-k}\int_A f(\lambda^{-1}\log(x))d\mathcal H^k(x)=\int_W f(y)d\mathcal H^k(y).$$
\end{defn}
We write $T_x A$ for the tangent space of $A$ at $x$. In intuitive terms, zooming in at $x$ makes the $k$-dimensional Hausdorff measure on $A$ look more and more like the $k$-dimensional Lebesgue measure on $\exp(W)$. By \cite[Chap. 3 Thm 1.6]{simongmt}, if $A$ is countably $k$-rectifiable and has locally finite $\mathcal H^k$-measure, then $\mathcal H^k$-almost every point of $A$ admits a unique tangent plane. These point are called the \emph{regular points} of $A$. For the parts of $A$ which are $C^1$-submanifolds (see \ref{eq-recfiabledecomp}) this definition of tangent space agrees with the ordinary definition of tangent space. In particular, the tangent space varies continuously after excluding a set of arbitrarily small $\mathcal H^k$-measure.

\begin{defn}[{\cite[Chap. 4]{simongmt}}]
Let $A\subset M$ be an $\mathcal H^k$-measurable, countably $k$-rectifiable set. Let $\theta$ be a non-negative $\mathcal H^k$-measurable function on $A$, such that $\int_{M\cap K}\theta d\mathcal H^k<\infty$ for any compact set $K.$ The rectifiable $k$-varifold $\underline{v}(A,\theta)$ attached to the pair $(A,\theta)$ is the equivalence class of pairs $(A',\theta')$ such that $A\setminus A', A'\setminus A$ are $\mathcal H^k$-measure $0$ and $\theta=\theta'$ $
\mathcal H^k$-almost everywhere. A varifold is integral if $\theta$ is integer valued $\mathcal H^k$-almost everywhere. 
\end{defn}

In this paper we work exclusively with integral rectifiable $k$-varifolds. We suppress the function $\theta$ from the notation and skip the adjective ``integral" from now on. Typically a rectifiable varifold will be denoted by the letter $S$. When integrating over a varifold $S$ represented by a pair $(A,\theta)$ we put 
$$\int_S f(x)d \mathcal H^k(x):= \sum_{i=1}^{\infty } \int_{N_i} \theta(x)f(x)\mathcal H^k(x) ,$$ 
for any $f$ for which the right hand side converges, where $A= A_0 \cup (\cup_{i=1}^\infty N_i)$ as in (\ref{eq-recfiabledecomp}.) 
Given a $C^1$-vector field $X$ on $M$ we define the associated flow $\Phi_t\colon M\to M, t\in \mathbb R$ by the condition $\frac{\partial}{\partial t}\Phi_t(x)=X(\Phi_y(x))$ and $\Phi_0(x)=x$ for all $x\in M$, see \cite[\S2.5]{perko2013differential} or \cite[Thm 17.8]{lee2012smooth} for the existence and uniqueness theorems.

Since a rectifiable $k$-varifold $S$ is $\mathcal H^k$-almost everywhere a finite union of $C^1$-submanifolds \cite[Chap 3. Lemma 1.2]{simongmt}, we use the standard differential geometry conventions with the understanding that they are well defined $\mathcal H^k$-almost everywhere on $S$. 

\begin{defn} \cite[\S4 Def 2.4]{simongmt} A rectifiable $k$-varifold $S\subset M$
is \emph{stationary} if for every compactly supported $C^1$-vector field $X$ on $M$ with the associated flow $\Phi_t$, we have $\frac{\partial}{\partial t}\mathcal H^k(\Phi_t(S))=0.$
\end{defn}

We now define the divergence of a vector field relative to $S$. Choose $C^1$-vector fields $e_1, \ldots, e_k$ on $M$ that restrict to an orthonormal frame in $TS$ at every regular point of $S$. For any $C^1$ vector field $X$ on $M$, put
\begin{equation}\label{eq-divergence}
{\rm div}_S X=\sum_{i=1}^k \langle \nabla_{e_i} X,e_i\rangle. 
\end{equation}
Note that when $X=\grad f, f\in C^2(M)$, then 
\begin{equation}\label{eq-Hess}
{\rm div}_S \grad f=\sum_{i=1}^k \langle \nabla_{e_i} \grad f,e_i\rangle = \sum_{i=1}^k \Hess f(e_i,e_i)=\tr \Hess f|_{TS}. 
\end{equation}
By \cite[Chap 2 \S5, Chap 4 \S2]{simongmt},
\begin{equation}\label{eq-FV}
\frac{\partial}{\partial t}\mathcal H^k(\Phi_t(S))=\int_S {\rm div}_S X d\mathcal H^k 
\end{equation}
When $S$ is stationary, the first variation of the Hausdorff measure is zero, so
\begin{equation}\label{eq-statFV}
\int_S {\rm div}_S X d\mathcal H^k= 0, 
\end{equation}
for any compactly supported $C^1$-vector field $X$.
In the particular case $X=\grad f, f\in C^2(M)$ we have ${\rm div_S} \grad f(x)=\tr \Hess f|_{T_xS}$ for $\mathcal H^k$-almost every $x\in S$. We can therefore write
\begin{equation}\label{eq-statHessianBd}
\frac{\partial}{\partial t} \mathcal H^k(\Phi_t(S))=\int_S {\rm div}_S \grad f d\mathcal H^k=\int_S \tr \Hess f|_{TS}\, d\mathcal H^k=0.
\end{equation}
We finish this section by giving an abstract form of the monotonicity estimate. 
Let us fix a smooth non-increasing function $\chi\colon \mathbb R\to \mathbb R_{\geq 0}$ such that $\chi(t)=1$ for $t<0$, $\chi(t)=0$ for $t>1$ and $\chi'(t)\geq -2$ for all $t\in\mathbb R$. We put $\chi_r(t):=\chi(t-r).$ The function $\chi_r$ will serve as smooth cut-off throughout the paper. The exact choice of $\chi$ is not important, as long as it satisfies the properties listed above. 
\begin{lem}\label{lem-abstmonotonicity}
Let $f\colon M\to\mathbb R_{\geq 0}$ be a proper $C^2$ function with $\|\grad f\|\leq 1,$ with the property that $\tr \Hess f(x)|_W\geq \kappa$ for every $x\in M$ and every $k$-dimensional subspace $W\subset T_x M$. For any stationary $k$-varifold $S\subset M$, we have 
$$ \int_S \chi_r(f)d\mathcal H^k \geq e^{\kappa(r-s)}\int_S \chi_s(f)d\mathcal H^k.$$
\end{lem}
\begin{proof}
Note that at any regular point of $S$ we have 
\begin{align*} {\rm div}_S (\chi_r(f)\grad f)=& \chi_r(f){\rm div}_S\grad f +\chi'_r(f)\|\grad_S f\|^2\\
=& \chi_r(f) \tr \Hess f|_{TS} + \chi_r'(f)\|\grad_S f\|^2.
\end{align*}
We apply the first variation formula to the vector field $\chi_r(f)\grad f.$
\begin{align*}
\int \chi_r(f(x))\tr \Hess f|_{T_xS}d\mathcal H^k(x)=& - \int_S \chi_r'(f(x))\|\grad_S f\|^2 d\mathcal H^k(x)\\
\kappa \int_S \chi_r(f(x))d\mathcal H^k(x)\leq& - \int_S \chi_r'(f(x))d\mathcal H^k(x)\\
=& \frac{\partial}{\partial r}\int_S \chi_r(f(x))d\mathcal H^k(x).
\end{align*}
Note that for the second inequality we used $\chi_r'\leq 0$ and $\|\grad_S f\|\leq \|\grad f\|\leq 1.$ The conclusion of the lemma then follows by applying Gr\"onewall's inequality \cite{gronwall1919note} to the function $r\mapsto \int_S \chi_r(f(x))d\mathcal H^k(x).$
\end{proof}
We will need one more fact about stationary varifolds in symmetric spaces.
\begin{lem}\label{lem-smallballmono}
Let $X$ be a symmetric space. Then for any $r$ there is $\delta(r)>0$ depending on $X$ so that any stationary varifold $S\subset X$ and any regular point $x\in S$ satisfies
$$\mathcal H^k(S\cap B(x,r))\geq \delta(r).$$

\end{lem}
\begin{proof}

By \cite[Section 6]{montezumanotes}, there is $r_0, C$ depending only on $X$ so that any stationary $k$-varifold $S\subset X$, $x \in S$ regular and $r<r_0$ we have
\[
\frac{d}{dr}\left( \log \frac{e^{Cr^2}\mathcal H^k(B(x,r) \cap S)}{r^k} \right) \geq 0.  
\]
Since the limit $\lim_{r \to 0} \mathcal H^k(B(x,r)\cap S)/r^k$ is a positive integer multiple of the volume of the k-dimensional unit ball, the conclusion follows.  

\end{proof}
\subsection{Symmetric spaces and spherical functions}\label{sec-symspaces}
We follow the standard notations of \cite{knapp1996lie,gv88}, with the exception of Iwasawa decomposition. We will be using $NAK$ instead of $KAN$, which is the order of factors in \cite{gv88} and most sources \footnote{The reason is that we want the central Iwasawa coordinate $H\colon G\to \mathfrak a$ to be a right $K$-invariant function.}. Throughout the paper, $G$ will be reserved for a semisimple real Lie group. We let $K$ be a maximal compact subgroup of $G$ fixed by a Cartan involution $\Theta$. The symmetric space of $G$ is the quotient $X=G/K$. Write $\mathfrak g, \mathfrak k$ for the Lie algebras of $G,K$ respectively and let $\mathfrak s=\{Y\in \mathfrak g\mid \Theta Y=-Y\}$. Let $A$ be a maximal split torus of $G$ stabilized by $\Theta$ and let $\mathfrak a$ be the Lie algebra of $A$. We necessarily have $\mathfrak a\subseteq \mathfrak s$. The tangent space $T_KX$ (at the identity coset of $K$) is $\mathfrak g/\mathfrak k$. The latter is identified with $\mathfrak s$ via the map $\mathfrak g/\mathfrak k\ni Y+\mathfrak k\mapsto \frac{1}{2}(Y+\Theta Y)\in \mathfrak s$. The Killing form on $\mathfrak g$ is given by 
$$B(Y,Z)=\tr \ad Y\ad Z.$$
It descends to a positive definite symmetric bilinear form 
$$\langle Y,Z\rangle_{\rm Kill}:= -B(Y,Z) \hspace{3mm} Y,Z\in \mathfrak s.$$
This extends to a unique right-$G$ invariant Riemannian metric on  $X$, denoted $g_{\rm Kill}$.  We will write $\langle \cdot,\cdot\rangle, g$ for unspecified (but fixed) real multiples of $\langle \cdot,\cdot\rangle_{\rm Kill}, g_{\rm Kill}.$ The proofs and statements in the paper are written in a way that does not depend on the chosen normalization. For example, in the case of rank one symmetric spaces, we can choose $g$ so that the maximal negative curvature is $-1$ as opposed to keeping the metric induced by the Killing form.

We will make no distinction between functions on $X$ and right-$K$ invariant functions on $G$. In this way, given $f\colon X\to \mathbb R$ we write $f(g)$ for $f(gK), gK\in X$. 

Let $P$ be a minimal parabolic subgroup of $G$ containing $A$. The group $P$ decomposes as $P=MAN$ with $M=P\cap K=Z(A)\cap K$ and $N$ being the unipotent radical of $P.$ This decomposition is often called the Langlands decomposition. Let $\mathfrak p, \mathfrak m, \mathfrak n$ be the Lie algebras of $P,M,N$ respectively.
The non-negative and positive Weyl chambers are defined in the standard way.
\begin{align*}
A^+=&\{a\in A\mid a^{-t} n a^t \to 1 \text{ as } t\to\infty \text{ for all }n\in N\},\\
\mathfrak a^+=&\log A^+\subset \mathfrak a \quad\text{ (the non-negative Weyl chamber)},\\
\mathfrak a^{++}=&\,{\rm int}\, \mathfrak a^+ \quad\text{ (the positive Weyl chamber)}.
\end{align*}
Let $\mathfrak a^*:=\Hom(\mathfrak a,\mathbb R)$ and $\mathfrak a^*_\mathbb C:=\Hom(\mathfrak a,\mathbb C).$ Let $\lambda\in \mathfrak a^*$. We put 
$$\mathfrak g_\lambda:=\{Y\in \mathfrak g\mid [H,Y]=\lambda(H)Y \text{ for } H\in \mathfrak a\}.$$
If $\mathfrak g_\lambda\neq 0$, then we call $\lambda$ a \emph{root} and put $m_\lambda:=\dim \mathfrak g_\lambda.$ The number $m_\lambda$ is called the \emph{multiplicity} of $\lambda$ and $\mathfrak g_\lambda$ is the \emph{root subspace} of $\lambda.$ We typically reserve letters $\alpha,\beta$ for the roots and $\lambda,\xi$ for general elements of $\mathfrak{a}^*$.
Let $\Phi$ be the set of roots and let $\Phi^+:=\{\alpha\in \Phi\mid \mathfrak g_\alpha\subseteq \mathfrak n\}$ be the set of positive roots. We write $\tilde \Phi,\tilde \Phi^+$ for the multisets of roots and positive roots, where $\alpha$ has multiplicity $m_\alpha.$
Let $\rho$ be the half-sum of positive roots $$2\rho:=\sum_{\alpha\in \Phi^+}m_\alpha \alpha= \sum_{\alpha\in\tilde\Phi^+}\alpha \in \mathfrak a^*.$$
The \emph{simple roots} are those positive roots which cannot be expressed as positive integer combinations of other positive roots. The set of simple roots will be denoted $\Delta$, their number is equal to the rank of $G$ and they form a basis of $\mathfrak a^*$. Moreover, any positive root can be expressed as a positive integer combination of simple roots. We say that a root $\alpha$ is \emph{reduced} if $\alpha/2\not\in\Phi.$
The \emph{Weyl group} is the quotient $N(A)/A$ of the normalizer in $G$ of $A$ by $A$. It acts on $A,\mathfrak a$ and $\mathfrak a^*$. The non-negative Weyl chamber $\mathfrak a^+$ is a fundamental domain for the action on $\mathfrak a$. The action is free on the union $\bigcup_{w\in W}w\mathfrak a^{++}.$
\emph{Iwasawa decomposition}\footnote{In other sources one might find a different order of factors. We opt to have $K$ on the right so that the Iwasawa coordinate $H$ descends to a function on $X=G/K$.} is the identity $G=NAK$. This decomposition is unique and the map sending $(n,a,k)\to nak$ is a diffeomorphism \cite[p.63]{gv88}. We define the function $H\colon G\to \mathfrak a$ by 
$$g=ne^{H(g)}k, k\in K,n\in N.$$ We can integrate in Iwasawa coordinates using \cite[Prop. 2.4.2]{gv88}
$$\int_G f(g) dg=\int_N\int_{\mathfrak a}\int_K f(ne^Hk)dk dH dN,$$ for suitable choices of Haar measures on $K,\mathfrak a, N.$
 \emph{Cartan decomposition} is the identity $G=KA^+K$. We define the function $a\colon  G\to \mathfrak a^+$ by putting $$g=k_1e^{a(g)}k_2, k_1,k_2\in K.$$ The map $(k_1,a,k_2)\to k_1ak_2$ is analytic diffeomorphism on $\mathfrak a^{++}$ and the image is an open dense subset of $G$. To integrate in Cartan coordinates we use the formula \cite[Prop. 2.4.6]{gv88}
$$ \int_G f(g)dg =\int_K\int_K\int_{\mathfrak a^+}f(k_1e^Hk_2) J(H)dHdk_1dk_2,$$ where $J(H):=\prod_{\alpha\in\Phi^+}\sinh(\alpha(H)))^{m_\alpha}$ and $dk_1,dk_2,dH$ are suitable Haar measures on $K,K,\mathfrak a.$

A function $\phi\colon G\to\mathbb R$ is a \emph{spherical function} if $\phi(k_1gk_2)=\phi(g)$ for all $g\in G, k_1,k_2\in K.$ Let $V$ be a Hilbert space and let $\pi\colon G\to \mathscr U(V)$ be a unitary representation of $G$. Let 
$$V^K:=\{v\in V\mid \pi(k)v=v, k\in K\},$$ be the subspace of $K$ fixed vectors. For $v,w\in V^K$, the matrix coefficient
$g\mapsto \langle \pi(g)v,w\rangle$ is a spherical function. If $(\pi, V)$ is irreducible and $V^K\neq 0$, we know that $\dim V^K=1$ \cite[Prop. 1.5.8]{gv88}.
The elementary spherical functions are the matrix coefficients $g\mapsto \langle \pi(g)v,v\rangle$, where $(\pi,V)$ is irreducible, $v\in V^K$ and $\|v\|=1.$
Thanks to Harish-Chandra, we have an explicit formulae for the elementary spherical functions.
\begin{thm}\label{thm-ElmSph}\cite[Prop 3.1.4 and Thm. 3.2.3]{gv88}
Let $(\pi, V)$ be an irreducible unitary representation with $V^K\neq 0$. Let $v\in V^K$, $\|v\|=1$. There exists a $\lambda\in \mathfrak a^*_\mathbb C$ such that 
$$\langle \pi(g)v,v\rangle =\int_K e^{(\lambda-\rho)(H(kg))}dk=:\varphi_\lambda(g).$$
The character $\lambda$ is determined uniquely up to the $W$-action and called the \emph{infinitesimal character} of $\pi.$ Characters $w\lambda, w\in W$ give rise to the same elementary spherical function.
\end{thm}
We will make use of the following identities and estimates for elementary spherical functions \cite[Thm 4.6.4]{gv88}. Let $d$ be the number of reduced positive roots counted without multiplicity. It is always bounded by $|\Phi^+|\leq |\tilde \Phi^+|=\dim X-\dim \mathfrak a.$
\begin{align} \varphi_0(g)\leq& e^{-\rho(a(g))}(1+\|a(g)\|)^{d}, \label{eq-basicMCbounds}\\
\varphi_{w\rho}(g)=&1 \text{ for all }w\in W, g\in G. \nonumber
\end{align}
The second formula can be recovered from the integral representation in Theorem \ref{thm-ElmSph} for $w=1$ and the fact that $\varphi_{w\lambda}=\varphi_{\lambda}$ for all $\lambda\in \mathfrak a^*_\mathbb C, w\in W$ \cite[Prop. 3.2.2]{gv88}. 
\begin{lem}\label{lem-MCconvexity}
The map $\lambda\mapsto \log \varphi_\lambda(g)$ is convex on $\mathfrak a^*.$
\end{lem}
\begin{proof}
To prove our claim we just need to show that the Hessian of $\log \varphi_\lambda(g)$ is non-negative definite.  Let $\xi\in \mathfrak a^*\cong T_\lambda \mathfrak a^*$. By the Cauchy-Schwartz inequality,
\begin{align*}\frac{\partial^2}{\partial \xi\partial \xi}\log \varphi_\lambda(g)=&\varphi_\lambda(g)^{-2}\left(\int_K e^{(\lambda-\rho)(H(gk))}dk\int_K \xi(H(gk))^2 e^{(\lambda-\rho)(H(gk))}dk\right.\\ & \left.-\int_K \xi(H(gk))e^{(\lambda-\rho)(H(gk))}dk\int_K \xi(H(gk))e^{(\lambda-\rho)(H(gk))}dk\right)\geq 0.
\end{align*} 
\end{proof}
\subsection{Decay of matrix coefficients in rank one groups}\label{sec-decay}
Throughout this section we let $G$ be a rank one simple real Lie group. Let $\alpha$ be the simple root. Then the set of positive roots is either $\{\alpha\}$ or $\{\alpha, 2\alpha\}.$ The multiplicities are listed in the table below, covering all isogeny classes of real rank one groups
\begin{table}[ht]
\centering
\begin{tabular}{ c c c c c}
$G$ & $X$ & $\dim X$ & $m_\alpha$ & $m_{2\alpha}$ \\ \hline
$\SO(n,1)$ & $\mathbb H^n_\mathbb R$ & $n$ & $n-1$ & $0$\\
$\SU(n,1)$ & $\mathbb H^n_\mathbb C$ & $2n$ & $2n-2$ & $1$\\
$\Sp(n,1)$ & $\mathbb H^n_\mathbb H$ & $4n$ & $4n-4$ & $3$\\
$F_4^{(-20)}$ & $\mathbb H^2_\mathbb O$ & $16$ & $8$ & $7$\\
\end{tabular}
\caption{Multiplicities of roots in rank one groups.}
\label{tab-multiplicities}
\end{table}
\begin{thm}[{\cite{kostant1969existence}}] Let $G$ be a rank one simple real Lie group and let $(\pi,V)$ be a non-trivial irreducible unitary representation with $V^K\neq 0.$ Then, the infinitesimal character $\lambda$ satisfies 
$$\left| \Re \frac{\langle \lambda,\alpha\rangle}{\langle \alpha, \alpha\rangle} \right|\leq \begin{cases} \frac{m_\alpha}{2}+1 & \text{if} \hspace{1mm} m_{2\alpha}\neq 0\\ \frac{m_\alpha}{2} & \text{ if }m_\alpha=0 \end{cases}.$$
Furthermore, if $\lambda\not\in i\mathfrak a^*$ then $\lambda\in \mathfrak a^*$ (it's either real or purely imaginary). 
\end{thm}
\begin{lem}\label{lem-rankonedecay} Let $G$ be either $\Sp(n,1)$ or $F_4^{(-20)}$. Let $\varphi_\lambda$ be an elementary spherical function corresponding to a non-trivial unitary representation. Then
$$| \varphi_\lambda(k_1e^Hk_2)|\leq e^{(1-m_{2\alpha})\alpha(H)}(1+\|H\|), \text{ for } H\in \mathfrak a^+, k_1,k_2\in K.$$
\end{lem}
\begin{proof}
The function is spherical, so we might as well assume $k_1=k_2=1$. If $\lambda\in i\mathfrak a^*$ then by (\ref{eq-basicMCbounds})
$$|\varphi_\lambda(e^H)|\leq |\varphi_0(e^H)|\leq e^{-\rho(H)}(1+\|H\|)\leq e^{(1-m_{2\alpha})\alpha(H)}(1+\|H\|).$$

Otherwise, by Kostant's theorem, $\lambda= t_1\alpha,$ with $|t_1|\leq \frac{m_\alpha}{2}+1.$ Applying the non-trivial Weyl group element to $\lambda$ (which in rank one case is the multiplication by $-1$) does not change the spherical function, so we can assume $t_1\geq 0$. We have $\lambda=\frac{(\rho-\lambda)(H)}{\rho(H)}0+\frac{\lambda(H)}{\rho(H)}\rho$, for any non-zero $H\in \mathfrak a$. Lemma \ref{lem-MCconvexity} yields
\begin{align*}
\log|\varphi_\lambda(e^H)|\leq& \frac{\lambda(H)}{\rho(H)}\log |\varphi_\rho(e^H)|+\frac{(\rho-\lambda)(H)}{\rho(H)}\log |\varphi_0(e^H)|\\
\leq& \frac{(\rho-\lambda)(H)}{\rho(H)}\log(e^{-\rho(H)}(1+\|H\|))\\
\leq& (\lambda-\rho)(H)+ \log (1+\|H\|).
\end{align*}
We have $\rho=(\frac{m_\alpha}{2}+m_{2\alpha})\alpha,$ so $(\lambda-\rho)(H)\leq (1-m_{2\alpha})\alpha(H).$
\end{proof} 
\section{Federer-Fleming deformation}\label{sec-Fed}
\subsection{Federer-Fleming deformation relative to a uniform triangulation}\label{sec-FF}
In this section we define the Federer-Fleming deformation. It is a homotopy taking a rectifiable $k$-varifold $S$ in a triangulated  manifold $M$ to a subset of the $k$-skeleton of the triangulation. These types of deformations were originally considered by Federer-Fleming in \cite{federer1960normal}. Federer-Fleming deformation serves as one of the key ingredients of the proof of Theorem \ref{mt-collapse}. It is obtained by successively collapsing $S$ onto lower dimensional skeleta, each time taking care to not increase the $k$-volume by more than a fixed multiplicative constant. One technical difficulty in applying this deformation is the dependence of the implicit constants on the triangulation. We wish for our arguments to work uniformly across the different families of locally symmetric spaces so we must control how the quality of the deformation depends on the triangulation and the geometry of the manifold. 

We recall that a triangulation of $M$ is a simplicial complex $\mathcal T$ together with a homeomorphism $\mathcal T\cong M.$ We shall identify the cells $\sigma\in \mathcal T$ with their image in $M.$ An $m$-simplex will refer to a simplex of dimension $m$ which will be closed unless otherwise mentioned.
We write $M^{[m]}, (M^{[\bullet]})$ for the set of (oriented) $m$-simplices (resp. all simplices) and $M^{(m)}$ for the $m$-skeleton on $M$. An \emph{affine simplex} of dimension $m$ is a convex hull of $m+1$ points in $\mathbb R^m$ with the standard Euclidean metric. An affine $m$-simplex is \emph{non-degenerate} if it is not contained in an $(m-1)$-dimensional hyperplane. 

\begin{defn}\label{def-UnifTrian} Let $M$ be a complete Riemannian manifold of dimension $n$. A triangulation $\mathcal T$ of $M$ is $(r,\delta)$-uniform if the following conditions hold. 
\begin{enumerate}
\item every cell $\sigma\in M^{[\bullet]}$ satisfies $\diam(\sigma)\leq r,$
\item every cell $\sigma\in M^{[\bullet]}$ is $(1+\delta)$-isometric\footnote{A map between metric spaces is a $(1+\delta)-$isometry if it is a bi-Lipschitz bijection with Lipschitz constants for the map and its inverse bounded by $(1+\delta)$. } to an affine simplex in $\mathbb R^{\dim \sigma},$\
\item any cell $\sigma\in M^{[\bullet]}$ satisfies $\mathcal H^k(\sigma)\geq \delta r^{\dim\sigma}.$
\end{enumerate}
\end{defn}

The construction of the Federer-Fleming deformation is based on the following lemma (see \cite[Section 8]{guth2010metaphors}.)  

\begin{figure}[h!]
		\includegraphics[width=100mm]{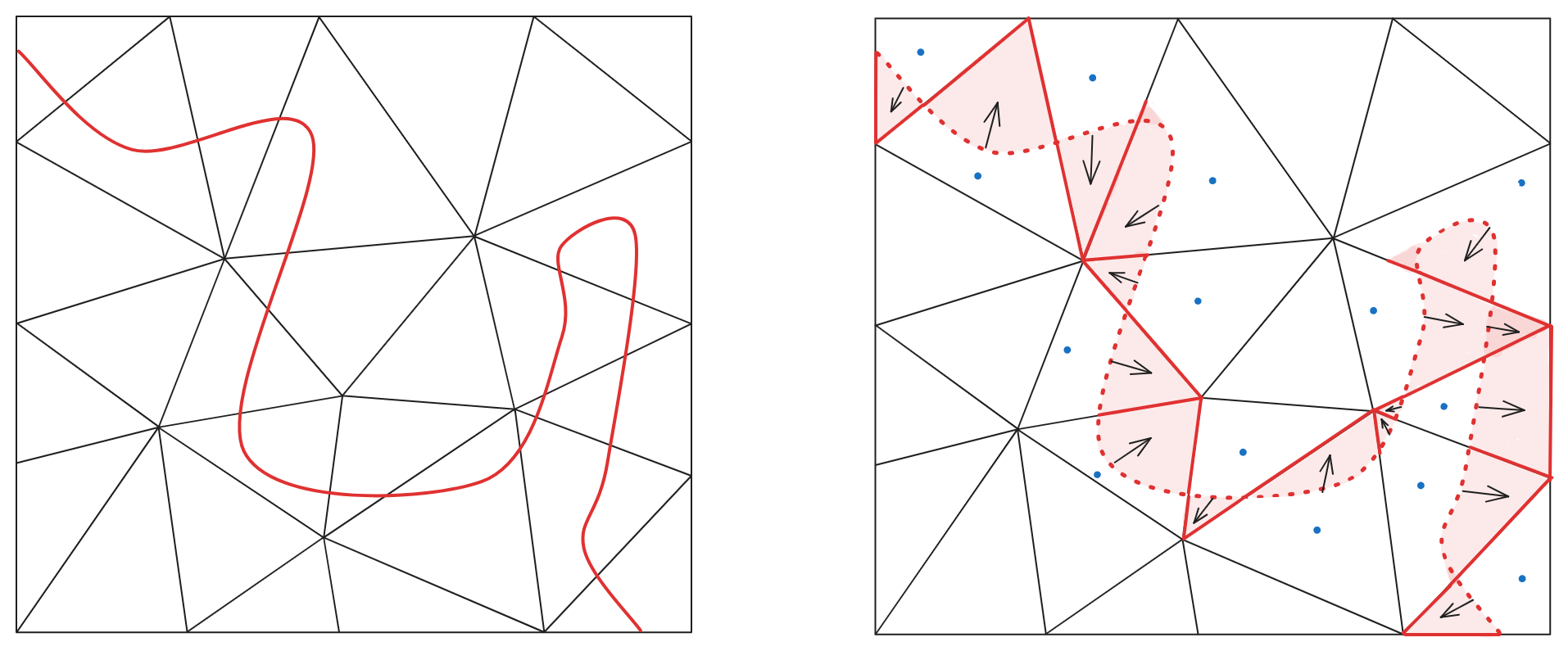}
		\caption{$\FF^{1}(S)$ when $k=1,\dim M=2.$ } \label{fig-FF} with $h^1(S,[0,1])$ highlighted.
	\end{figure}

\begin{lem}\label{lem-FFlocal} Let $\sigma\subset \mathbb R^m$ be an affine $m$-simplex satisfying $\diam \sigma\leq r$ and $\mathcal H^m(\sigma)\geq \delta r^m$. For a point $x\in {\rm int}(\sigma),$ let $p_x\colon\sigma\setminus \{x\}\to \partial \sigma$ be the map sending $y\in \sigma\setminus\{x\}$ to the unique point $p_x(y)\in \partial \sigma$ such that $y$ lies on the line segment connecting $x$ to $p_x(y).$ Let $h_x\colon (\sigma\setminus\{x\})\times [0,1]\to \sigma$ be the homotopy $h_x(y,t):=(1-t)y+t p_x(y)$. Let $S\subset \sigma$ be a rectifiable $k$-varifold. 
\begin{enumerate}
\item if $k<m$, there exists a point $x_0\in \sigma\setminus S$, such that $\mathcal H^k(p_{x_0}(S))\leq c \mathcal H^k(S)$ and $\mathcal H^{k+1}(h_{x_0}(S,[0,1]))\leq c\mathcal H^k(S)$ where $c=c(k,m,\delta)>0.$
\item if $k=m$ and $x_0\in {\rm int}\,\sigma\setminus S$, then $\mathcal H^{k+1}(h_{x_0}(S,[0,1]))=\mathcal H^k(p_{x_0}(S))=0$.
\end{enumerate}
\end{lem}
 Let $k<\dim M$ and let $S\subset M$ be a closed rectifiable $k$-varifold. 
 We now define the Federer-Fleming deformation of $S$ relative to the triangulated manifold $M$.  We will repeatedly use Lemma \ref{lem-FFlocal} to collapse $S$ onto lower skeleta of $M$, simplex by simplex, until we reach the dimension $k$ or $k-1$. 
As the first step, we fix a family of $(1+\delta)$-isometries $\iota_\sigma\colon \sigma \to \mathbb R^{\dim \sigma}, \sigma\in M^{[\bullet]}$ such that $\iota_\sigma(\sigma)$ is an affine simplex. We construct the deformation in steps. Let $\FF^n$ be the identity map. Suppose $k<m\leq n$, and that we have already defined $\FF^{m}$ with the following properties:
\begin{itemize}
\item $\FF^{m}$ is the identity on $S$ intersected with $M^{(m)},$
\item $\FF^{m}(S\cap \tau)$ is contained in the $m$-skeleton of the simplex $\tau$, when $\dim \tau>m$.
\end{itemize}
We now explain how to construct $\FF^{m-1}$. For each $m$-simplex $\sigma$ we choose a point $x_\sigma\in \iota_\sigma(\sigma)$ such that the conclusions of Lemma \ref{lem-FFlocal} (1) hold for the varifold $\iota_\sigma(\FF^{m}(S)\cap \sigma)$. Put 
$$\FF^{m-1}(s):=\iota_\sigma^{-1}(p_{x_\sigma}(\iota_\sigma(\FF^m(s)))),$$ for any $s\in (FF^{m})^{-1}(\sigma).$
The maps $\FF^{m}$ and $\FF^{m-1}$ are connected by the homotopy 
$$h^m(s,t)=\iota_\sigma^{-1}(h_{x_\sigma}(\iota_\sigma(\FF^m(s)),t)), s\in (\FF^m)^{-1}(\sigma), t\in[0,1].$$ Indeed, we have $h^m(s,0)=\FF^m(s), h^m(s,1)=\FF^{m-1}(s), s\in S.$
The map $\FF^{m-1}$ is the identity on the $(m-1)$-skeleton of $\mathcal T$, and it takes $S\cap \sigma$ to the $(m-1)$-skeleton of $\sigma$. Figure \ref{fig-FF} illustrates the process in case $\dim M=2.$
 
We proceed until we reach $m=k$. When $m=k$, for all $k$-cells $\sigma$ such $\FF^k(S)\cap \sigma\neq \sigma$ we choose a point $x_0\in {\rm int}\,\sigma \setminus S$ (note $\FF^k(S)$ is compact so we will be able to find $x_0$ in the interior of $\sigma$) and define $\FF^{k-1}(s)$ and $h^{k}(s)$ identically as before, using Lemma \ref{lem-FFlocal} (2). For $k$-cells $\sigma$ with $\FF^k(S)\cap\sigma=\sigma$ we put $\FF^{k-1}=\FF^{k}$ and $h^k(s,t)=\FF^{k}(s), t\in [0,1]$. Finally, we put $\FF:=\FF^{k-1}$ and let $h\colon S\times [0,1]\to M$ be the concatenation of $h^n,h^{n-1},\ldots,h^k.$ 

\begin{lem}\label{lem-FFproperties}
\begin{enumerate}
\item For any $\sigma\in M^{[\bullet]}, m=0,\ldots,n$ we have $\FF(S\cap \sigma)\subset \sigma$. 
\item For any $\sigma\in M^{[n]}$ we have $$\mathcal H^k(\FF(S\cap \sigma))\leq c \mathcal H^k(S\cap\sigma),\quad c=c(\delta,k,n)>0.$$
\item There is a homotopy $h\colon S\times[0,1]\to M$ with $h(s,0)=s$ and $h(s,1)=\FF(s)$. Moreover, $$\mathcal H^{k+1}(h(S,[0,1]))\leq c \mathcal H^k(S), \quad c=c(\delta,k,n)>0.$$
\end{enumerate}
\end{lem}
\begin{proof} All statements are proved by reverse induction on $m$ for the sequence of maps $\FF^m$ used in the construction of the Federer-Fleming deformation. Each step of the deformation increases $k$-dimensional volume by at most a multiplicative constant, by Lemma \ref{lem-FFlocal}.
\end{proof}
We will also need a tool to recognize which parts of $S$ are collapsed to the $(k-1)$-skeleton under the map $\FF.$
\begin{lem}\label{lem-FFvanishing}
There exists a constant $\eta=\eta(\delta,r,n,k)>0$ with the following property. Any point $x\in S$ such that $\mathcal H^k(S\cap B_M(x,r))\leq \eta$ has $\FF(x)\in M^{(k-1)}$.
\end{lem}
\begin{proof}
Let $x\in S$ and suppose $\FF(x)\not\in M^{(k-1)}$. Let $\sigma_0$ be the $k$-cell containing $\FF(x)$. Since $\diam(\tau)\leq r$ for all $\tau\in M^{[\bullet]}$, we have $\cup_{\sigma_0\subset \tau}\tau \subset B(x,r).$ The fact that $\FF(S)\cap\sigma_0$ is not collapsed to the boundary of $\sigma_0$ means that $\sigma_0\subset \FF(S)$ so
$$\mathcal H^k(\FF^k(S)\cap \sigma_0)=\mathcal H^k(\sigma_0)\geq \min_{\sigma\in M^{[k]}}\mathcal H^k(\sigma).$$
By Lemma \ref{lem-FFproperties} (1),(2) we have 
$${n+1\choose k+1}c \mathcal H^k(S\cap B(x,r))\geq \sum_{\substack{\sigma_0\subset \tau \\ \tau\in M^{[n]}}} c\mathcal H^k(S\cap \tau)\geq \mathcal H^k(\FF^k(S)\cap \sigma_0).$$
Choosing $\eta:=c^{-1}{n+1\choose k+1}^{-1}\min_{\sigma\in M^{[k]}}\mathcal H^k(\sigma),$ we see that the condition $\FF(x)\not\in \partial \sigma_0$ forces $\mathcal H^k(S\cap B(x,r))\geq \eta$. Note that the traingulation is $(r,\delta)$-uniform so $\eta$ depends only on $(r,\delta)$ but not on $M$. 
\end{proof}
We end this section with a simple but useful observation. 
\begin{lem}\label{lem-FFcodim3remainder}
Let $S\subset M$ be a closed rectifiable $k$-varifold. Then $\FF(S)=N\cup Q$ where $N\subset M^{k}$ is a closed sub-complex and $Q\subset M^{(k-1)}.$
\end{lem}
\begin{proof}
$\FF(S)\subset M^{k}$ and any simplex $\sigma\in M^{[k]}$ with interior intersecting $S$ must be contained in $\FF(S)$, because of last step of the construction. If we put $N:=\bigcup_{\sigma\in M^{[k]}, \sigma\subset \FF(S)}\sigma,$ then $Q:=\FF(S)\setminus N\subset M^{(k-1)}.$
\end{proof}
\subsection{Existence of uniformly non-degenerate triangulations}
 We explain two ways of constructing $(r,\delta)$-uniform triangulations of locally symmetric spaces. The first is very simple.
\begin{lem}\label{lem-UTcovers} Let $M$ be a compact smooth Riemannian manifold. For all $r>0$ there exists a $\delta>0$ such that all finite covers of $M$ admit $(r,\delta)$-uniform triangulations. 
\end{lem}
\begin{proof}
Choose any triangulation $\mathcal T_0$ of $M$ with no degenerate simplices and with all simplices of diameter at most $r$. Let $\delta>0$ be such that $\mathcal T_0$ is $(r,\delta)$-uniform. Then the lift of $\mathcal T_0$ to any finite cover of $M$ is still $(r,\delta)$-uniform. 
\end{proof}
The above Lemma would be enough to prove Theorem \ref{mt-collapse} for families of finite covers of a fixed manifold. While this isn't entirely unsatisfactory, it feels like an unnecessary restriction for locally symmetric spaces, where the uniformity should depend on the injectivity radius alone. We provide a way to construct $(r,\delta)$-uniform triangulations on arbitrary families of non-positively curved Riemannian manifolds of bounded geometry. Although the existence of such triangulation is intuitively clear, we weren't able to find a straightforward argument. We sketch a probabilistic construction based on Delaunay triangulations and Local Lovasz Lemma. Our construction is partly inspired by \cite{boissonnat2018delaunay}, where the goal was to construct Delaunay triangulations on any smooth compact manifold. We note that our approach uses non-positive curvature so we do not provide an alternative proof of \cite{boissonnat2018delaunay}.

\begin{defn}
Let $\Omega$ be a set of Riemannian $n$-manifolds. We say $\Omega$ is \emph{bounded geometry} if the sectional curvatures are uniformly bounded at any point of any manifold in $\Omega$ and the injectivity radii of $M\in \Omega$ are uniformly bounded away from zero.
\end{defn}

\begin{thm}\label{thm-UTboundedgeo}
Let $\Omega$ be a family of non-positively curved $n$-manifolds of bounded geometry. There exists $\delta>0$ and $r_0>0$, such that for any $r<r_0$, each $M\in \Omega$ admits a $(r,\delta)$-uniform triangulation.
\end{thm}
We need several preliminary results before we can prove it.
\begin{defn} Let $Q$ be a discrete subset of $\mathbb R^n$, with convex hull $W$. A triangulation $\mathcal T$ of $W$ with vertex set $Q$ is a \emph{Delaunay triangulation} if for any $n$-simplex $\sigma\in \mathcal T$, the smallest closed ball $B_\sigma$ containing $\sigma$ contains no points of $Q$ besides the vertices of $\sigma$.
\end{defn}
We say that a subset $Q\subset \mathbb R^n$ is in \emph{general position} if no $n+2$ points lie on a proper hyperplane or a sphere.
\begin{lem}\label{lem-DelExistence} \cite{delone1934sphere}
Let $Q\subset \mathbb R^n$ be a discrete subset in general position. The Delaunay triangulation with vertex set $Q$ exists and is unique.
\end{lem}

\begin{defn}\label{def-barysimplex}
Let $M$ be a non-positively curved Riemannian $n$-manifold and let $E\subset M$ be a subset with $|E|\leq n+1$ and $3\diam(E)< \injrad(M).$ Let $\Delta\subset \mathbb R^{|E|-1}$ be a non-degenerate affine simplex with vertex set $\{v_e\}_{e\in E}$ indexed by $E$. We define the barycentric map $b \colon \Delta\to M$ as follows. Every point $v\in \Delta$ can be written uniquely as a convex combination 
$$v=\sum_{x\in E}t_x v_x, \quad \sum_{x\in E} t_x=1,\quad 0\leq t_x, x\in E.$$
The point $b(v)$ is then defined as the unique minimum of the function 
$$y\mapsto \sum_{x\in E} t_x d(x,y)^2$$ on the $\diam(E)$-neighborhood of $E$. The uniqueness of the minimum requires the non-positive curvature assumption, which guarantees that the maps $y\mapsto d(x,y)^2$ is strictly convex \cite[Lemma 12.15]{lee2018introduction}.  One can also check that the map $b$ is differentiable. We will call $b(\Delta)$ the \emph{simplex spanned by $E$}. It does not depend on the choice of $\Delta.$ 
\end{defn}
One nice feature of the above definition is that the face of $b(\Delta)$ corresponding to a subset $E'\subset E$ is the simplex spanned by $E'$. When $M=\mathbb R^n$, then the simplex spanned by an $(n+1)$-tuple of points is just the convex hull.
\begin{lem}\label{lem-barysimplex}
Let $M$ be a nonpositively curved $n$-manifold with bounded geometry. For any $\delta>0$ there is an $\varepsilon>0$ depending only on the injectivity radius and curvature tensor of $M$ so that the following holds.  Let $E\subset M$ be a subset with $|E|=n+1$ so that $E$ is contained in an embedded geodesic ball of radius $\varepsilon$.  Let $\sigma_E$ the simplex spanned by $E$. Then the barycenter map $b\colon \Delta \to \sigma_E$ is a $(1+\delta)$-isometry, for some choice of $\Delta.$

\end{lem}

\begin{proof}

Any tuple of points $p=(p_1,..,p_{n+1})$ contained in an injectivity ball $B$ at some point $q\in M$  as above gives a barycenter map $b_{p}: \Sigma \rightarrow M$. We choose $\Sigma$ to be the affine simplex in $T_q M$ spanned by $\log_q p_i, i=1,\ldots,n+1,$ where $\log_q \colon B\to T_q M$ is the inverse of the exponential map. This map has the following continuity properties:

\begin{enumerate}
\item $b_p$ depends continuously on $p$ (i.e., on the $p_i$)
\item $b_p$ and its differential depend continuously on the metric in the $C^2$-topology
\end{enumerate}

Assume for contradiction that the conclusion of the theorem fails for some $\delta>0$, and take a sequence of bounded geometry nonpositively curved $M_k$ as in the statement of the Lemma with points $p^k=(p_1^k,..,p_{n+1}^k)$ contained in balls $B_{q_k}(\epsilon_k)$ for which the conclusion of the lemma fails.  Rescaling the metrics on the $B_{q_k}(\epsilon_k)$ so that they become balls $\overline{B}_{q_k}(\epsilon_k)$ with unit radius, it will then be the case that the metrics on the $\overline{B}_{q_k}(\epsilon_k)$ are converging to the flat metric on the n-dimensional unit disk $D$ in the $C^2$ topology.  Note that scaling the metric does not affect the barycenter map. 

By passing to a subsequential limit of the $p^k$, we obtain an $(n+1)$-tuple of points in $D$.  Since for this tuple of points the barycenter map is an isometry, this gives a contradiction with our assumption that all of the $b_{p^k}$ failed to be $(1+\varepsilon)$-isometries, in view of the continuity properties listed above.

\end{proof}


\begin{lem}[{Local Lovasz Lemma \cite{LLLspencer1977asymptotic}}]\label{lem-LLL}
Let $0<p<1, d\geq 0$. Let $A_1,\ldots, A_N$ be random events, such that $\mathbb P(A_i)\leq p$ for all $i=1,\ldots, N$ and that each $A_i$ is independent of all but $d$ other events. Suppose that $ep(d+1)<1$ where $e=\exp(1)=2.71\ldots$ is the Euler's constant. Then, the probability that none of $A_i$ occur is positive. 
\end{lem}
\begin{proof}[Proof of Theorem \ref{thm-UTboundedgeo}]
Let $\Omega$ be a bounded geometry family of Riemannian $n$-manifolds and let $M\in \Omega$. Let $\delta>0$ be a small positive parameter, to be fixed later. Choose $\varepsilon >0$ so that the conclusion of Lemma \ref{lem-barysimplex} holds. We choose $0<r_0\leq 1$ so that for any $r_1<r_0/20$ the exponent map 
$\exp_x\colon T_x \to M$ is an $(1+\varepsilon)$-isometry on the ball $\{ v\in T_x M\mid \|v\|\leq 32r_1\}, x\in M$ and $\injrad(M)\geq 64r_1$.
We choose a maximal $4r_1$-separated set $P\subset M$. Then, by maximality, $P$ is a $8r_1$-net. 

We construct a new random set $\hat P$ by independently replacing each $p\in P$ by a uniformly randomly chosen element $\hat p\in B(p,r_1).$ The new set $\hat P$ is $2r_1$-separated $10r_1$-net.
We will say that a subset $E\subset \hat P$ of $(n+1)$-points is \emph{Delaunay} if \begin{enumerate}
\item there is a unique minimal closed ball $B_E$ containing $E$.
\item the ball $B_E$ contains no other points of $\hat P.$
\end{enumerate}
Note that condition (1) is automatically satisfied for non-positively curved manifolds and subsets of diameter below the injectivity radius. Any lower cardinality subset $F\subset \hat Q$ will be declared Delaunay if it is contained in some $(n+1)$-element Delaunay set. 

\emph{ Claim 1.} Let $E$ be a Delaunay set. Then, the diameter of $E$ is bounded by $20r_1$. 

For the proof we may assume $|E|=n+1.$ The minimal closed ball $B_E$ containing $E$ contains no other points so its radius is bounded by $10r_1$. Otherwise, the center of $B_E$ would more than $10r_1$-away from $\hat P$, contradicting the fact that $\hat P$ is a $10r_1$-net. The diameter of $E$ is therefore bounded by $20r_1$, proving the claim.

For each Delaunay set $E$ let $\sigma_E\subset M$ be the simplex spanned by $E$. The exponential maps are $(1+\varepsilon)$-isometries on the balls of radius $32r_1$, which is bigger than the diameter of any Delaunay set. Hence, by Lemma \ref{lem-barysimplex}, each $\sigma_E$ is $(1+\delta)$-isometric to an affine simplex. Write ${\rm Del}(\hat P)$ for the set of simplices in $M$ spanned by the Delaunay subsets.

We will prove that for $\delta$ small, with positive probability ${\rm Del}(\hat P)$ is an $(20r_1,\delta)$-uniform triangulation of $M$. For that, we need to check three conditions.
\begin{enumerate}
\item the interiors of the simplices $\sigma_F$ are pairwise disjoint.
\item the simplices $\sigma_F\in {\rm Del}(\hat Q)$ cover entire $M$.
\item every simplex of ${\rm Del}(\hat Q)$ is $(20r_1,\delta)$-uniform. 
\end{enumerate}

For any $p\in P$ (the original set, before the random resampling) let $A_p$ be the event that one of the conditions (1) and (2) fails inside $B_M(p,8r_1)$ or that $B_M(p,28r_1)$ contains an $(n+1)$-tuple of points spanning a non-$(20r_1,\delta)$-uniform simplex. Note that any non $(20r_1,\delta)$-uniform simplex in ${\rm Del}(\hat P)$ would trigger at least one of the events $A_p$. Similarly, a failure to be a triangulation would also trigger some of the $A_p$'s, as $B_M(p,8r_1),p\in P$ cover entire $M$. Therefore, to check if ${\rm Del}(\hat P)$ is an $(20r_1,\delta)$-uniform triangulation of $M$ it is enough to show that none of the $A_p$'s occur. 
The bound on the diameter of $\sigma\in {\rm Del}(\hat Q)$ ensures that the events $A_p, A_q$ are independent if $d(p,q)> 2(20+8+1)r_1=58r_1$ and $p,q\in P$. The $8r_1$'s come from the radii, $20r_1$ from the a priori diameter bounds and $r_1$'s from the resampling radius. 

Since $P$ is $4r_1$-separated $8r_1$-net, the size of the set $\{q\in P\mid d(p,q)\leq 58r_1\}$ is bounded uniformly for all $p\in P$, by a constant $d$ dependent only on the bounds on geometry of $M$, but not on $r_1$. Therefore, each $A_p$ is independent of all but $d$ events $A_q, q\in P.$ By Lemma \ref{lem-APprob}, $\mathbb P(A_p)$ goes to zero as $\delta\to 0$. For small enough $\delta$, the assumptions of Local Lovasz Lemma (Lemma \ref{lem-LLL}) will be satisfied, so none of the events $A_p$ occur with positive probability. This proves the existence of an $(20r_1,\delta)$-uniform triangulation on $M$. Since $r_1\leq r_0/20$ was arbitary, the theorem is proved.
\end{proof}
\begin{lem}\label{lem-APprob}
\begin{enumerate}
\item The probability of the events $A_p, p\in P$ goes to $0$ uniformly as $\delta\to 0.$
\item Let $g$ be a non-positively curved Riemannian metric on $B_{\mathbb R^n}(28)$ (standard Euclidean ball around $0$) which is $(1+\delta)$-isometric to the standard flat metric. Let $Q\subset B_{\mathbb R^n}(28)$ be a $4$-separated $8$-net (with respect to the metric $g$). Let $\hat Q$ be the random set obtained by replacing each $q\in Q$ independently by a uniform random point $\hat q\in B_{(\mathbb R^n, g)}(q,1)$. Let $A$ be the event that ${\rm Del}(\hat Q)$ (defined identically as before) is not a triangulation on $B_{(\mathbb R^n,g)}(8)$ or that one of the simplices spanned by an $(n+1)$-tuple in $\hat Q\cap B_{\mathbb R^n}(28)$ has volume less than $\delta.$ Then $\mathbb P(A)\to 0$ as $r,\delta\to 0$.
\end{enumerate}
\end{lem}
\begin{proof}
In (1), everything happens inside the ball $B_M(p,28r)$. Using the inverse of the exponential map and rescaling the metric we can reduce (1) to (2).  
For (2), we only sketch the proof. Suppose that we have a sequence of metrics $g_i$ on $B_{\mathbb R^n}(28)$ which are $(1+\delta_i)$-isometric to the standard metric. Let $A^i$ be the event $A$ for $B_{\mathbb R^n}(28)$ endowed with the metric $g_i$. Consider any limiting distribution ${\rm Del}(\hat Q_\infty)$ of the collection ${\rm Del}(\hat Q_i)$, which becomes a random family of affine simplices in $B_{\mathbb R^n}(28)$, by Lemma \ref{lem-barysimplex}. In the standard metric on $\mathbb R^n$, the event that one of the simplices spanned by an $(n+1)$-tuple in $\hat Q\cap B(28)$ has volume zero is zero. Similarly, the event that $\hat Q$ is not in general position (i.e. assumptions of Lemma \ref{lem-DelExistence} are not satisfied) is zero. By Lemma \ref{lem-DelExistence}, ${\rm Del}(\hat Q_\infty)$ is a triangulation in $B_{\mathbb R^n}(8)$ almost surely. A priori, it could happen that a family of Delaunay simplices converges to an affine triangulation despite not forming a triangulation along the way. However, one can show that it could only happen if the limiting triangulation contains degenerate simplices. We know that this is a zero probability event. We can therefore apply a continuity argument to show that the probability of $A_i$ goes to $0$ as $\delta_i\to 0$. 
\end{proof}
We record the following corollary of Theorem \ref{thm-UTboundedgeo}. 
\begin{cor}\label{cor-LSFF}
Let $X$ be a symmetric space of non-compact type and let $\varepsilon>0.$ There exists $\delta>0$ such that for any small enough $r$, any compact locally symmetric space $\Gamma\bs X$ with $\injrad(\Gamma\bs X)>\varepsilon$ admits an $(r,\delta)$-uniform triangulation.
\end{cor}
\section{Hessians in symmetric spaces}\label{sec-Hess}
We compute the Hessians of the Iwasawa coordinate $H$ and the spherical functions of the form $f\circ a,$ where $a$ is the Cartan coordinate and $f\colon \mathfrak a\to \mathbb R$ is a $C^2$-function invariant under the Weyl group. 
\subsection{Iwasawa coordinate $H(g)$}\label{sec-IwasHess}
The goal of this section is to determine the Hessian of the function $H\colon G\to \fa,$ which was defined by $g=ne^{H(g)}k, n\in N, H(g)\in\fa, k\in K$. In this section we treat it as a function on $X=G/K$ and compute the Hessian on $X$. 
The list of eigenvalues of Hessian of $H$ was implicit in the computations of \cite{cmw23}. We provide the computation for completeness, as we were unable to locate it in the literature. 
\begin{prop}\label{prop-Hessian}
For any $x\in X$ there is an orthonormal basis of $T_x X$ such that the matrix of $\Hess H$ at $x$ is a diagonal $\fa$-valued matrix with diagonal entries $0$ repeated $\rk G$ times and $-H_\alpha$ for $ \alpha\in \tilde\Phi^+,$ where $H_\alpha\in \fa$ is the vector satisfying $\langle H_\alpha, Y\rangle=\alpha(Y), Y\in \fa.$  
\end{prop}
Recall that we identify $T_K X=\mathfrak g/\mathfrak k\cong \mathfrak s$ via the map $Y+\mathfrak k\mapsto \frac{1}{2}(Y-\Theta Y)$. For the remainder of this section we fix an orthonormal basis $H_1,\ldots, H_{\rk G}$ of $\fa$ and for each $\alpha\in \tilde\Phi^+$ we choose a vector $E_\alpha\in \fg_\alpha\subseteq \fn$  such that the set $\{\frac{1}{2}(E_\alpha-\Theta E_\alpha)\,|\,\alpha\in\tilde\Phi^+\}$ is an orthonormal basis of the orthogonal complement of $\mathfrak a$ in $\mathfrak s$. 

Consider the map $\psi\colon NA \to X$ given by $\psi(na)=naK\in X.$ This map is a surjective diffeomorphism. We will use $\psi$ to transport the vector fields on $NA$ to vector fields on $X$. We adopt the following notational convention. Whenever $Y\in \fn+\fa$, we will write $\mathbf Y^\circ$ for the unique $NA$ left invariant vector field on $NA$ with $\mathbf Y^\circ(1)=Y,$ and $\mathbf Y:=d\psi(\mathbf Y^\circ)$. We note that if $[Y,Z]=W$ for some $ Y,Z,W\in\fa+\fn$ then $[\mathbf Y,\mathbf Z]=\mathbf W.$ 
\begin{lem} The fields $\mathbf H_i, \mathbf E_\alpha$ form an orthonormal basis of $T_xX$ at every $x\in X$.
\end{lem}
\begin{proof}
$NA$ acts transitively and isometrically on $X$, so we only need to verify the lemma at the point $K$. 
\begin{align*}
    \mathbf E_\alpha(K)=&\frac{1}{2}(E_\alpha-\Theta E_\alpha)\in \mathfrak s,\\
    \mathbf H_i(K)=&H_i\in \mathfrak s. 
\end{align*} The Lemma follows from the conditions imposed on the vectors $H_i, E_\alpha$.
\end{proof}

The vector fields defined above will turn out to be the orthonormal eigenbasis of $\Hess H$ at every point of $X$. Since $H$ is an $\mathfrak a$ valued function, the gradient of $H$ is  an element of $C^\infty(TX\otimes \mathfrak a)$ defined by the condition 
$\langle \grad H, Y\rangle= Y H \in C^\infty(X)\otimes \fa,$ for every smooth vector field $Y\in C^\infty(TX).$
\begin{lem}\label{lem-Hgradient}
$\grad H=\sum_{i=1}^{\rk G}\mathbf H_i\otimes H_i$ and $dH=\sum_{i=1}^{\rk G} \mathbf H_i^*\otimes H_i$ where $\mathbf H_i^*(Y):= \langle H_i,Y\rangle.$ 
\end{lem}
\begin{proof}
    The two formulas are equivalent so we will just prove the second one. Note that the function $H$ restricted to $A\subset NA$ is the inverse of the exponential map, so the derivative restricted to $\mathfrak a$ is the identity map, i.e. $dH(\mathbf H_i)=H_i$ for all $i=1,\ldots, \rk G.$ As $H$ is left $N$ invariant we immediately get that $dH(\mathbf E_\alpha)=0$. This proves the lemma because $\mathbf H_i, \mathbf E_\alpha$ form a basis of $TX$ at every point.  
\end{proof}

We are now ready to compute $\Hess H\in C^\infty( T^*X\otimes T^*X\otimes \fa).$ 
\begin{lem}\label{lem-hessianinonb} $\Hess H$ is diagonalized in the basis $\mathbf H_i, \mathbf E_\alpha$ at each point of $X$. We have
\begin{align*}
&\Hess H(\mathbf H_i, \mathbf H_j)=0&\text{ for }i,j=1,\ldots, \rk G,\\
&\Hess H(\mathbf H_i,\mathbf E_\alpha)=0&\text{ for }i=1,\ldots,\rk G,\, \alpha\in\tilde\Phi^+,\\
&\Hess H(\mathbf E_\alpha, \mathbf E_\beta)=-\delta_{\alpha,\beta}H_\alpha& \text{ for }\alpha,\beta\in \tilde\Phi^+.
\end{align*}
\end{lem}
\begin{proof}
Let $Y,Z$ be vector fields on $X$. Using Koszul's identity for the Levi-Civita connection on $X$ and Lemma \ref{lem-Hgradient}, we get \begin{align*}\langle \Hess H(Y,Z), H_i\rangle=&\Hess(\langle H,H_i\rangle) (Y,Z)=\langle \nabla_Y (\grad \langle H,H_i\rangle),Z\rangle=\langle \nabla_Y \mathbf H_i, Z\rangle\\=& \frac{1}{2}\left(Y\langle \mathbf H_i,Z\rangle +\mathbf H_i\langle Y,Z\rangle-Z\langle Y, \mathbf H_i\rangle +\right.\\
&\left. \langle [Y,\mathbf H_i], Z\rangle -\langle [\mathbf H_i,Z],Y\rangle -\langle [Y,Z],\mathbf H_i\rangle\right).\end{align*}
When we apply this formula for $Y,Z\in \{ \mathbf H_j,\mathbf E_\alpha\}$, the left $NA$-invariance of the vector fields implies that all the inner products involved in the first three terms of the right hand side are constant, hence their derivatives vanish. 
\begin{align*}
\langle \Hess H(\mathbf H_j, \mathbf H_k),H_i\rangle=&\frac{1}{2}(\langle 0,\mathbf H_k\rangle - \langle 0, \mathbf H_j\rangle -\langle 0, \mathbf H_i\rangle)=0,\\
\langle \Hess H(\mathbf E_\alpha,\mathbf H_j),H_i\rangle=&\frac{1}{2}(\langle [\mathbf E_\alpha,\mathbf H_i],\mathbf H_j\rangle -\langle [\mathbf H_i,\mathbf H_j],\mathbf E_\alpha\rangle -\langle [\mathbf E_\alpha, \mathbf H_j],\mathbf H_i\rangle)\\
=&\frac{1}{2}(-\alpha(H_i)\langle \mathbf E_\alpha, \mathbf H_j\rangle- 0+\alpha(H_j)\langle \mathbf E_\alpha, \mathbf H_i\rangle)=0,\\
\langle \Hess H(\mathbf E_\alpha, \mathbf E_\beta),H_i\rangle=&\frac{1}{2}(\langle [\mathbf E_\alpha, \mathbf H_i],\mathbf E_\beta\rangle -\langle [\mathbf H_i,\mathbf E_\beta],\mathbf E_\alpha\rangle -\langle [\mathbf E_\alpha,\mathbf E_\beta], \mathbf H_i\rangle)\\
=& \frac{1}{2}(-\alpha(H_i)\langle \mathbf E_\alpha, \mathbf E_\beta\rangle -\beta(H_i)\langle \mathbf E_\beta,\mathbf E_\alpha\rangle- 0)\\
=& -\delta_{\alpha,\beta}\alpha(H_i).
\end{align*}
\end{proof}
\begin{cor}\label{cor-CompositionHessian}
Let $f\colon \fa\to \mathbb C$ be a $C^2$ function and define  $F\colon X\to\mathbb C$ by $F(g):=f(H(g))$.  Then,
\begin{align*}
&\Hess F(gK)(\mathbf H_i, \mathbf H_j)=\frac{\partial^2 f}{\partial H_i\partial H_j}(H(gK)) &\text{ for }i,j=1,\ldots, \rk G,\\
&\Hess F(gK)(\mathbf H_i,\mathbf E_\alpha)=\Hess F(gK)(\mathbf E_\alpha,\mathbf H_i)=0&\text{ for }i=1,\ldots,\rk G,\, \alpha\in\tilde\Phi^+,\\
&\Hess F(gK)(\mathbf E_\alpha, \mathbf E_\beta)=-\delta_{\alpha,\beta} \frac{\partial f}{\partial H_\alpha}(H(gK)) & \text{ for }\alpha,\beta\in \tilde\Phi^+.\\
\end{align*}
\end{cor}

\begin{cor} \label{cor-HessianofH} For any $\xi\in \fa_{\mathbb C}^*$ we have 
\begin{align*}
&\Hess (\xi H)(\mathbf H_i, \mathbf H_j)=0&\text{ for }i,j=1,\ldots, \rk G,\\
&\Hess (\xi H)(\mathbf H_i,\mathbf E_\alpha)=0&\text{ for }i=1,\ldots,\rk G,\, \alpha\in\tilde\Phi^+,\\
&\Hess (\xi H)(\mathbf E_\alpha, \mathbf E_\beta)=-\delta_{\alpha,\beta}\langle \alpha, \xi\rangle.& \text{ for }\alpha,\beta\in \tilde\Phi^+,\\
&\Hess (e^{\xi H})(\mathbf H_i, \mathbf H_j)=e^{\xi H}\xi(H_i)\xi(H_j)&\text{ for }i,j=1,\ldots, \rk G,\\
&\Hess (e^{\xi H})(\mathbf H_i,\mathbf E_\alpha)=\Hess H(e^{\xi H})(\mathbf E_\alpha,\mathbf H_i)=0&\text{ for }i=1,\ldots,\rk G,\, \alpha\in\tilde\Phi^+,\\
&\Hess (e^{\xi H})(\mathbf E_\alpha, \mathbf E_\beta)=-e^{\xi H}\delta_{\alpha,\beta}\langle \alpha, \xi\rangle.& \text{ for }\alpha,\beta\in \tilde\Phi^+
\end{align*}
\end{cor}

\subsection{Cartan coordinate $a(g)$}

The goal of this section is to explicitly compute the Hessian of $f\circ a$ for any $W$-invariant $C^2$ function $f$ on $\fa$. 
We start by defining several auxiliary vector fields on the open dense subset $KA^{++}K/K\subset X$. We keep the notation from the previous section. Let $Y\in \mathfrak k$ and $i=1,\ldots,\rk G$. For any $gK=kaK, a\in A^{++}$ put 
\begin{align*}\mathbf{W}_Y(kaK):=&\frac{\partial}{\partial t}(ke^{tY}aK)|_{t=0},\quad 
\mathbf U_i(kaK):=\frac{\partial}{\partial t}(kae^{tH_i}K)|_{t=0}.
\end{align*}
These fields are left $K$-inavaraint and a routine computation yields 
\begin{equation}\label{eq-LieBrackets}
[\mathbf W_Y, \mathbf W_Z]=\mathbf W_{[Y,Z]},\quad [\mathbf W_Y,\mathbf U_i]=0,\quad [\mathbf U_i,\mathbf U_j]=0.
\end{equation}
For any $H\in \mathfrak a^{++}$ we have $\mathbf W_Y(e^HK)=\frac{\partial}{\partial t}(e^{tY}e^HK)|_{t=0}.$ It follows that for any $k\in K$ 
\begin{align}\label{eq-WYformula}(e^{-H}k^{-1}\mathbf W_Y)(K)=(e^{-H}\mathbf W_Y)(K)=&\frac{\partial}{\partial t} (e^{-H}e^{tY}e^HK)|_{t=0}=\Ad(e^H)Y+\frak k\\
=&\frac{1}{2}(\Ad(e^H)Y-\Theta \Ad(e^H)Y)\in \mathfrak s.\nonumber\end{align} 
For the last identity we used the identification $\mathfrak g/\mathfrak k\ni Y+\mathfrak k\mapsto \frac{1}{2}(Y-\Theta Y)\in \mathfrak s.$ Similarly, for any $H\in \fa^{++},k\in K,$
\begin{align*}
(e^{-H}k^{-1}\mathbf U_i)(K)=(e^{-H}\mathbf U_i)(K)&=\frac{\partial}{\partial t}(e^{-H}e^{tH_i}e^HK )|_{t=0}=H_i. 
\end{align*}
Put $T_\alpha:=\frac{1}{2}(E_\alpha+E_{-\alpha})=\frac{1}{2}(E_\alpha+\Theta E_\alpha), \alpha\in \tilde\Phi^+.$ It is fixed by $\Theta$, so it is an element of $\fk$. By (\ref{eq-WYformula}), $$(e^{-H}k^{-1}\mathbf W_{T_\alpha})(K)=\sinh(\alpha(H))(E_\alpha-\Theta E_{\alpha}).$$ Recall that $E_\alpha$ were chosen so that $\frac{1}{2}(E_\alpha-\Theta E_\alpha),\alpha\in \tilde\Phi^+$ are pairwise orthogonal of length one. We derive the identities
\begin{align*}
\langle \mathbf W_{T_\alpha},\mathbf W_{T_\beta}\rangle(ke^H)=&4\delta_{\alpha,\beta}\sinh^2(\alpha(H)), \\
\langle \mathbf W_{T_\alpha},\mathbf U_i\rangle(ke^H)=&0,\\
\langle [\mathbf W_{T_\alpha},\mathbf W_{T_\beta}],\mathbf U_i\rangle(ke^H)=&0,\\
 \langle \mathbf U_i,\mathbf U_j\rangle(ke^H)=&\delta_{ij},
\end{align*}
for $\alpha,\beta\in\tilde\Phi^+$ and $i,j=1,\ldots, \rk G$. Finally, put $\mathbf V_\alpha:=\mathbf W_{T_\alpha}/ \|\mathbf W_{T_\alpha}\|$. We can now state the main result of this section. 
\begin{prop}\label{prop-CartanHess}
Let $F(gK):=f(a(g)).$ For any $H\in \mathfrak a^{++}, k\in K$
\begin{align*}
&\Hess F(ke^H K)(\mathbf U_i, \mathbf U_j)=\frac{\partial^2 f}{\partial H_i\partial H_j}(H)&\text{ for }i,j=1,\ldots, \rk G,\\
&\Hess F(ke^HK)(\mathbf U_i,\mathbf V_\alpha)=0&\text{ for }i=1,\ldots,\rk G,\, \alpha\in\Delta^+,\\
&\Hess F(ke^HK)(\mathbf V_\alpha, \mathbf V_\beta)=\delta_{\alpha,\beta}\frac{\cosh(\alpha(H))}{\sinh(\alpha(H))}\frac{\partial f}{\partial H_\alpha}(H)& \text{ for }\alpha,\beta\in \Delta^+.\\
\end{align*}
\end{prop}

\begin{proof}
We compute the gradient of $F(gK):=f(a(g)).$ Let $F_i(ke^H):=\frac{\partial f}{\partial H_i}(H)$.

\begin{align*}
dF(ke^H)(\mathbf W_Y)&=dF(e^H)(\mathbf W_Y)=\frac{\partial}{\partial t}f(a(e^{tY}e^H))=0,\\
dF(ke^H)(\mathbf U_i)&=dF(e^H)(\mathbf U_i)=\frac{\partial}{\partial t}f(a(e^{H+tH_i}))=\frac{\partial}{\partial H_i}f(H)=F_i(H).
\end{align*}
Hence, \begin{equation}\label{eq-gradientformula}\grad F(ke^H)=\sum_{i=1}^{\rk G}\frac{\partial f}{\partial H_i}(H) \mathbf U_i(ke^H)=\left(\sum_{i=1}^{\rk G} F_i \mathbf U_i\right)(ke^H).\end{equation}
We use Koszul's formula to evaluate the Hessian of $F.$
\begin{align*}
\Hess F(K)(\mathbf Y,\mathbf Z)=&\langle \nabla_{\mathbf Y}\grad F, \mathbf Z\rangle\\
=& \sum_{i=1}^{\rk G}\frac{1}{2}\left(\mathbf Y\langle F_i\mathbf U_i,\mathbf Z\rangle +F_i\mathbf U_i\langle \mathbf Y,\mathbf Z\rangle-\mathbf Z\langle \mathbf Y, F_i\mathbf U_i\rangle +\right.\\
&\left. \langle [\mathbf Y,F_i\mathbf U_i], \mathbf Z\rangle -\langle [F_i\mathbf U_i,\mathbf Z],\mathbf Y\rangle -\langle [\mathbf Y,\mathbf Z],F_i\mathbf U_i\rangle\right)\\
\end{align*}
We now substitute $\{\mathbf Y,\mathbf Z\}\subseteq \{\mathbf W_{T_\alpha},\mathbf U_j\}, \alpha\in \Delta^+, j=1,\ldots,\rk G.$ Write $F_{i,j}(ke^H):=\frac{\partial^2 f}{\partial H_i\partial H_j}(H).$ 
\begin{align*}
\Hess F(\mathbf U_j,\mathbf U_l)=& \sum_{i=1}^{\rk G}\frac{1}{2}\left(\mathbf U_j\langle F_i\mathbf U_i,\mathbf U_l\rangle +F_i\mathbf U_i\langle \mathbf U_j,\mathbf U_l\rangle-\mathbf U_l\langle \mathbf U_j, F_i\mathbf U_i\rangle +\right.\\
&\left. \langle [\mathbf U_j,F_i\mathbf U_i], \mathbf U_l\rangle -\langle [F_i\mathbf U_i,\mathbf U_l],\mathbf U_j\rangle -\langle [\mathbf U_j,\mathbf U_l],F_i\mathbf U_i\rangle\right)\\
=&\sum_{i=1}^{\rk G}\frac{1}{2}\left(F_{i,j}\langle \mathbf U_i,\mathbf U_l\rangle -F_{i,l}\langle \mathbf U_i,\mathbf U_j\rangle+ F_{i,j}\langle \mathbf U_i,\mathbf U_l\rangle +F_{i,l}\langle \mathbf U_i,\mathbf U_j\rangle \right)\\
=&F_{j,l}.\\
\Hess F(\mathbf U_j,\mathbf W_\alpha)=& \sum_{i=1}^{\rk G}\frac{1}{2}\left(\mathbf U_j\langle F_i\mathbf U_i,\mathbf W_\alpha\rangle +F_i\mathbf U_i\langle \mathbf U_j,\mathbf W_\alpha\rangle-\mathbf W_\alpha\langle \mathbf U_j, F_i\mathbf U_i\rangle +\right.\\
&\left. \langle [\mathbf U_j,F_i\mathbf U_i], \mathbf W_\alpha\rangle -\langle [F_i\mathbf U_i,\mathbf W_\alpha],\mathbf U_j\rangle -\langle [\mathbf U_j,\mathbf W_\alpha],F_i\mathbf U_i\rangle\right)\\
=&0.\\ \end{align*}
\begin{align*}
\Hess F(\mathbf  W_\alpha,\mathbf W_\beta)(ke^H)=& \sum_{i=1}^{\rk G}\frac{1}{2}\left(\mathbf  W_\alpha\langle F_i\mathbf U_i,\mathbf W_\beta\rangle +F_i\mathbf U_i\langle \mathbf  W_\alpha,\mathbf W_\beta\rangle-\mathbf W_\beta\langle \mathbf  W_\alpha, F_i\mathbf U_i\rangle +\right.\\
&\left. \langle [\mathbf  W_\alpha,F_i\mathbf U_i], \mathbf W_\beta\rangle -\langle [F_i\mathbf U_i,\mathbf W_\beta],\mathbf  W_\alpha\rangle -\langle [\mathbf  W_\alpha,\mathbf W_\beta],F_i\mathbf U_i\rangle\right)\\
=&\sum_{i=1}^{\rk G}\frac{1}{2}F_i\mathbf U_i\langle \mathbf  W_\alpha,\mathbf W_\beta\rangle=2\delta_{\alpha\beta}\sum_{i=1}^{\rk G}F_i\mathbf U_i(\sinh(\alpha(H))^2)\\
=&4\delta_{\alpha\beta}\sum_{i=1}^{\rk G}F_i(ke^H) \alpha(H_i)\sinh(\alpha(H))\cosh(\alpha(H))\\
=&4\delta_{\alpha,\beta}\frac{\partial f}{\partial H_\alpha}(H)\sinh(\alpha(H))\cosh(\alpha(H)).
\end{align*}
We deduce that 
\begin{align*}
\Hess F(\mathbf U_i, \mathbf U_j)(ke^HK)=&\frac{\partial^2 f}{\partial H_i\partial H_j}(H),\\
\Hess F(\mathbf U_i, \mathbf V_\alpha)(ke^HK)=&\frac{1}{\|\mathbf W_{\alpha}\|}\Hess F(\mathbf U_i, \mathbf W_\alpha)(ke^H)=0,\\
\Hess F(\mathbf V_\alpha, \mathbf V_\beta)(ke^HK)=&\frac{1}{\|\mathbf W_{\alpha}\|\|\mathbf W_\beta\|}\Hess F(\mathbf W_\alpha,\mathbf W_\beta)(ke^H)\\
=& \delta_{\alpha\beta}\frac{\partial f}{\partial H_\alpha}(H)\frac{\cosh(\alpha(H))}{\sinh(\alpha(H))}.
\end{align*}
\end{proof}
\begin{cor} The eigenvalues of $\Hess F(ke^H)$ are 
$$\left\{\lambda_1,\ldots,\lambda_{\rk G},\frac{\partial f}{\partial H_\alpha}(H)\frac{\cosh(\alpha(H))}{\sinh(\alpha(H))},\alpha\in\tilde\Phi^+\right\}$$ where $\lambda_1,\ldots,\lambda_{\rk G}$ are the eigenvalues of $\Hess f(H).$ In case $H$ is singular, and $\alpha(H)=0$, the value $\frac{\partial f}{\partial H_\alpha}(H)\frac{\cosh(\alpha(H))}{\sinh(\alpha(H))}$ should be replaced by $\frac{2}{\|\alpha\|^2}\frac{\partial^2 f}{(\partial H_\alpha)^2}(H).$
\end{cor}
\begin{proof} If $H\in\mathfrak a^{++}$, the corollary follows directly from Proposition \ref{prop-CartanHess}. To handle the singular cases we use the fact that $F$ is $C^2$ and take the limit as $\alpha(H)\to 0.$
\end{proof}

\subsection{Plurisuperharmonic Functions} \label{section:naturalvectorfields}

In this section we define, for $\Lambda$ an infinite covolume discrete subgroup of $G$, $\Lambda$-equivariant functions on $X$ whose gradient flows contract volume in low codimensions.  These flows were recently constructed by Connell-McReynolds-Wang, who called them natural flows \cite{cmw23}.  It would be possible to use their work as a black box, but for the sake of completeness we give the full construction here. These flows will be crucial in the proof of the topological waist inequalities (Theorems \ref{intro:critical exponent}, \ref{mt-topwaist}, and \ref{mt-finiteorbit}).

\begin{defn}\label{defn-criticalexp} Let $G$ be a semisimple real Lie group and let $\Lambda \subset G$ be a discrete subgroup. First assume that $G$ has rank one.  The critical exponent $\delta(\Lambda)$ is defined as the exponent of convergence of the Poincare series 
$$P_\Lambda(s):=\sum_{g\in \Lambda}e^{-s \alpha(a(g))},$$ where $a(g)$ is the central Cartan coordinate and $\alpha$ is the simple root of $G$. 
\end{defn}
The invariant $\delta(\Lambda)$ measures the rate of growth of $\Lambda$ inside $G$. Lattices in $G$ have always the same critical exponent, equal to $\delta(G):=m_\alpha+2m_{2\alpha}$. This is the maximal critical exponent achieved by any discrete subgroup of $G$. For $G=\Sp(n,1), F_4^{(-20)}$ the maximal critical exponent is achieved only by lattices and there is a gap between $\delta(G)$ and the critical exponents of infinite covolume discrete subgroups. This was proved by Corlette \cite{corlette1990}, using quantitative property (T) for $\Sp(n,1)$ and $F_4^{(-20)}$ proved by Kostant \cite{kostant1969existence}. 
\begin{thm}[{\cite[Thm 4.4]{corlette1990}}] \label{asdfj}
\begin{enumerate} 
\item Suppose $\Lambda$ is a discrete subgroup of $G=\Sp(n,1)$. Then either $\delta(\Lambda)\leq 4n$ or $\Lambda$ is a lattice and $\delta(\Lambda)=4n+2$.
\item Suppose $\Lambda$ is a discrete subgroup of $G=F_4^{(-20)}$. Then either $\delta(\Lambda)\leq 16$ or $\Lambda$ is a lattice and $\delta(\Lambda)= 22$.
\end{enumerate}
\end{thm}
In higher rank groups the analogue of critical exponent is the \emph{growth indicator function} defined by Quint \cite{quint2002}. 
\begin{defn}
For any open cone $\mathcal C\subset \mathfrak a^+$ let $\tau_\mathcal C$ be the exponent of convergence of the series 
$$\sum_{\substack{g\in \Lambda\\ a(g)\in\mathcal C}}e^{-s\|a(g)\|}.$$
The value of the growth indicator function $\Psi_\Lambda\colon \mathfrak a^+\to \mathbb R\cup\{-\infty\}$ at a point $H\in \mathfrak a^+$ is defined as 
$$\psi_\Lambda(H):=\|H\|\inf_{\mathcal C \ni H}\tau_\mathcal C.$$
\end{defn}
The infimum here is taken over all open cones $\mathcal{C}$ containing $H$. 
 The growth indicator function is homogeneous, that is, $\psi(tH)=t\psi(H)$ for any $t>0$ and concave upper semi-continuous \cite{quint2002}. It is non-negative on the \emph{limit cone} $\mathcal L$ \cite{benoist97}. Moreover, if $\Lambda$ is Zariski dense, then the limit cone $\mathcal L$ has non-empty interior \cite{benoist97}. In the rank one cases, the critical exponent and the growth indicator function are related by the identity
 
 \begin{equation} \label{criticalexponentgrowthindicator}
 \psi_\Lambda=\delta(\Lambda)\alpha,
 \end{equation}
 where $\alpha$ is the simple root.
We say that an element $\xi\in \mathfrak a^*$ is \textit{tangent to} the growth indicator function if $\xi\geq \psi_\Lambda$ on $\mathfrak{a}^+$ and  there exists an $H\in\mathfrak a^{++}$ with $\psi_{\Lambda}(H)=\xi(H).$ 

Let $\mathcal{F} := G/P = K /M $ be the Furstenberg boundary of $G$, and let $\Lambda < G$ be a discrete subgroup. At this point we need to remind that our convention with $H$ is different from \cite{Quint2002Mesures} or \cite{edwards2023temperedness}. If we define $\tilde H(g)$ by $g=ke^{\tilde H(g)}n$, then the Iwasawa coordinate of \cite{Quint2002Mesures,edwards2023temperedness} is $\tilde H$. To obtain the formulas in terms of ``our" $H$ one needs to use the formula $H(g)=-\tilde H(g^{-1}).$

\begin{defn} Let $\xi\in\frak a^*$. A probability measure $\nu$ on $\mathcal{F}$ is $(\Lambda,\xi)$-conformal if 
$$\frac{d \gamma _*\nu}{d\nu}(kP)=e^{\xi(H(k^{-1}\gamma))},$$ for all $\gamma \in \Lambda.$
\end{defn}

Here $\gamma _*\nu(B)= \nu (\gamma^{-1}(B))$ for any Borel subset $B \subset \mathcal{F}$. We will need the following result from Quint \cite{Quint2002Mesures}.  
\begin{thm} \label{thm-existenceofpattersonsullivan}
Let $\Lambda$ be a discrete subgroup of $\Isom(X)$. Then if $\xi$ is tangent to $\psi_\Lambda$, there exists a $(\Lambda,\xi)$-conformal measure.  
\end{thm}

  We have to recall one more result due to Quint and Lee-Oh. We say that a subset $\Omega\in \Phi^+$ is strongly orthogonal if for every $\alpha,\beta\in \Omega$ we have $\pm \alpha\pm \beta\not\in \Omega.$ 
\begin{thm}[{\cite[Thm 7.1]{LeeOh24}}, \cite{oh02}] \label{thm-indfctiongap}
Let $G$ be semisimple Lie group with no rank one factors. Let $\Theta$ be the half sum of the roots in a strongly orthogonal subset of $\Phi^+$. Then, for any discrete $\Lambda\subset G$ which is not a lattice 
$$\psi_\Lambda\leq 2\rho-\Theta.$$
\end{thm}

For $\xi \in \mathfrak{a}^*$ and $\nu$ a $(\Lambda, \xi)$-invariant conformal measure the function $\EO_\nu$ associated to $\nu$ is defined by the formula 
\begin{equation} \label{formula:br}
\EO_{\nu}(g):=\int_K e^{\xi(H(k^{-1}g))}d\nu(kP)
\end{equation}
for $g\in G$ (see \cite[Section 3]{edwards2023temperedness}.) Since $\EO_\nu$ is left $\Lambda$-invariant and right $K$-invariant, it descends to a function on $\Lambda \backslash X$ which by abuse of notation we also denote by $\EO_{\nu}$. It is the density of the Burger-Roblin measure projected from $\Lambda\bs G$ to $\Lambda\bs X$ (see \cite[Lemma 9.1]{edwards2023temperedness}). One can check that $\EO_{\nu}$ is an eigenfunction for the Laplacian (\cite[Lemma 3.4]{edwards2023temperedness}). We will be interested in determining ranges of dimensions $k$ for which the gradient of $\EO_{\nu}$ contracts the $k$-dimensional volume. 

To this end, for $\xi \in \mathfrak{a}^*$, we define $\tau_k(\xi)$ to be the sum of the k-largest elements (possibly with multiplicity) of the multiset  
\begin{equation} \label{listofevs}
\{ ||\xi||^2,0,..,0\} \cup \{ -\langle \alpha , \xi \rangle | \alpha \in \Phi^+ \}, 
\end{equation}
\noindent where the number of 0's in the first set in the union is $\text{rank}(G)-1$. We will see in a moment that this multiset is exactly the eigenvalues of the Hessian of $g \mapsto e^{\xi(H(k^{-1}g))}$. Let $\Omega_d= \{\xi \in \mathfrak{a}^* | \tau_{\text{dim}(X)-d}(\xi) <0 \}$.  We do not undertake here a close study of the sets $\Omega_d$.  The following lemma records a few basic properties and will be more than sufficient for our purposes.   

\begin{lem} \label{lem-omegakconvex}
$\Omega_d$ is a convex open subset of $\mathfrak{a}^*$.  Moreover, $\Omega_d \cup \{0\}$ is star-shaped with respect to the origin $0$.  
\end{lem}

\begin{proof}
That $\Omega_d$ is open and that $\Omega_d\cup \{0\}$ is star-shaped with respect to $0$ is clear.  We only need to check that $\Omega_d$ is convex.  That $\Omega_d$ is convex follows from the fact that $$||\lambda \xi_1 + (1-\lambda) \xi_2||^2 \leq \max (||\xi_1||^2, ||\xi_2||^2).$$ 

and 

\[
-\langle \alpha, \lambda \xi_1 + (1-\lambda)\xi_2 \rangle \leq \max( -\langle \alpha, \xi_1\rangle, -\langle \alpha, \xi_2 \rangle) 
\]

for $\lambda \in (0,1)$.  
\end{proof}

The \textit{k–trace} $\tr_k(Q)$ of a bilinear form Q is the sum of the k largest eigenvalues of Q. The $k$-Jacobian ${\rm Jac}_k(T)$ of a linear map $T:V\rightarrow W$ between inner product spaces is the largest value, over all $k$-dimensional subspaces $V_0$ of $V$  for the absolute value of the Jacobian of the map $T|_{V_0} : V_0 \rightarrow T(V_0)$.  We will use the following result from \cite{cmw23}, which says that flowing along the gradient of a function whose Hessian has negative $k$-traces exponentially decreases the volume of any k-dimensional submanifold. 

\begin{thm} \cite[Theorem 1.1]{cmw23} \label{contraction}
 \noindent Let $\phi\colon M \rightarrow \mathbb{R}$ be any $C^2$–smooth function
whose associated $C^1$–smooth vector field $\grad \phi$ is complete. Let $\Phi_t: M \rightarrow M$ be the associated one-parameter subgroup of self-diffeomorphisms that the vector field integrates to. If $k$ is a positive integer and 
\begin{equation} \label{defn:psph}
\tr_k(\Hess (\phi)) \leq -\varepsilon 
\end{equation} 
\noindent for all $x \in M$, then
\[
{\rm Jac}_k(\Phi_t(x)) \leq e^{-\varepsilon t}
\]
for all $t>0$ and $x \in M$. In particular, if $\epsilon>0$ then $\Phi_t$ uniformly exponentially contracts the k–volume
as t increases.
\end{thm}

We note that the fact that these flows decrease the volumes of $k$-submaniflolds exponentially fast is also a consequence of the first variation formula for gradient flows (\ref{eq-statHessianBd}) and Gr\"onewall's inequality. 
\noindent We call a function that satisfies (\ref{defn:psph})  $(k,\epsilon)$\textit{-plurisuperharmonic} ($(k,\epsilon)$-PSpH) or $k$-plurisuperharmonic ($k$-PSpH) if we just require the negativity of $k$-traces. The set  of $k$-plurisuperharmonic functions is a convex cone. On general Riemannian manifolds the related $k$-plurisubharmonic functions (which are minus $k$-plurisuperharmonic) were introduced and studied by Harvey and Lawson Jr (see \cite{harvey2012geometric}).

\begin{prop} \label{prop-brpsph}
Suppose that $\nu$ is a $(\Lambda,\xi)$-conformal measure on $\mathcal{F}=G / P$ with $\tau_k(\xi)<0.$  Then the function $\EO_{\nu}$ is $k$-PSpH.  
\end{prop}
\begin{proof}
Note that for any fixed $f\in G$ the eigenvalues for the Hessian of $gK \mapsto e^{\xi(H(f^{-1}g))} $  are given by the list (\ref{listofevs}) by the computations in Section \ref{sec-IwasHess} and Corollary \ref{cor-HessianofH} (see the proofs of Propositions \ref{prop-octonionicr(x)} and \ref{prop-novereight} further down for more details on the octonionic hyperbolic and $\SL_n(\mathbb R)$ cases). Then the statement follows from the fact that $\EO_{\nu}$ is obtained by averaging functions of this form.  

\end{proof}

\begin{defn} \label{defn-CXk}
For rank one symmetric spaces $X$, we define the quantity $C_X(d)$ to be the supremal $\delta>0$ so the function $E_\nu$ 
for every discrete subgroup $\Lambda$ of $\Isom(X)$ with critical exponent smaller than $\delta$ is $(n-d)$-PSpH. It is straightforward to get a lower bound for $C_X(d)$ by computing the Hessians of exponentials of negative multiples of Busemann functions for $X$ as in the proof of Proposition \ref{prop-octonionicr(x)} below. See Table \ref{tab-cxd} for explicit estimates. For the real hyperbolic space $\mathbb H^n_\mathbb R$ one can check that $C_X(d)=n-d-1$, which matches the exponent of the growth for the volumes of balls in $(n-d)$-dimensional hyperbolic space (with the constant curvature $-1$ metric).  
\end{defn}

\begin{table}[ht]
\centering
\begin{align*}
C_{\mathbb H^n_\mathbb R}(d)=&\begin{cases}n-1-d & d=0,\ldots, n-1\\ 0 & d=n,\end{cases}\\
C_{\mathbb H^n_{\mathbb C}}(d)\geq &\begin{cases}2n-2d & d=0,1\\ 2n-1-d & d=2,\ldots, 2n-1\\
0 & d=2n,\end{cases}\\
C_{\mathbb H^n_\mathbb H}(d)\geq &\begin{cases} 4n+2-2d & d=0,1,2,3\\ 
4n-1-d & d=4,\ldots, 4n-1\\ 0 & d=4n,\end{cases}\\
C_{\mathbb H^2_\mathbb O}(d)\geq &\begin{cases} 22-2d & d=0,1,\ldots,7\\ 
15-d & d=8,\ldots, 15\\ 0 & d=16,\end{cases}\\
\end{align*}
\caption{$C_X$ for rank one spaces.}
\label{tab-cxd}
\end{table}
\subsection{Dimension gap $r(X)$} \label{computingr(x)}

Let $X$ be a symmetric space. The quantity $r(X)$ is defined to be the greatest $d$ so that every infinite co-volume quotient $\Gamma \backslash X$ admits an $(n-d)$-PSpH function. If $\Lambda$ has torsion we ask for a $\Lambda$-invariant $k$-PSpH function on $X$, to avoid dealing with technicalities at the singular points of the orbifold $\Lambda\bs X$). In this section, we obtain lower bounds for $r(X)$ in the octonionic hyperbolic and $\SL_n(\mathbb R)$ cases. Similar computation could be carried for classical groups like $\SO(p,q)$, but we limit ourselves to those two cases in the interest of space.

\begin{prop} \label{prop-octonionicr(x)}
$r(\mathbb H^2_\mathbb O)\geq 2.$
\end{prop}

\begin{proof}
Recall that if $\alpha$ is the simple root of $\rm{F}_4^{(-20)}$, then the positive roots are $\{\alpha,2\alpha\}$ and the multiplicities are $m_{\alpha}=8$ and $m_{2\alpha}=7$.  
By Theorem \ref{asdfj}, we know that the critical exponent $\delta(\Lambda)$ of an infinite covolume discrete subgroup $\Lambda$ of $G=F_4^{(-20)}$ is at most 16.  

For $s \in \mathbb{R}$, define $F:X = G/K \rightarrow \mathbb{R}$ by 
\[
F_s(gK):= e^{s \alpha (H(g))}
\]

\noindent where recall $g= n e^{H(g)}k$.  Then relative to the orthonormal basis from Section \ref{sec-IwasHess}, it follows from Corollary \ref{cor-HessianofH} that $\Hess(F_s)/(s F_s)$  is equal to the diagonal matrix  
\[
 \begin{bmatrix}
    s & & & & & & \\ 
    & -1 &&&&& \\ 
    && \ddots &&&&\\ 
    &&& -1 &&& \\ 
    &&&& -2 && \\ 
   &&&&&\ddots &\\
   &&&&&& -2  \\ 
  \end{bmatrix},  
\]
\noindent where there are eight -1's and seven -2's on the diagonal.  For $s \leq 16$, one can check that the $k$-trace is negative for $k=16,15,14$. We now let $\nu$ be a $(\delta(\Lambda),\Lambda)$-conformal density, which exists by Theorem \ref{thm-existenceofpattersonsullivan}. The  function $\EO_\nu$ (\ref{formula:br}) is obtained by averaging translates of $F$, so it is $(n-d)$-PSpH, for $d=0,1,2.$ This shows that $r(\mathbb H^2_\mathbb O) \geq 2$.

\end{proof}

 \begin{prop} \label{prop-novereight}
 $r(\SL_n(\mathbb R)/ \SO(n)) \geq  \lfloor n/8 \rfloor -1 $.  
\end{prop}

\begin{proof}
Denote the positive roots of $\SL_n(\mathbb R)$  by $\Phi^+=\{\alpha_{i,j} | 1 \leq i < j \leq n \}$  and the simple roots by $\Delta^+= \{ \alpha_{i,i+1}| i =1,..,n-1 \}$.  Here $\alpha_{i,j}$ evaluated at an element 
\[
\begin{bmatrix}
    e_1 & &\\
   
    & \ddots &  \\ 
    && e_n \\ 
  \end{bmatrix}  
\]
\noindent of $\mathfrak{a}$ is equal to $e_i-e_j$.   For $\xi \in \mathfrak{a}^*$ a nonnegative linear combination of simple roots, define $F\colon X \rightarrow \mathbb{R}$ by 
\[
F(gK):= e^{\xi (H(g))}
\]

Then relative to the orthonormal basis from Section \ref{sec-IwasHess}, by Corollary \ref{cor-HessianofH}  $\Hess(F)/F$  is equal to the diagonal matrix  
\[
 \begin{bmatrix}
    ||\xi||^2 & & & & & &\\
    & 0 && & &&\\
    && \ddots &&&& \\ 
    &&& 0 &&& \\ 
    &&&& -\langle \alpha_1 , \xi \rangle  && \\ 
    &&&&& \ddots &\\ 
    &&&&&& -\langle \alpha_{|\Phi^+|}, \xi \rangle \\ 
  \end{bmatrix}.  
\]
Here there are $n-2$ upper-left diagonal entries equal to zero, and the $\alpha_1,..,\alpha_{|\Phi^+|}$  are the elements of $\Phi^+$.  

Recall that $\rho$ is the half sum of the multiset of positive roots, and $\Theta$ is the half sum of the roots in a maximal strongly orthogonal subset of $\Phi^+$, which we choose to be $\{\alpha_{1,n},\alpha_{2,n-1},\ldots\}$.    By Theorem \ref{thm-indfctiongap}, the growth indicator function $\psi_{\Lambda}$ of every infinite covolume discrete subgroup $\Lambda \subset \Isom(X)$ satisfies $\psi_{\Lambda} \leq 2\rho -\Theta$. By Proposition \ref{prop-brpsph}, the formula (\ref{formula:br}) for the  functions $E_\nu$, and the definition of $\Omega_d$ (see (\ref{listofevs})), it is therefore enough to check that $\xi = 2\rho - \Theta$ is in $\Omega_d$ for $d \leq \lfloor n/8 \rfloor -1$.

In our case, we have that 
\[
2\Theta= \sum_{i=1}^{\lfloor n/2 \rfloor} \alpha_{i,n-i}= e_1 +.. + e_{\lfloor n/2 \rfloor} - e_{\lfloor n/2 \rfloor +1 }-..-e_n
\]
and 
\[
2 \rho = (n-1)e_1 + (n-3)e_2 + ...+ (1-n) e_n.  
\]

We give the computation in the case when $n$ is even; the computation in the case when $n$ is odd, which we omit, is similar and yields the same lower bound. Let $n=2m$.  We have 

\begin{align*}
	||2 \rho - \Theta||^2 &= 2((n-3/2)^2 + (n-7/2)^2 +  .. + (1/2)^2) \\
	&=  2 \sum_{k=0}^{m-1} \frac{1}{4}(1+4k)^2 = \frac{1}{2}\left( m + 8\left(\frac{m(m-1)}{2}\right) + 16\frac{(m-1)m(2m-1)}{6}\right).   
	\end{align*}

  \noindent We also have
\[
\langle \alpha_{i,j}, 2\rho -\Theta \rangle = 2(j-i) -  \mathbb{I}_{\substack{i< \lfloor n/2 \rfloor,\\j> \lfloor n/2 \rfloor}}(i,j)  
\]
and 
\begin{align*}
\sum_{j>i} \langle \alpha_{i,j}, 2\rho -\Theta \rangle &= \langle 2\rho, 2\rho-\Theta \rangle \\ 
&= \langle 2\rho, 2\rho \rangle - \langle 2\rho, \Theta \rangle \\
&= 2(1^2 + 3^2 +.. + (n-1)^2) - 2( \frac{1}{2} + .. + \frac{n-3}{2} + \frac{n-1}{2} )) \\ 
&= 2 \left(\frac{m(2m-1)(2m+1)}{3}\right) - m^2. 
\end{align*}

The trace of $\Hess(F)/F$ is then equal to 
\begin{align*}
&|| 2\rho - \Theta||^2 -\sum_{j>i} \langle \alpha_{i,j}, 2\rho -\Theta \rangle  \\
&=\left(\frac{m}{2} + 2m(m-1) + \frac{4(m-1)m(2m-1)}{3} \right) - \left(\frac{2m(2m-1)(2m+1)}{3} - m^2\right)  \\
&= 2m^2 - \frac{3}{2}m + \frac{4}{3}(2m^3 - 3m^2 + m ) - \frac{2}{3} (4m^3 - m )  + m^2\\ 
&=  - (2-4 + 1)m^2 + (10/3 - 3/2) m= -m^2 + \frac{11}{6}m. 
\end{align*}

We now need to check what is the largest $d$ is so that if we omit the $d$ largest possible values of  $\langle \alpha_{i,j}, 2\rho -\Theta \rangle$
from $\sum_{j>i} \langle \alpha_{i,j}, 2\rho -\Theta \rangle$ we still get a negative value when we subtract it from $|| 2\rho - \Theta||^2$.  It will then be the case that the trace of $\Hess(F)/F$ through every $(n-d)$-dimensional tangent plane is negative.  

The largest possible value of $\langle \alpha_{i,j}, 2\rho -\Theta \rangle$ is $2(n-1) -1 = 2n-3=4m-3$.  We can thus omit at least 
\[
\left\lfloor \frac{m^2 - \frac{11}{6}m}{4m-3} \right\rfloor \geq \left\lfloor \frac{m^2 - \frac{11}{6}m}{4m} \right\rfloor \geq \left\lfloor \frac{1}{4}m - \frac{11}{24}  \right\rfloor = \left\lfloor \frac{1}{8}n -\frac{11}{24} \right\rfloor \geq \left\lfloor \frac{1}{8}n  \right\rfloor -1 
\]
terms from the sum and still obtain a negative value.

\end{proof}

\section{Monotonicity estimates in rank one symmetric spaces}\label{sec-mono}

In Section \ref{sec-varifolds}, we have fixed a smooth non-increasing function $\chi\colon \mathbb R\to \mathbb R_{\geq 0}$ such that $\chi(t)=1$ for $t<0$, $\chi(t)=0$ for $t>1$ and $\chi'(t)\geq -2$ for all $t\in\mathbb R$. We have also defined $\chi_r(t):=\chi(t-r).$ We keep this choice and notation fixed in the remainder of the paper.

Let $G$ be a rank one simple real Lie group. We keep the notations from Section \ref{sec-symspaces} and \ref{sec-decay}. 
Let  $\kappa\colon \{1,\ldots, \dim X-1\}\to \mathbb R_{\geq 0}$, be given as
\begin{table}[ht]
\centering
\begin{tabular}{ l l |c |c }
 & & $X$ & $\dim X$\\ \hline
 $\kappa(k)=k-1$ & & $\mathbb H_\mathbb R^n$ & $n$\\
 $\kappa(k)=\begin{cases}k-1 \\ 2n \end{cases}$ & $\begin{array}{l} k=1,\ldots, 2n-1 \\  k=2n-1 \end{array}$ & $\mathbb H_\mathbb C^n$ & $2n$\\
 $\kappa(k)=\begin{cases}k-1\\ 4n-4+2i \end{cases}$ & $\begin{array}{l} k=1,\ldots, 4n-3 \\ k=4n-3+i  \end{array}$ &$\mathbb H_\mathbb H^n$ & $4n$\\
$\kappa(k)=\begin{cases}k-1\\ 8+2i\end{cases}$ & $\begin{array}{l} k=1,\ldots,9 \\ k=9+i  \end{array}$ & $\mathbb H_\mathbb O^2$ & $16$\\
\end{tabular}
\caption{$\kappa$ in rank one groups.}
\label{tab-kappa}
\end{table}

We define a spherical function $\ff\colon G\to \mathbb R_{\geq 0}$ by 
\begin{equation}\label{defn-ffdef}
\ff(g):= \frac{1}{2\|\alpha\|}\log (2\cosh(2\alpha(a(g)))), 
\end{equation} 
where $\alpha\in \Phi^+$ is the simple root. This function is a smooth approximation to the distance to the origin $K\in X$, which is given by 
$$d(gK,K)=\|a(g)\|=\frac{\alpha(a(g))}{\|\alpha\|}.$$

\begin{lem}\label{lem-HessianRk1}
For any $x\in X$, any $k$-dimensional subspace $W\subset T_x X,$ we have 
$$ \tr \Hess \ff|_W\geq \kappa(k)\|\alpha\|.$$
\end{lem}
\begin{proof}
Let $H_1\in \mathfrak a^+$ be the element with $\|H_1\|=1.$ By Proposition \ref{prop-CartanHess}, the eigenvalues of $\Hess \ff$ at a point $ke^HK\in X, k\in K, H=tH_1\in \mathfrak a^{++}$ are 
\begin{align*}
\frac{\partial^2}{(\partial H_1)^2}\frac{1}{2\|\alpha\|}\log(2\cosh(2\alpha(H)))=& \frac{\|\alpha\|}{2}\frac{\partial^2}{(\partial t)^2}\log\cosh(2t)= \|\alpha\|\frac{\partial}{\partial t}\frac{\sinh(2t)}{\cosh(2t)}\\
=&2\|\alpha\|(1-\tanh^2(2t))>0, \text{ with multiplicity }1,\\
\frac{\cosh(\alpha(H))}{\sinh(\alpha(H))}\frac{\partial}{\partial H_\alpha}\frac{1}{2}\log\cosh(2\alpha(H))=&\|\alpha\| \frac{\cosh(\alpha(H))\sinh(2\alpha(H))}{\sinh(\alpha(H))\cosh(2\alpha(H))}\\
\geq& \|\alpha\| \text{ with multiplicity }m_\alpha,\\
\frac{\cosh(2\alpha(H))}{\sinh(2\alpha(H))}\frac{\partial}{\partial H_\alpha}\frac{1}{2}\log\cosh(2\alpha(H))=& 2\|\alpha\| \text{ with multiplicity }m_{2\alpha}.\\
\end{align*}
The trace $\tr \Hess \ff(x)|_W$ is bounded from below by the sum of bottom $k$ eigenvalues of $\Hess \ff(x)$. Therefore, for $k=1$ the lower bound is $0$, for $k=1+i, i\leq m_\alpha$, the lower bound is $i\|\alpha\|$ and for $k=1+m_{\alpha}+i, i<m_{2\alpha}$, the lower bound is $(1+m_\alpha+2i)\|\alpha\|$. Analyzing these bounds case by case, we arrive at the lower bound $\tr \Hess \ff(x)|_W\geq \kappa(k)\|\alpha\|$.
\end{proof}
Let $S\subset X$ be a varifold, let $gK\in G/K=X$ and let $r\geq 0$. Put 
\begin{equation}
v_S(gK,r):= \int_S \chi_r(\ff(s^{-1}g))d\mathcal H^k(sK)=\int_S \chi_r(\ff(g^{-1}s))d\mathcal H^k(sK).
\end{equation}
The value $v_S(gK,r)$ roughly corresponds to the mass of $S$ in the $r$-ball around $gK$. Because of the smoothing effect of $\chi$ it is much more convenient to work with than $\mathcal H^k(S\cap B(gK,r)),$ which is the actual mass in the $r$-ball.
\begin{lem}\label{lem-VarifoldMono}
\begin{enumerate}
\item Let $S\subset X$ be a stationary $k$-varifold. For any $x\in X$, $r>s>0$
$$v_S(gK,r)\geq e^{\kappa(k)\|\alpha\|(r-s)}v_S(gK,s).$$
\item There is $\delta>0$ such that $v_S(x,0)\geq \delta$ for any stationary $k$-varifold $S\subset X$ and $\mathcal H^k$-almost every point of $S$.
\end{enumerate}
\end{lem}
\begin{proof}
Part (1) follows from Lemma \ref{lem-abstmonotonicity} applied to function $\mathfrak f.$ Part (2) follows from Lemma \ref{lem-smallballmono} and our choice of $\chi.$
\end{proof}
In \cite[Theorem 1]{anderson1982complete} Anderson proved that on any complete simply connected negatively curved manifold $N$ of sectional curvature at most $-a$ the volume of the intersection of $B(x,r)\cap S,$ where is $S$ is a stationary $k$-varifolds and $x$ is a regular point of $S$, is at least equal to the volume of an $r$-ball in the $k$-dimensional contractible space with constant sectional curvature $-a.$ Up to a multiplicative constant, our estimate is comparable with Anderson's for $\mathbb H^n_\mathbb R$ but it is strictly stronger for $\mathbb H^n_\mathbb C, \mathbb H^n_\mathbb H, \mathbb H^2_\mathbb O$ for $k> 2n-1, 4n-3, 9$ respectively. This seemingly modest improvement is at the core of our arguments in Sections \ref{sec-VarifoldsMonoMC} and \ref{sec-collapsing}. 

\section{Stationary Varifolds and matrix coefficients}\label{sec-VarifoldsMonoMC}
Let $\Gamma$ be a discrete subgroup of $G$. In this section we will use the monotonicity estimates to control the decay of certain matrix coefficients of the right regular representation $L^2(\Gamma\bs G)$. The discussion applies to all discrete subgroups of $G$ although the corollaries seem to be most interesting in the case of lattices. We use the tension between the monotonicity estimates and the decay of matrix coefficients to prove Theorem \ref{mt-stationary}.

\subsection{Convolutions, unitary representations and spherical functions.} 
The discussion in this part applies to all semi-simple real Lie groups with no restrictions on the rank. Let us recall the definition of the convolution for functions on $C(G)$. For $\varphi_1,\varphi_2\in C_c(G)$ we put
$$ \varphi_1\ast \varphi_2(h)=\int_G \varphi_2(g^{-1}h)\varphi_1(g)dg=\int_G \varphi_1(hg^{-1})\varphi_2(g)dg.$$
The same definition applies when only one of the functions is in $C_c(G)$ and the other is only locally $L^1$-integrable or when both functions are $L^2$-integrable. In particular, if $f\in L^2(\Gamma\bs G)$, then $f\ast \varphi\in L^2(\Gamma\bs G)$ for any $\varphi\in C_c(G).$ Let us explain the interplay between convolutions and right regular unitary representations in more detail.
Let $\lambda_{\Gamma\bs G}\colon G\to \mathscr U (L^2(\Gamma\bs G))$ be the right regular representation 
$$\lambda_{\Gamma\bs G}(g)f(h):=f(hg), \quad f\in L^2(\Gamma\bs G).$$
Note that $\lambda_{\Gamma\bs G}$ is a left action $\lambda_{\Gamma\bs G}(g_1)\lambda_{\Gamma\bs G}(g_2)f=\lambda_{\Gamma\bs G}(g_1g_2)f.$
For any $\varphi\in L^1(G)$ we define 
$$[\lambda_{\Gamma\bs G}(\varphi)f](h) := \int_{G}\varphi(g)(\lambda_{\Gamma\bs G}(g)f(h))dg, \quad f\in L^2(\Gamma\bs G), h\in 
 G.$$ The resulting function $\lambda_{\Gamma\bs G}(\varphi)f$ is left $\Gamma$-invariant and descends to a square integrable function on $\Gamma\bs G.$
 We have  
\begin{equation}\lambda_{\Gamma\bs G}(\varphi)f(h)=\int_G\varphi(g)f(hg)dg=f\ast \check{\varphi}(h), \text{ where } \check{\varphi}(g):=\varphi(g^{-1}).\end{equation}\label{eq-convcheck}
Convolutions play nicely with the operators defined just above
\begin{equation}\label{eq-convrreg}\lambda_{\Gamma\bs G}(\varphi_1)\lambda_{\Gamma\bs G}(\varphi_2)f=\lambda_{\Gamma\bs G}(\varphi_1\ast \varphi_2)f, \quad \varphi_1,\varphi_2\in C_c(G).\end{equation}
In particular $\varphi\mapsto \lambda_{\Gamma\bs G}(\varphi)$ is an algebra homomorphism from the $C_c(G)$ with convolution to the algebra of bounded operators on $L^2(\Gamma\bs G).$ We also have 
\begin{equation}\label{eq-convassoc}
 (f\ast \varphi_1)\ast \varphi_2= f\ast (\varphi_1\ast \varphi_2), \quad f\in L^2(\Gamma\bs G), \varphi_1,\varphi_2\in C_c(G). 
\end{equation}
\noindent and 
\begin{equation}\label{eq-adjoint}
    \langle f_1\ast \varphi, f_2\rangle=\langle f_1, f_2\ast \check{\varphi}\rangle.
\end{equation}
Recall that $K$ was a fixed maximal compact subgroup of $G$. A function $f\colon G\to \mathbb C$ is called spherical if $f(k_1gk_2)=f(g)$ for all $k_1,k_2$. We write $L^p(G\sslash K), C(G\sslash K)$ for the spaces of $L^p$-integrable and continuous spherical function respectively. In the sequel we will heavily rely on the fact that convolution of spherical functions is commutative. 

\begin{lem}
Let $\varphi_1,\varphi_2\in L^1(G\sslash K)$ be spherical functions. Then $\varphi_1\ast \varphi_2=\varphi_2\ast \varphi_1.$
\end{lem}
\begin{proof}
It follows from \cite[Prop. 1.5.2]{gv88} combined with the definition \cite[Def. 1.5.1]{gv88}. 
\end{proof}
 For any $k$-varifold $S\subset \Gamma\bs X$ write $\tilde S$ for the lift to $X$. 
$$\tilde S=\Gamma S=\{x\in X\mid \Gamma x\in S\}.$$
Since being stationary is a local property, $S$ is stationary if and only if $\tilde S$ is.
\begin{defn} For $\varphi$ a spherical function the formula 
$$(\tilde S\ast \varphi)(h):=\int_{\tilde S} \varphi(g^{-1}h)d\mathcal H^k(gK)$$ defines a left $\Gamma$-invariant, right $K$-invariant function on $G$ which descends to a right $K$-invariant function in $L^2(\Gamma\bs G)$, denoted $S\ast \varphi$.
\end{defn}

From now on let us fix a smooth compactly supported spherical test function $\varphi\in C_c^\infty(G)$ such that $\varphi(g)\geq 0, \int_G \varphi(g)dg=1$. The constants appearing in the results below will depend on our choice of $\varphi$ but not on $S$ or $\Gamma$. 
We put \[\mathfrak u:=S\ast \varphi\in L^2(\Gamma\bs G)\text{ and }\Psi(g):=\langle \mathfrak u,\lambda_{\Gamma\bs G}(g^{-1})\mathfrak u\rangle=\langle \mathfrak u, \lambda_{\Gamma\bs G}(g)\mathfrak u\rangle.\]
Choose a Borel measurable fundamental domain $S_0\subset \tilde S$ of $S$, so that $\tilde S=\bigsqcup_{\gamma\in\Gamma}\gamma S_0.$

\begin{lem}\label{lem-BallIntForm} For any $r\geq 0$ and any  smooth compactly supported spherical function $\Phi\colon G\to\mathbb R$
$$\int_{G}\Phi(g) \Psi(g) dg= \int_{S_0}\int_G \left(\int_{\tilde S}\Phi(s^{-1}uh)d\mathcal H^k(sK)\right)(\varphi\ast\check{\varphi})(h)dh d\mathcal H^k(uK). $$
\end{lem}

\begin{proof}
By the bi-linearity of the inner product we have that the left side equals  
$$\int_{G}\Phi(g)\langle S\ast \varphi, \lambda_{\Gamma\bs G}(g^{-1})(S\ast \varphi)\rangle dg=\langle S\ast \varphi, \lambda_{\Gamma\bs G}(\check{\Phi})(S\ast \varphi) \rangle=\langle S\ast \varphi, (S\ast \varphi)\ast \Phi\rangle.$$
Let us compute the convolution on the right side of the product. Let $h\in \Gamma\bs G.$
\begin{align*}
    (S\ast \varphi)\ast \Phi (h)=&\int_G \int_{\tilde S} \varphi(s^{-1}hg^{-1})d\mathcal{H}^k(sK) \Phi(g)dg\\
    =& \int_G\int_{\tilde S} \varphi(s^{-1}hg^{-1}))\Phi(g)d\mathcal H^k(sK)\\
    =&\int_{\tilde S}(\varphi\ast \Phi)(s^{-1}h)d\mathcal H^k(sK)\\
    =&\int_{\tilde S}(\Phi\ast \varphi)(s^{-1}h)d\mathcal H^k(sK)\\
    =&\int_G\int_{\tilde S} \Phi(s^{-1}hg^{-1})\varphi(g)d\mathcal H^k(sK)dg\\
    \end{align*}
The second to last equality above is where we use the assumption that $\Phi$ is spherical, hence $\varphi\ast \Phi=\Phi\ast \varphi$. We can now evaluate the inner product in question
\begin{align*}
\langle S\ast \varphi, S\ast \varphi\ast\Phi \rangle=& \int_{\Gamma\bs G}(S\ast \varphi)(h)\int_G\int_{\tilde S} \Phi(s^{-1}hg^{-1})\varphi(g)d\mathcal H^k(sK)dgdh\\
=& \int_{\Gamma\bs G}\int_{\tilde S}\varphi(u^{-1}h)d\mathcal H^k(uK)\int_G\int_{\tilde S} \Phi(s^{-1}hg^{-1})\varphi(g)d\mathcal H^k(sK)dg dh\\
=& \int_{\Gamma\bs G}\left(\sum_{\gamma\in\Gamma}\int_{S_0}\varphi(u^{-1}\gamma h)d\mathcal H^k(uK)\right)\\ &\left(\int_G\int_{\tilde S}\Phi(s^{-1}hg^{-1})\varphi(g)d\mathcal H^k(sK)dg\right) dh\\
=& \int_{\Gamma\bs G}\left(\sum_{\gamma\in \Gamma}\int_{S_0}\int_G\int_{\tilde S}\Phi(s^{-1}\gamma hg^{-1})\varphi(g)\varphi(u^{-1}\gamma h) d\mathcal H^k(sK)dg d\mathcal H^k(uK)\right) dh\\
=& \int_{S_0}\int_G\int_G\int_{\tilde S}\Phi(s^{-1}hg^{-1})\varphi(g)\varphi(u^{-1}h)d\mathcal H^k(sK)dhdg d\mathcal H^k(uK)\\
=& \int_{S_0}\int_G\int_G\int_{\tilde S}\Phi(s^{-1}uh)\varphi(g)\varphi(hg)d\mathcal H^k(sK)dhdg d\mathcal H^k(uK)\\
=& \int_{S_0}\int_G \left(\int_{\tilde S}\Phi(s^{-1}uh)(\varphi\ast \check{\varphi})(h)d\mathcal H^k(sK)\right)dh d\mathcal H^k(uK).
\end{align*}

To go from the third-last to the second-last line, we do change of variables on $h$ first by left multiplication by $u$ and then by right multiplication by $g$.   

\end{proof}
\subsection{Monotonicity and the decay of matrix coefficients}\label{sec-stationaryproof}
From now on we restrict to the rank one case. Let $\ff$ be the function defined in (\ref{defn-ffdef}). Put
\begin{align*}
v_S(g,r):=& \int_{\tilde S} \chi_r(\ff(s^{-1}g))d\mathcal H^k(sK)\quad \text{ (same as before)},\\
w_S(g,r):=& \int_G v_S(gh)(\varphi\ast \check\varphi)(h)dh= \int_G \left(\int_{\tilde S}\chi_r(\ff(s^{-1}gh))(\varphi\ast \check{\varphi})(h)d\mathcal H^k(sK)\right)dh.
\end{align*}

Applying Lemma \ref{lem-BallIntForm} to $\Phi(g)=\chi_r(\ff(g))$, we get the identity 
\begin{equation}\label{eq-PsiIntegralformula}\int_G \chi_r(\ff(g))\Psi(g)dg=\int_{S_0}\int_G v_S(sh,r)(\varphi\ast\check{\varphi})(h)dh d\mathcal H^k(sK)=\int_{S_0}w_S(s,r)d\mathcal H^k(sK).
\end{equation}
Recall $\mathfrak u:=S\ast \varphi\in L^2(\Gamma\bs G)$. Put $\mathfrak u_0:=\mathfrak u-\frac{\mathcal H^k(S)}{\vol(\Gamma\bs G)}$. Then, since $\int_G\varphi(g)dg=1$, we have $\int_{\Gamma\bs G}\mathfrak u_0(g)dg=0.$ Put $\mathfrak u_1:=\frac{\mathcal H^k(S)}{\vol(\Gamma\bs G)}.$ Let 
$$\Psi_i(g):=\langle \mathfrak u_i,\lambda_{\Gamma\bs G}(g^{-1})\mathfrak u_i\rangle, i=0,1.$$
Then $\Psi(g)=\langle \mathfrak u,\lambda_{\Gamma\bs G}(g^{-1})\mathfrak u\rangle= \Psi_0(g)+\Psi_1(g)$. It is not hard to compute $\Psi_1(g)=\frac{\mathcal H^k(S)^2}{\vol(\Gamma\bs G)},$ so 
\begin{equation}\label{eq-PsiDecomp}
\Psi(g)=\Psi_0(g) + \frac{\mathcal H^k(S)^2}{\vol(\Gamma\bs G)}.
\end{equation}
The function $\Psi$ plays a key role in the proof of Theorem \ref{mt-stationary}. 

\begin{lem}\label{lem-matrixmono} Suppose $S\subset \Gamma\bs X$ is a stationary $k$-varifold. Then
\begin{enumerate}
\item $\Psi(1)\geq c \mathcal H^k(S)$,
\item $\int_G \chi_0(\ff(g))\Psi(g)dg\geq c \Psi(1)$,
\item  for any $r\geq 0$  $$\int_{G}\chi_r(\ff(g))\Psi(g)dg\geq c e^{\kappa(k)\|\alpha\|r}\Psi(1).$$
\end{enumerate}
\end{lem}
\begin{proof}
\begin{enumerate}
\item $$\Psi(1)=\langle S\ast\varphi,S\ast\varphi\rangle=\int_{S_0}\int_{\tilde S}(\varphi\ast\check\varphi)(s^{-1}u)d\mathcal H^k(uK)d\mathcal H^k(sK),$$
Let $\varepsilon>0$ be such that $\varphi\ast\check\varphi(g)\geq \varepsilon$ for all $g\in B_G(\varepsilon)$. By Lemma \ref{lem-smallballmono}, $\mathcal H^k(s^{-1} \tilde S\cap B(\varepsilon))>\delta=\delta(\varepsilon)>0$ for $\mathcal H^k$-almost all $s\in \tilde S$. Therefore $\Psi(1)\geq \varepsilon \delta \mathcal H^k(S).$ The constants $\varepsilon, \delta$ are independent of $S$ and $\Gamma$. 
\item 
Since $\varphi, \chi(\ff)$ are both non-negative continuous and positive at $1$, there is a constant $c>0$ such that $\varphi\ast \chi_0(\ff)\geq c \varphi.$
We have 
$$\int_G \chi_0(\ff(g))\Psi(g)dg=\langle S\ast \varphi, S\ast\varphi\ast \chi_0(\ff)\rangle\geq c\langle S\ast \varphi, S\ast\varphi\rangle=c \Psi(1).$$
\item By (\ref{eq-PsiIntegralformula}), Lemma \ref{lem-VarifoldMono} and point (2) proved just above
\begin{align*}\int_{G}\chi_r(\ff(g))\Psi(g)dg=&\int_{S_0}\int_G v_S(uhK,r)(\varphi\ast\check\varphi)(h)dh d\mathcal H^k(uK) \\
\geq& e^{\kappa(k)\|\alpha\|r}\int_{S_0}\int_G v_S(uhK,0)(\varphi\ast\check\varphi)(h)dh d\mathcal H^k(uK)\\
=& e^{\kappa(k)\|\alpha\|r}\int_G \chi_0(\ff(g))\Psi(g)dg\geq c e^{\kappa(k)\|\alpha\|r}\Psi(1).
\end{align*}
\end{enumerate}
\end{proof}
We have now enough tools to prove Theorem \ref{mt-stationary}. We reproduce the statement below. 
\begin{reptheorem}{mt-stationary}
Let $M$ be an octonionic hyperbolic $16$-manifold. There is a constant $c>0$ such that any stationary $k$-varifold $S\subset M$ satisfies $\mathcal{H}^{k}(S)\geq c \vol(M)$ for $k=14,15$.  
\end{reptheorem}
\begin{proof}[Proof of Theorem \ref{mt-stationary}]
Any finite volume octonionic manifold is of the form $M=\Gamma\bs X$ where $\Gamma$ is a lattice in $F_4^{(-20)}$ and $X=\mathbb H_\mathbb O^2$. In that case $\vol(M)=\vol(\Gamma\bs G).$
By Lemma \ref{lem-matrixmono}, for any $r\geq 0$
$$\int_G\chi_r(\ff(g))\Psi(g)dg\geq c e^{\kappa(k)\|\alpha\|r}\Psi(1).$$
On the other hand 
\begin{align*}
\int_G\chi_r(\ff(g))\Psi(g)dg=&\int_G\chi_r(\ff(g))\Psi_1(g)dg+ \int_G\chi_r(\ff(g))\Psi_0(g)dg\\
=&\frac{\mathcal H^k(S)^2}{\vol(\Gamma\bs G)}\int_G\chi_r(\ff(g))dg+ \int_G\chi_r(\ff(g))\Psi_0(g)dg.
\end{align*}
Therefore 
\begin{align}\label{eq-VolumeMonoDecay}
\frac{\mathcal H^k(S)^2}{\vol(\Gamma\bs G)}\int_G\chi_r(\ff(g))dg\geq c e^{\kappa(k)\|\alpha\|r}\Psi(1) - \int_G\chi_r(\ff(g))\Psi_0(g)dg.
\end{align}
After unfolding the definition of $\ff,\chi_r$, integrating in spherical coordinates and changing variable $t:=\alpha(H)$ we get:
\begin{align*}\int_G\chi_r(\ff(g))dg \geq& c \int_0^{(r+1)\|\alpha\|} \sinh(t)^{m_\alpha}\sinh(2t)^{m_{2\alpha}}dt\geq c e^{r(m_\alpha+2m_{2\alpha})\|\alpha\|}.
\end{align*}
The function $\Psi_0$ is a spherical matrix coefficient of the unitary representation $$L^2_0(\Gamma\bs G):=\left\{f\in L^2(\Gamma\bs G)\mid \int_{\Gamma\bs G}fdg=0\right\}.$$ The latter has no fixed vectors, so after decomposing this representation into irreducible unitary representations Lemma \ref{lem-rankonedecay} yields 
$$|\Psi_0(k_1e^Hk_2)|\leq c e^{(1-m_{2\alpha})\alpha(H)}(1+\|H\|)\Psi_0(1)\leq  c e^{(1-m_{2\alpha})\alpha(H)}(1+\|H\|)\Psi(1).$$
Using spherical integration once again, we get
\begin{align*}\left|\int_G\chi_r(\ff(g))\Psi_0(g) dg \right|
\leq& c \Psi(1)\int_0^{(r+1)\|\alpha\|} \sinh(t)^{m_\alpha}\sinh(2t)^{m_{2\alpha}}e^{(1-m_{2\alpha})t}(1+t)dt\\
\leq& c \Psi(1) e^{r(m_\alpha+m_{2\alpha}+1)\|\alpha\|}(1+r).
\end{align*}
Plugging the previous inequalities back to (\ref{eq-VolumeMonoDecay}) we get 
\begin{align*} \frac{\mathcal H^k(S)^2}{\vol(\Gamma\bs G)}e^{(m_\alpha+2m_{2\alpha})\|\alpha\|r}\geq& \Psi(1)(c_1 e^{\kappa(k)\|\alpha\|r}-c_2 e^{(m_\alpha+m_{2\alpha}+1)\|\alpha\|r}(r+1))\\
\frac{\mathcal H^k(S)^2}{\vol(\Gamma\bs G)}\geq& \Psi(1)(c_1 e^{(\kappa(k)-m_\alpha-2m_{2\alpha})\|\alpha\|r}-c_2 e^{(1-m_{2\alpha})\|\alpha\|r}(r+1))\\
\end{align*}
Put $A(r):=(c_1 e^{(\kappa(k)-m_\alpha-2m_{2\alpha})\|\alpha\|r}-c_2 e^{(1-m_{2\alpha})\|\alpha\|r}(r+1))$. Our conditions on $k$ (i.e. $k=14,15$) are chosen exactly to guarantee that $\kappa(k)-m_\alpha-2m_{2\alpha}> 1-m_{2\alpha}$ (i.e. $\kappa(k)> 1+m_\alpha+m_{2\alpha}=\dim X$). We choose big enough $r$ so that $A(r)>0$. Then, by Lemma \ref{lem-matrixmono},
\begin{align*}\frac{\mathcal H^k(S)^2}{\vol(\Gamma\bs G)}\geq & \Psi(1)A(r)\geq c \mathcal H^k(S) A(r)\\
\mathcal H^k(S)\geq&  c A(r)\vol(\Gamma\bs G).\\
\end{align*}
Theorem \ref{mt-stationary} is proved.
\end{proof}

\section{Collapsing Small Area Subsets to Lower Dimensions}\label{sec-collapsing}
This section is devoted to the proof of Theorem \ref{mt-collapse}. Let us recall the statement. 
\begin{reptheorem}{mt-collapse}
Let $\varepsilon>0$ and let $M$ be  an octonionic hyperbolic $16$-manifold with injectivity radius at least $\varepsilon.$ There exists an $\alpha=\alpha(\varepsilon)>0$ with the following property. Any rectifiable $k$-varifold $S$  with $k=14,15$ admits a homotopy $h\colon S\times [0,1]\to M$ such that
\begin{enumerate}
\item $h(s,0)=s$ for all $s\in S$,
\item $\mathcal H^{k-1}(h(S,1)))<\infty,$
\item $\mathcal H^{k+1}(h(S,[0,1])\leq c \mathcal H^k(S),$ with $c=c(\varepsilon).$
\end{enumerate}
\end{reptheorem}

Stated more plainly, Theorem \ref{mt-collapse} says that a small codimension-$2$ subset of $M$ can be efficiently deformed into a codimension-$3$ set. The fact that we can perform an actual deformation, not just fill the cycle in homology, will be of crucial importance in Section \ref{sec-branched} where we use it to prove a branched cover stability theorem for octonionic hyperbolic manifolds.

\subsection{Shortcut to an interesting corollary} \label{simplerproofofmt-collapse}

Before heading toward the proof of Theorem \ref{mt-collapse}, we sketch here a simpler proof of a consequence of Theorem \ref{mt-collapse}: every integral cycle $F$ in codimension one or two with area smaller than some constant  times the volume of the ambient space can be homotoped to a lower-dimensional set. First, we can assume that $F$ is connected and that $\mathcal H^k(F)\leq c \vol(M),$ where $c$ is the constant from Theorem \ref{mt-stationary}. 

We consider two cases according to whether $\pi_1(F)$ includes to a finite index subgroup of $\pi_1(M)$ or not. 

In the first case, we pass to a finite cover $M'$ such that a connected component $F'$ of the lift of $F$ to $M'$ has $\pi_1(F')$ surjecting onto $\pi_1(M')$. Let $F''$ be an integral cycle homologous to $F'$ mod $2$, with the minimal $k$-dimensional volume. Such representatives exists by \cite[Theorem 1.4]{Almgren1968ExistenceAR} and are stationary. By Theorem \ref{mt-stationary}, this cycle must be trivial, since there are no stationary $k$-varifolds in $M'$ of $k$-volume less than $c\vol(M')$. We deduce that $F$ is null-homologous mod $2$, and conclude by \cite{white1984mappings} that it can be homotoped to sets of arbitrary low $k$-dimensional volume. Then an application of Federer-Fleming deformation (see Section \ref{sec-FF}) homotopes it into a $(k-1)$-dimensional set. 



In the second case, we can lift $F$ to an infinite sheeted cover $\Lambda\bs X,$ where $\Lambda$ is the inclusion of $\pi_1(F)$ to $\pi_1(M)$. Then \cite[Theorem 1.1]{cmw23} (see also section \ref{section:naturalvectorfields}) implies that $F$ can be homotoped to a set of arbitrarily small $k$-volume $F'$. Applying the Federer-Fleming deformation to $F'$ (see Section \ref{sec-Fed}) and projecting down to $M$ gives the desired homotopy of $F$ to a lower dimensional set.

\subsection{Sketch of the argument} 

Let $X= \mathbb H_\mathbb O^2$ and $G=F_4^{(-20)}$. Any compact octonionic hyperbolic manifold $M$ is of the form $\Gamma\bs X$, where $\Gamma$ is a cocompact lattice in $G$.  

We will construct the homotopy one step at a time, arguing that we can always decrease the $k$-dimensional volume a little bit, provided that the initial $k$-dimensional volume  (or just $k$-volume) of $S$ is not too big. We will either apply a suitable flow squeezing $k$-volume or the Federer-Fleming deformation (see Section \ref{sec-FF}). 
The rough idea to construct the flow is to run the proof of monotonicity estimate from Lemma \ref{lem-abstmonotonicity} at each point of $S$ and obtain the following dichotomy. Either the monotonicity estimate holds at a point $sK\in S$ or one of the vector fields that appears in the proof of Lemma \ref{lem-abstmonotonicity} infinitesimally squeezes the $k$-dimensional volume of $S$. In the proof, the points admitting such a field will be called heavy points. 
Such fields will induce the flow deforming $S$ to lower $k$-dimensional volume. To exclude the scenario where the monotonicity estimate holds at all points of $S$, we argue as in the proof of Theorem \ref{mt-stationary}. The proof would be literally the same if we knew that the analogue of Lemma \ref{lem-smallballmono} held for $S$ or even for a fixed positive proportion of $S$. 
Unfortunately, there is no reason to assume $S$ satisfies a uniform monotonicity bound at small scales. For example $S$ could be a very thin tube around a $(k-1)$ or lower dimensional manifold. At this point we are saved by the Federer-Fleming deformation. We prove that if Lemma \ref{lem-smallballmono} fails at most points of $S$, then an application of Federer-Fleming deformation significantly decreases the $k$-dimensional volume. By repeatedly applying the squeezing flows or the Federer-Fleming deformations we make the $k$-dimensional volume of $S$ small enough that a single application of Federer-Fleming collapses it into a $(k-1)$-dimensional set.  

This sketch omits several important points. For example what if the flows that infinitesimally decrease the $k$-dimensional volume only do it for a very short time? We address this issue and others using the second variation formula and fine properties of the squeezing fields. The proof is constructive, and gives an explicit algorithm to deform a small $k$-varifold  to a $(k-1)$-dimensional set.

\subsection{Thin part of $S$ and the Federer-Fleming deformation}

We keep the notation from Section \ref{sec-VarifoldsMonoMC}. In particular, $\varphi$ will be a fixed smooth compactly supported spherical function with $\int_G\varphi(g)dg=1$. For technical reasons we choose $\varphi$ so that $\|\varphi\|_2^2=(\varphi\ast\check{\varphi})(1)>1$ and $\supp \varphi\subset B_G(1)$. The functions $\mathfrak u:=S\ast \varphi\in L^2(\Gamma\bs G)$ are defined as before and $\Psi(g):=\langle \mathfrak u,\lambda_{\Gamma\bs G}(g^{-1})\mathfrak u\rangle.$ In this section the rectifiable varifold $S$ will evolve and the same will hold for $\Psi$ and other functions implicitly dependent on $S$. Still, we shall work with at most one $S$ at a time so hopefully it won't cause any confusion. We write $\tilde S$ for the lift of $S$ to $X$ and $S_0\subset \tilde S$ for a fundamental domain. The latter can be chosen unambiguously by intersecting $\tilde S$ with a fixed fundamental domain of $\Gamma\bs X$ in $X$.

Using Theorem \ref{thm-UTboundedgeo}, we choose $\delta>0$ such that all compact quotients $\Gamma\bs X$ with injectivity radius at least $\varepsilon$ admit an $(r,\delta)$-uniform triangulation for any small enough $r$. We choose $r_0>0$ such that $\varphi\ast\check\varphi(h)\geq 1$ for all $h\in B_G(3r_0)$ and let $\mathcal T$ be an $(r_0,\delta)$-uniform triangulation on $\Gamma\bs X$. Write $\FF$ for the Federer-Fleming deformation relative to $\mathcal T$ and $S$. The constants controlling the deformation in Lemmas \ref{lem-FFproperties}, \ref{lem-FFvanishing} depend only on $r_0,\delta$, so they are independent of $\Gamma$ and $S$. We will write $M=\Gamma\bs X$ for the locally symmetric manifold and $M^{(m)}, M^{[m]}$ for the $m$-skeleton and the set of $m$-simplices respectively.

One of the major difficulties in the proof of Theorem \ref{mt-collapse} compared to Theorem \ref{mt-stationary} is the lack of the monotonocity estimate from Lemma \ref{lem-smallballmono}.
Recall from Section \ref{sec-stationaryproof} that $w_S(\Gamma gK, r):=\int_G v_S(ghK,r)(\varphi\ast\check\varphi)(h)dh$ and $v_S(\Gamma g K,r):= \int_{\tilde S} \chi_r(\ff(s^{-1}g))d\mathcal H^k(sK)$. To lighten notation, especially inside the integrals, we will sometimes write $gK$ for points in $\Gamma\bs X=\Gamma\bs G/K$ with the understanding that we mean the coset $\Gamma g K$.
For $\eta>0$ define 
$$S_{<\eta}:=\{\Gamma sK\in S\mid w_S(sK,0)<\eta\}, \quad S_{\geq \eta}=S\setminus S_{<\eta}.$$

\begin{lem}\label{lem-ThinFFdichotomy}
There are constants $\eta=\eta(\varepsilon)>0, c=c(\varepsilon)>0$ such that either 
\begin{enumerate}
\item $\mathcal H^k(S_{\geq\eta})\geq c \mathcal H^k(S),$ or 
\item $\mathcal H^k({\rm FF}(S))\leq \frac{1}{2}\mathcal H^k(S).$
\end{enumerate}
\end{lem}
\begin{proof}
By Lemma \ref{lem-FFvanishing}, there exists an $\eta>0$ dependent only on $\varepsilon$ (because $\delta, r_0$ are fixed depending only on $\varepsilon$), such that any point $\Gamma sK\in S$ with $\mathcal H^k(B(\Gamma sK,2r_0)\cap S)\leq \eta $ is mapped to $M^{(k-1)}$ under the Federer-Fleming map. Since $(\varphi\ast \check\varphi)(h)\geq 1$ for all $h\in B(3r_0)$, we conclude that $r_0$-neighbourhood of any point $\Gamma sK\in S$ with $w_S(sK,0)<\eta$ must be mapped to $M^{(k-1)}.$ Let $Z_0$ be the union of all simplices of $\mathcal T$ intersecting $S_{<\eta}$ and let $Z_1$ be the union of all simplices disjoint from $S_{<\eta}$ (note that these two sets might intersect because a simplex in the first set and a simplex in the second set might share a common face). 
By our previous observations, together with the fact that the diameters of all simplices in $\mathcal T$ are bounded by $r_0$, we get 
$${\rm FF}(S\cap Z_0)\subset M^{(k-1)}, \text{ so } \mathcal H^k({\rm FF}(S \cap Z_0))=0.$$
On the other hand, by Lemma \ref{lem-FFproperties}, $\mathcal H^k({\rm FF}(S\cap Z_1))\leq c_0 \mathcal H^k(S\cap Z_1)$, with $c_0=c_0(\varepsilon).$ Together we get 
$$\mathcal H^k({\rm FF}(S))\leq c_0 \mathcal H^{k}(S\cap Z_1)\leq c_0 \mathcal H^k(S_{\geq \eta}).$$
The Lemma now follows. If $\mathcal H^k(S_{\geq \eta})\leq \frac{1}{2c_0}\mathcal H^k(S)$, then the second case occurs. If not, then we are in the first case, with $c=\frac{1}{2c_0}.$ 
\end{proof}
\begin{cor}\label{cor-FFalternative} There are constants $\eta=\eta(\varepsilon)>0, c=c(\varepsilon)>0$ such that either 
\begin{enumerate}
\item $\int_{S_{\geq\eta}}w_S(sK,0)d\mathcal H^k(sK) \geq c \Psi(1)$ and $\Psi(1)\geq \eta c \mathcal H^k(S)$ or 
\item $\mathcal H^k({\rm FF}(S))\leq \frac{1}{2}\mathcal H^k(S).$
\end{enumerate}
\end{cor}
\begin{proof} By Lemma \ref{lem-ThinFFdichotomy}, we can assume $\mathcal H^k(S_{\geq\eta})\geq c_0 \mathcal H^k(S),$ with $c_0=c_0(\varepsilon)$. By making $c_0$ smaller if necessary we can also assume $c_0<1$.
Then 
\begin{align*}\Psi(1)=&\int_{S}w_S(sK,0)d\mathcal H^k(sK)\leq \int_{S_{\geq \eta}}w_S(sK,0)d\mathcal H^k(sK)+ \eta \mathcal H^k(S_{<\eta})\\
\leq& \int_{S_{\geq \eta}}w_S(sK,0)d\mathcal H^k(sK)+ \eta(1-c_0)\mathcal H^k(S). 
\end{align*}
But $\int_{S_{\geq \eta}}w_S(sK,0)d\mathcal H^k(sK)\geq \eta c_0 \mathcal H^k(S)$, which implies 
\[
\eta (1-c_0) \mathcal H^k(S) \leq \eta (1-c_0) \frac{1}{\eta c_0} \int_{S_{\geq\eta}}w_S(sK,0)d\mathcal H^k(sK)
\]
so $$\int_{S_{\geq\eta}}w_S(sK,0)d\mathcal H^k(sK) \geq c_0  \Psi(1).$$ 
Similarly 
\begin{align*}\Psi(1)=&\int_{S}w_S(sK,0)d\mathcal H^k(sK)\geq  \int_{S_{\geq \eta}}w_S(sK,0)d\mathcal H^k(sK)\\ \geq& \eta \mathcal H^k(S_{\geq \eta}) 
\geq \eta c_0\mathcal H^k(S). 
\end{align*}
\end{proof}
\subsection{Heavy points and their properties}
In this section we assume that $\kappa(k)>\dim X$, which happens exactly for the range of $k$ allowed in Theorem \ref{mt-collapse}. We fix an auxiliary constant 
$$\kappa':=\frac{\kappa(k)+\dim X}{2}.$$ The only important thing about $\kappa'$ is that it is a fixed constant strictly between $\kappa(k)$ and $\dim X=m_\alpha+m_{2\alpha}+1.$ 
\begin{defn}
A point $\Gamma K\in S$ will be called \emph{$t$-heavy} if the two conditions are satisfied. 
\begin{enumerate}
\item $\Gamma sK\in S_{\geq \eta},$
\item $\frac{\partial}{\partial t}w_S(sK,t)\leq \kappa'\|\alpha\| w_S(sK,t)$ and $w_S(sK,t+8)\leq e^{8\kappa'\|\alpha\|}w_S(s K,t)$. 
\end{enumerate}
The second condition can be interpreted as a bound on the ordinary and discrete derivatives of $\log w_S(sK,t)$ at time $t.$ The number $8$ is not accidental, and is related to radii of the support of $\varphi$ and the function $\chi_r(\ff)$. 
\end{defn}
The name heavy is motivated by the following geometric interpretation. At the $t$-heavy points $\Gamma sK\in S$, the varifold $S$ has suspiciously large mass in a $t$-ball around $\Gamma sK$ compared to $(t+8)$-ball. This property will be used to construct flows squeezing the $k$-volume of $S$. 
We now adapt the proof of Theorem \ref{mt-stationary} from section \ref{sec-stationaryproof} to show the existence of heavy points. 
\begin{lem}\label{lem-heavypintdichotomy}
There are constants $\alpha, r_1>0$, dependent only on $\varepsilon,$ with the following property. Any rectifiable $k$-varifold $S\subset \Gamma\bs X$ with $\injrad(\Gamma\bs X)\geq \varepsilon$ and $\mathcal H^k(S)\leq \alpha \vol(\Gamma\bs X)$ satisfies either 
\begin{enumerate}
\item $\mathcal H^k({\rm FF}(S))\leq \frac{1}{2} \mathcal H^k(S),$
\item there is a $t$-heavy point $\Gamma sK\in S, t\leq r_1.$
\end{enumerate}
\end{lem}
\begin{proof}
By Corollary \ref{cor-FFalternative}, either (1) holds or $\int_{S_{\geq\eta}}w_S(sK,0)d\mathcal H^k(sK) \geq c_0 \Psi(1)$ and $\Psi(1)\geq \eta c_0\mathcal H^k(S)$. Assume we are in the second case.
For the sake of contradiction, suppose none of the points  $S_{\geq\eta}$ is $t$-heavy for $0\leq t\leq r_1$. This means that for $\Gamma sK\in S_{\geq \eta}$ and any $0<t<r_1$ we have either $\frac{\partial}{\partial t}\log w_S(sK,t) > \kappa'\|\alpha\|$ or $\log w_S(sK,t+8)-\log w_S(sK,t) > 8\kappa'\|\alpha\|$. Since $\log w_S(sK,t)$ is non-decreasing, it easily implies $\log w_S(sK,r_1)\geq \log w_S(sK,0)+\kappa'\|\alpha\| (r_1-8).$ By (\ref{eq-PsiIntegralformula}), we get 
\begin{align*}\int_G \chi_r(\ff(g))\Psi(g)dg \geq& \int_{S_{\geq\eta}}w_S(sK,r)d\mathcal H^k(sK)\\ \geq& c_1 e^{\kappa'\|\alpha\|r_1}\int_{S_{\geq\eta}}w_S(sK,0)d\mathcal H^k(sK) \geq c_1 e^{\kappa'\|\alpha\|r_1}\Psi(1).
\end{align*}
Now arguing exactly as in the proof of Theorem \ref{mt-stationary} we arrive at
$$\frac{\mathcal H^k(S)^2}{\vol(\Gamma\bs G)}\geq \Psi(1)(c_2e^{(\kappa'-m_\alpha-2m_{2\alpha})\|\alpha\|r_1}-c_3e^{(1-m_{2\alpha})\|\alpha\|r_1}(1+r_1)),$$
with $c_1,c_2$ depending only on $\varepsilon$. We put $$A(r_1):=c_2e^{(\kappa'-m_\alpha-2m_{2\alpha})\|\alpha\|r_1}-c_3e^{(1-m_{2\alpha})\|\alpha\|r_1}(1+r_1).$$ Since $\kappa'>\dim X=1+m_\alpha+m_{2\alpha}$, for $r_1$ big enough $A(r_1)>0$ (how big depending only on $\varepsilon$). Then, 
$$\frac{\mathcal H^k(S)^2}{\vol(\Gamma\bs G)}\geq \Psi(1)A(r_1)\geq \eta c_0 A(r_1)\mathcal H^k(S).$$
If we took $\alpha< \eta c_0 A(r_1)$ this would be a contradiction with our assumption. Therefore, for the right choice of $\alpha, r_1$ the dichotomy holds. 
\end{proof}
In the next lemma, we unwind the proof of Lemma \ref{lem-abstmonotonicity} to show that certain vector fields ``collapsing" the space towards a $t$-heavy point will squeeze the $k$-volume of $S$. 
\begin{lem}\label{lem-heavyfield}
Let $S\subset \Gamma\bs X$ be a rectifiable $k$-varifold and $\injrad(\Gamma\bs X)\geq \varepsilon.$ Suppose $S$ contains a $t$-heavy point $\Gamma sK$. Then, there is a smooth vector field $Y$ on $\Gamma\bs X$ with the associated flow $\Phi_v\colon (\Gamma\bs X)\to (\Gamma\bs X), \frac{\partial}{\partial v}\Phi_v(x)=Y(\Phi_v(x)), x\in \Gamma\bs X$ satisfying 
\begin{enumerate}
\item $\|Y\|,\|\nabla Y\|\leq c=c(\varepsilon)$ uniformly on $\Gamma\bs X$.
\item $Y$ is supported on $\Gamma B(sK, t+3)$ (seen as a subset of $\Gamma\bs X$).
\item $$\frac{\partial}{\partial v}\mathcal H^k(\Phi_v(S))|_{v=0}\leq -\left(\frac{\kappa(k)-\dim X}{2}\right)w_S(sK,t).$$
\end{enumerate}
\end{lem}
\begin{proof}
We start by defining fields $Y^t_{gK}, gK\in X, t>0$. We'll write $l(g) F(h):=F(g^{-1}h)$ for any function $F$ on $X$. The action naturally extends to vector fields on $X$, in such a way that $l(g)F_1\grad(l(g)F_2)=l(g)(F_1\grad F_2)$ for any $F_1,F_2\in C^1(X)$. 
Recall the definition of the function $\ff$ that smoothly approximates the distance function (\ref{defn-ffdef}). Put 
$$Y^t_{gK}:=-l(g)\left(\chi_t(\ff)\grad \ff\right).$$
This is a smooth vector field on $X$ supported on $B(gK,t+1)$, satisfying $\|Y_{gK}^t\|,$ $ \|\nabla Y_{gK}^t\|\leq c_0$ uniformly on $X$. One can imagine that this is a field pulling the interior of the $(t+1)$-ball around $gK\in X$ towards the center at roughly constant unit speed, suitably smoothed at $gX$ and the boundary of the ball.   
Let $\Gamma sK\in S$ be the $t$-heavy point. Define 
$$Y:=\sum_{\gamma\in\Gamma}\int_{G}(\varphi\ast\check\varphi)(h) Y^t_{\gamma s hK}dh.$$
The field is a left $\Gamma$ invariant smooth vector field on $X$ so it descends to a smooth vector field on $\Gamma\bs X$. For each $\gamma\in\Gamma$ the field $\int_{G}(\varphi\ast\check\varphi)(h) Y^t_{\gamma s h}dh$ is supported on $B(\gamma sK, t+3)$, since $\supp (\varphi\ast\check\varphi)\subset B(2).$ For each $\gamma\in \Gamma$, the field $\int_G (\varphi\ast \check\varphi)(h)Y^t_{\gamma shK}$ is a convex combination of the fields of type $Y^t_{gK}$ so its norm and the norms of its first two covariant derivatives are still bounded by $c_0.$ The sum over $\gamma$ is locally finite, with the number of non-vanishing summands over a point bounded uniformly in terms of $\varepsilon.$ It follows that the field $Y$ satisfies the conditions (1) and (2). It remains to verify (3). 
\begin{align*}\frac{\partial}{\partial v}\mathcal H^k(\Phi_v(S))|_{v=0}=&\int_S{\rm div}_S(Y)(uK)d\mathcal H^k(uK)\\
=& \sum_{\gamma\in\Gamma}\int_{S_0}\int_G (\varphi\ast\check\varphi)(h){\rm div}_S(Y^t_{\gamma sh})(uK) dh d\mathcal H^k(uK)\\
=& \int_{\tilde S}\int_G (\varphi\ast\check\varphi)(h){\rm div}_S(Y^t_{sh})(uK) dh d\mathcal H^k(uK)\\
\end{align*}
We now use the formula for divergence ${\rm div}_S Y_{sh}^t$ from Lemma \ref{lem-abstmonotonicity}, and repeat the computation therein to upper bound the last expression. 
\begin{align*}
\int_{\tilde S}\int_G (\varphi\ast\check\varphi)(h)&{\rm div}_S(Y^t_{sh}) dh d\mathcal H^k(uK)\\ \leq& -\kappa(k)\|\alpha\| \int_{\tilde S}\int_G \chi_t(\ff(u^{-1}sh))(\varphi\ast\check\varphi)(h) dh d\mathcal H^k(uK)\\ +& \int_{\tilde S}\int_G \frac{\partial}{\partial t} \chi_t(\ff(u^{-1}sh))(\varphi\ast\check\varphi)(h)dh d\mathcal H^k(uK)\\
=& -\kappa(k)\|\alpha\|w_S(sK,t)+\frac{\partial}{\partial t}w_S(sK,t)\\
\leq& -\left(\frac{\kappa(k)-\dim X}{2}\right)\|\alpha\|w_S(sK,t).
\end{align*} The last inequality used the assumption that $\Gamma sK$ was $t$-heavy: $\frac{\partial}{\partial t}w_S(sK,t)\leq \kappa'\|\alpha\| w_S(sK,t) =\frac{(\text{dim}(X) + \kappa)}{2} \|\alpha\| w_S(sK,t) .$
\end{proof}
\begin{lem}\label{lem-secondvariation}
Let $Y$ be a smooth vector field on $\Gamma\bs X$, supported on a closed set $B$, satisfying $\|Y\|, \|\nabla Y\|, \| \nabla^2 Y \| \leq A$ uniformly on $\Gamma\bs X$. Let $\Phi_v$ be the flow associated to $Y$. For any rectifiable $k$-varifold $S\subset \Gamma\bs X$ 
$$ \left|\frac{\partial^2}{(\partial v)^2}\mathcal H^k(\Phi_v(S))\right|\leq c \mathcal H^k(\Phi_v(S)\cap B).$$
\end{lem}
\begin{proof} 
Let $e_1,\ldots, e_k$ be $C^1$-vector fields forming an orthonormal basis of $T_{\Gamma sK}S$ at $\mathcal H^k(S)$-almost every point of $S$. Write $R$ for the curvature tensor of $\Gamma\bs X$. By the second variation formula (see \cite[(1.142)]{colding2011course} or \cite[Theorem 6.1]{montezumanotes})
\begin{align*}\frac{\partial^2}{(\partial v)^2}\mathcal H^k(\Phi_v(S))|_{v=0}=&\int_S \left( \sum_{i=1}^k(\|\nabla_{e_i}^{\perp} Y\|^2 - {\rm div}_S \nabla_Y Y-\sum_{i=1}^k R(Y,e_i,Y,e_i) \right. \\ &\left.+ \left(\sum_{i=1}^k \langle \nabla_{e_i}Y,e_i\rangle\right)^2 - \sum_{i,j=1}^k \langle \nabla_{e_j}Y,e_i\rangle \langle \nabla_{e_i}Y,e_j\rangle \right)d\mathcal H^k.\\
\left|\frac{\partial^2}{(\partial v)^2}\mathcal H^k(\Phi_v(S))|_{v=0}\right| \leq &\int_{S\cap B} \left( kA^2 + kA^2 +c_0 kA^2 + k^2A^2 + k^2A^2 \right)d\mathcal H^k\\
\leq& c \mathcal H^k(S\cap B).\end{align*} Same computation yields the bound at all times $v$. 
\end{proof}

\begin{lem}\label{lem-flowvolumebound}
Let $Y$ be a smooth vector field on $\Gamma\bs X$, supported on a closed set $B$, satisfying $\|Y\|\leq A$ uniformly on $\Gamma\bs X$. Let $\Phi_v$ be the flow associated to $Y$. For any rectifiable $k$-varifold $S\subset \Gamma\bs X$ 
and any $v_0>0$ 
$$\mathcal H^{k+1}\left(\bigcup_{v\in [0,v_0]}\Phi_v(S)\right)\leq A\int_0^{v_0}\mathcal H^k(\Phi_v(S)\cap B)dv.$$
\end{lem}
\begin{proof}

For each fixed $v$, $\Phi_v$ restricted to $S$ is a diffeomorphism on the subset of regular points. Denote by $g_v$ the induced metric on $\Phi_v(S)$.  We can define a metric on $S \times [0,v_0]$ by $\Phi_v(s)^* g_v + A dv^2$ at each point $(s,v)$. The volume of $S \times [0,v_0]$ in this metric is equal to $A\int_0^{v_0}\mathcal H^k(\Phi_v(S)\cap B)dv.$     Then since the Jacobian of the  map $\Phi\colon S \times [0,v_0] \rightarrow \Gamma \backslash X $ is pointwise at most one, this proves the lemma.

\end{proof}
\begin{lem}\label{lem-flowsqueezing}
Let $S\subset \Gamma\bs X$ be a rectifiable $k$-varifold and $\injrad(\Gamma\bs X)\geq \varepsilon.$ Suppose $S$ contains a $t$-heavy point $sK$. Let $Y$ be the vector field afforded by Lemma \ref{lem-heavyfield} and let $\Phi_v$ be the associated flow on $\Gamma\bs X$. Then there is a time $v_0=v_0(\varepsilon)>0$ and 
$\tau_0=\tau_0(\varepsilon)>0$ such that 
\begin{enumerate}
\item $\mathcal H^k(\Phi_{v_0}(S))\leq \mathcal H^k(S)-\tau_0,$
\item $\mathcal H^{k+1}(\bigcup_{v\in [0,v_0]}\Phi_v(S))\leq c (\mathcal H^k(S)-\mathcal H^k(\Phi_{v_0}(S))),$ with $c=c(\varepsilon)>0.$
\end{enumerate}
\end{lem}
\begin{proof} By Lemma \ref{lem-heavyfield}, $\frac{\partial}{\partial v}\mathcal H^k(\Phi_v(S))|_{v=0}\leq -\left(\frac{\kappa(k)-\dim X}{2}\right)\|\alpha\|w_S(sK,t).$ Let $B:=\Gamma B(sK,t+3)\subset \Gamma\bs X$. 
Put $v_1=\sup\{v\geq 0\mid \frac{\partial}{\partial v}\mathcal H^k(\Phi_v(S))< 0\}.$ By definition, $\mathcal H^k(\Phi_v(S))$ will be decreasing for $v\in [0,v_1]$. We know that $Y$ is supported on $B$ so the flow fixes the outside of $B$. In particular, $\mathcal H^k(\Phi_v(S)\cap B)$ is also decreasing on $[0,v_1]$. 

By Lemma \ref{lem-secondvariation}, for all $v\in[0,v_1]$ 
$$\left|\frac{\partial^2}{(\partial v)^2}\mathcal H^k(\Phi_v(S))\right|\leq c_0 \mathcal H^k(\Phi_v(S)\cap B)\leq c_0 \mathcal H^k(S\cap B).$$
By the definition of a $t$-heavy point, $w_S(sK, t+8)\leq e^{8\kappa'\|\alpha\|}w_S(sK,t)$. If we unfold the definition of $w_S(sK,t+8)$, we see that 
\begin{align*}w_S(sK,t+8)=&\int_{\tilde S}\int_G \chi_{t+8}(\ff(u^{-1}sh))(\varphi\ast\check\varphi)(h)dh d\mathcal H^k(S)\\ \geq& \int_{\tilde S}(\chi_{t+8}(\ff)\ast \varphi\ast\check\varphi)(u^{-1}s)d\mathcal H^k(uK).
\end{align*}
The function $\varphi\ast\check\varphi$ is supported on $B(2)$ and has total mass $1$, while $\chi_{t+8}(\ff)$ is equal to $1$ on $B(t+6)$. The convolution $\chi_{t+8}(\ff)\ast \varphi\ast\check\varphi$ is therefore equal $1$ on $B(t+4)$, so 
\begin{equation}\label{eq-ScapBmono}\mathcal H^k(S\cap B)\leq \mathcal H^k(\tilde S\cap B(sK,t+3))\leq w_S(sK,t+8)\leq c_1 w_S(sK,t), \text{ with } c_1=e^{8\kappa'\|\alpha\|}.\end{equation}
It follows that 
$$\left|\frac{\partial^2}{(\partial v)^2}\mathcal H^k(\Phi_v(S))\right|\leq c_0 c_1 w_S(sK,t), \text{ for all } v\in [0,v_1].$$
We can now deduce that 
\begin{align*}\frac{\partial}{\partial v}\mathcal H^k(\Phi_v(S))\leq& \left(c_0c_1v -\left(\frac{\kappa(k)-\dim X}{2}\right)\|\alpha\|\right)w_S(sK,t) \\
\end{align*}
for all $v\in [0,v_1]$. Taking $v_0:=\frac{(\kappa(k)-\dim X)\|\alpha\|}{4c_0c_1}$ we get 
\begin{equation} \label{eqn-firstassertion}\mathcal H^k(\Phi_{v_0}(S)))\leq \mathcal H^k(S)-\left(\frac{\kappa(k)-\dim X}{4}\right)\|\alpha\| w_S(sK,t)  v_0.\end{equation}
Since $w_S(sK,t)\geq \eta$, this proves the first part of the lemma with $\tau_0=\frac{(\kappa(k)-\dim X)\|\alpha\|\eta v_0}{4}.$ To prove the second assertion, we apply Lemma \ref{lem-flowvolumebound} and observe that 
\begin{align*}
\mathcal H^{k+1}(\{\Phi_v(S)\mid v\in [0,v_0]\})\leq & A\int_0^{v_0}\mathcal H^k(\Phi_v(S)\cap B)dv\leq Av_0 \mathcal H^k(S\cap B)\\
\leq & Av_0 c_1 w_S(sK,t) \leq c_2\left(\mathcal H^k(S) - \mathcal H^k(\Phi_{v_0}(S))\right),
\end{align*} by (\ref{eq-ScapBmono}), (\ref{eqn-firstassertion}), and the fact  that $\mathcal H^k(\Phi_v(S)\cap B)$ is decreasing on $[0,v_1].$
\end{proof}
bo\subsection{Proof of Theorem \ref{mt-collapse}}
\begin{proof} Let $M=\Gamma\bs X$ be a compact octonionic hyperbolic manifold of injectivity radius at least $\varepsilon$. We fix an $(r_0,\delta)$-uniform triangulation on $M$, with $r_0,\delta>0$ depending only on $\varepsilon$. Denote the initial $k$-varifold by $S_0:=S$, with $\mathcal H^k(S)<\alpha \vol(\Gamma\bs X)=\alpha\vol(M),$ where $\alpha$ is as in Lemma \ref{lem-heavypintdichotomy}. We will construct a sequence of varifolds $S_i,i\in\mathbb N$ together with homotopies $h_i\colon S_i\times [0,1]\to \Gamma\bs X$ between $S_i$ and $S_{i+1}$, such that 
\begin{align*}
\mathcal H^k(S_{i+1})\leq& \max\left\{\frac{1}{2}\mathcal H^k(S_i), \mathcal H^k(S_i)-\tau_0\right\},\\
\mathcal H^{k+1}(h_i(S_i,[0,1]))\leq& c (\mathcal H^k(S_i)-\mathcal H^k(S_{i+1})).
\end{align*}
The construction is inductive. Suppose we have constructed $S_i$. Then $\mathcal H^k(S_i)<\alpha \vol(\Gamma\bs X)$ so by Lemma \ref{lem-heavypintdichotomy}, either $S_i$ contains a $t$-heavy point, for $t\leq r_1=r_1(\varepsilon)$, or the Federer-Fleming deformation decreases the $k$-volume by at least a half. In the second case we put $S_{i+1}:={\rm FF}(S_i)$ and let $h_i$ be the homotopy given by the Lemma \ref{lem-FFproperties}.

In the first case, we use the $t$-heavy point to construct the vector field $Y$ from Lemma \ref{lem-heavyfield} and the flow $\Phi_v$ associated to $Y$. By Lemma \ref{lem-flowsqueezing}, there is a time $v_0>0$ such that $\mathcal H^k(\Phi_{v_0}(S_i))\leq \mathcal H^k(S_i)-\tau_0.$ We put 
$$h_i(s,t):=\Phi_{tv_0}(s),\quad s\in S_i, t\in [0,1].$$
By Lemma \ref{lem-flowsqueezing}, $\mathcal H^{k+1}(h_i(S_i,[0,1]))\leq c(\mathcal H^k(S_i)-\mathcal H^k(S_{i+1})).$
It is clear that the sequence of volumes $\mathcal H^k(S_i)$ tends to $0$, so at some point $\mathcal H^k(S_{i_0})\leq \eta,$ where $\eta$ is the constant from Lemma \ref{lem-FFvanishing}. Since there can be no more heavy points, the next step will be the Federer-Fleming deformation. By Lemma \ref{lem-FFvanishing}, $S_{i_0+1}={\rm FF}(S_{i_0})\subset M^{(k-1)}$, so $S_{i_0+1}$ has finite $(k-1)$-dimensional volume. 
We obtain the desired homotopy $h$ from $S=S_0$ to $S_{i_0+1}$ by gluing together the homotopies $h_i, i=0,\ldots, i_0.$ We have 
$$\mathcal H^{k+1}(h(S,[0,1]))\leq \sum_{j=0}^{i_0} \mathcal H^{k+1}(h_i(S_i,[0,1]))\leq  \sum_{j=0}^{i_0} c (\mathcal H^k(S_i)-\mathcal H^k(S_{i+1}))\leq c \mathcal H^k(S).$$
\end{proof}
\section{Stratified Fibrations}\label{sec-stratified}

The following concept will be useful for us \cite[pg. 762]{gr09}.  

\begin{defn} \label{sf}
	Let $X$ and $Y$ be locally compact simplicial complexes.  We say that $f\colon X \rightarrow Y$ is a \textit{stratified fibration} if it satisfies the following. Let $Y^{(i)}$, $i=0,..,d$, be the $i$-skeleta of $Y$.  
	\begin{enumerate} 
		\item $f$ restricted to the preimage of each open face of $Y^{(i)} - Y^{(i-1)}$ defines a fiber bundle. 
		\item There are homotopies $h_i(t): Y^{(i)} \rightarrow Y^{(i)} $ of the identity map, $t\in [0,1]$, so that $h_i$ maps some neighborhood $U^{i-1}$ of $Y^{(i-1)}$ in $Y^{(i)}$ to $Y^{(i-1)}$.  The restriction of the homotopy $h_i(t)$ to $U^{i-1}$ moreover lifts to a homotopy of the identity map of  $f^{-1}(U^{i-1})$: i.e., there is a homotopy $h_i^{\text{lift}}: f^{-1}(U^{i-1}) \times I \rightarrow X$ such that $h_i^{\text{lift}}(\cdot, 0)= \text{Id}$ and $f(h_i^{\text{lift}}(y,t)) = h_i (f(y),t)$.  
	\end{enumerate}
\end{defn} 


Let $M$ be a smooth manifold that has been given a triangulation $\mathcal{T}$. 
A triangulation of $M$ is a finite simplicial complex $\mathcal{T}$ together with a homeomorphism $\mathcal{T} \rightarrow M$ that restricts to a smooth embedding on the interior of each simplex in $\mathcal{T}$. We will treat simplices in $\mathcal{T}$ as interchangeable with their embedded images in $M$.  Every smooth manifold has a triangulation \cite[Theorem 10.6]{munkres2016elementary}. As in previous sections we write $M^{(m)}, M^{[m]}$ for the $m$-skeleta and sets of $m$-simplices respectively.   

In the interests of convenience and keeping the paper maximally self-contained, we adopt the following definition of piecewise-linear (PL) map, which differs slightly from the standard definition \cite{rspltopology}. Every map that is PL in our sense is also PL according to the standard definition.  


\begin{defn} \label{PLdefn}
	We say that a map $F$ from a triangulated manifold $(M,\mathcal{T})$  to $\mathbb{R}^d$ is piecewise linear (PL) if there is a map $\Phi:M \rightarrow \mathbb{R}^n$ so that the following holds.  First we require that $\Phi$ maps each simplex $\sigma$ in $\mathcal{T}$ homeomorphically onto a simplex $\sigma'$ in $\mathbb{R}^n$. Second we require that, for every simplex $\sigma$ in $(M,\mathcal{T})$, that the map $F \circ (\Phi|_{\sigma})^{-1}: \Phi(\sigma) \rightarrow \mathbb{R}^d $  is the restriction of an affine linear map $\mathbb{R}^n \rightarrow \mathbb{R}^d$.  
\end{defn}
We note that it is possible to construct many PL maps $M \rightarrow \mathbb{R}^d$ by fixing a triangulation $\mathcal{T}$ of $M$, choosing some map $\Phi$ as above, choosing where the vertices of $\mathcal{T}$ (i.e. $M^{(0)}$) map to, and then extending linearly over each simplex. We also note that it is possible to construct many such maps $\Phi$ by choosing where in $\mathbb{R}^n$ to send the vertices of $\mathcal{T}$ and then extending inductively over skeleta.



\begin{defn}
	Let $M$ be a smooth manifold that has been given a triangulation $\mathcal{T}$.  Then we say that a piecewise linear map $M \rightarrow \mathbb{R}^d$ is \textit{generic} if for no $d$-dimensional simplex is it the case that its $d+1$ vertices  map to a set of points that are contained in an affine hyperplane in $\mathbb{R}^d$.  
\end{defn}

Gromov states that a generic smooth or piecewise linear map from a closed smooth manifold $M$ to a lower dimensional Euclidean space $\mathbb{R}^d$ is a stratified fibration \cite[2.4,3]{gr09}.  This fact for piecewise linear maps will be important for us, so we give a proof.

\begin{lem} \label{generic}
	Let $M$ be a closed smooth n-manifold equipped with a triangulation.  Then a  generic piecewise linear map $f$ from $M$ to $\mathbb{R}^d$, $d<n$, is a stratified fibration for some triangulation of $\mathbb{R}^d$.  
\end{lem}

\begin{proof}
	
	Since $f$ is generic, it maps each $d$-dimensional simplex $\sigma\in M^{[d]}$ by an affine linear isomorphism onto a $d$-dimensional simplex $\Delta$ in $\mathbb{R}^d$.  Here we are using the map $\Phi$ from Definition \ref{PLdefn} to identify $\sigma$ with a simplex in Euclidean space. The union of all such $\sigma$ defines a polyhedral complex structure on $\mathbb{R}^d$, which is contained in some simplicial complex structure on $\mathbb{R}^d$ \cite[Theorem 2.2]{rspltopology}.  
	
	Fix an $i$-dimensional open face $\Delta^i$ of this simplicial complex structure on $\mathbb{R}^d$.  Let $F_1,..,F_m$ be the set of all faces of the triangulation of $M$ that contain points of $f^{-1}(\Delta^i)$. Note that the dimension of each $F_j$ is greater than or equal to $i$ because $f$ is generic.  Note also that by how we constructed the triangulation of $\mathbb{R}^d$, for each $j$ it is the case that $f(F_j)$ contains $\Delta^i$.  Finally note that the restriction of $f$ to each $f^{-1}(\Delta^i)\cap F_j$ is topologically conjugate to the normal projection of some simplex in $\mathbb{R}^n$ onto the affine subspace determined by one of its faces,  restricted to the preimage of an  open set of the latter homeomorphic to a ball (i.e., there are homeomorphisms that conjugate the  restriction $f|_{f^{-1}(\Delta^i)\cap F_i}$ to a map of this kind).   Hence for any point $p$ of $\Delta^i$, the preimage of $p$ in any of the $F_j$ is an affine subspace of $F_j$ (i.e., its image under the coordinate map $\Phi$ of Definition \ref{PLdefn} is the intersection of a subspace of the ambient Euclidean space, with $\Phi(F_j)$) and the restriction of $f$ to $f^{-1}(\Delta^i) \cap F_j$ defines a fiber bundle over the image of $f^{-1}(\Delta^i) \cap F_j$ under $f$.  By inducting on skeleta one can show that the restrictions of $f$ to $f^{-1}(\Delta^i) \cap F_j$  piece together to give a trivialization of $f$ over the interior of $\Delta^i$. Therefore $f$ defines a fiber bundle over the interior of $\Delta^i$.  
	
	The second condition of Definition \ref{sf} can also be checked directly but its verification is a bit more involved.  One takes the sets $U^{i-1}$ to be small tubular neighborhoods of the $(i-1)$-skeleton in the $i$-skeleton, for which there is a homotopy $h_i(s,t)$ from $U^{i-1}$ to the $(i-1)$-skeleton that extends to a homotopy of the identity map on the $i$-skeleton.  Note that the $f$-preimage of each open simplex $\Delta^i$ intersected with the interior of any $s$-dimensional face with which it has non-empty intersection is homeomorphic to $\Delta^i \times D^{s-i}$ via a homeomorphism $\Psi$, where $D^j$ is the j-dimensional open disk. We can moreover take the homeomorphism $\Psi$ to satisfy $f(\Psi(s,d))= f(\Psi(s,d'))$ for $(s,d), (s,d') \in \Delta^i \times D^{s-i}$. In order to construct the lift of $h_i$, we will need the homeomorphisms $\Psi$ to satisfy some additional properties.  We construct them inductively simplex by simplex in $M$ as follows.  

  \begin{figure}[h!] \label{figure:stratfibration}
		\includegraphics[width=80mm]{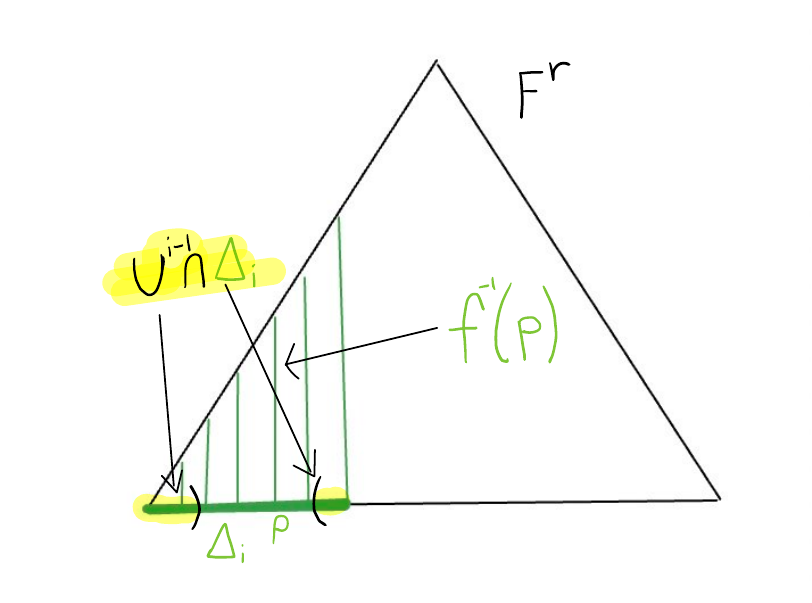}
		\caption{ } 
	\end{figure}
	
	It is enough to construct the lift of $h_i$ over the intersection of $U^{i-1}$ with each closed simplex $\overline{\Delta}^i$ in the triangulation of $\mathbb{R}^d$, so fix such a $\Delta^i$. We construct the lift of $h_i$ over $\overline{\Delta}^i$ simplex by simplex by induction on the skeleta of $M$, starting with the $i$-skeleton.  As part of the induction we will also construct, for any $r$-dimensional simplex $F^r$ with $f^{-1}(\Delta^i)\cap F^r$ nonempty, a continuous map $\Psi_{F^r}:\overline{\Delta}^i \times \overline{D}^{r-i}$ that satisfies some conditions we now describe.    
	
	 Let $\mathcal{S}_{F^r}$ be the set of all simplices in the triangulation of $M$ whose interiors intersect $\overline{f^{-1}(\Delta^i)\cap F^r}$ nontrivially.  Choosing an arbitrary but fixed once and for all $p\in \Delta^i$ and a homeomorphism $H_{F^r}$ between $\partial D^{r-i}$ and $\partial({f^{-1}(p)\cap F^r})$,  we obtain for each $F^{r'} \in \mathcal{S}_{F^r}$ a map $\phi_{F^{r'},F^r}: \overline{D}^{r'-i} \rightarrow \partial D^{r-i}$ by using the map $\Psi_{F^{r'}}: \overline{\Delta}^i \times \overline{D}^{r'-i} \rightarrow F^{r'} $ (which we assume we have already constructed in a previous stage of the induction) to identify $\overline{D}^{r'-i} \cong \{p\} \times \overline{D}^{r'-i}$ with $\overline{f^{-1}(p)\cap F^{r'}}$ and then applying $H^{-1}_{F^r}$.  We will require $\Psi_{F^r}$ to satisfy the following conditions.

	\begin{enumerate}
		\item $\Psi_{F^r}|_{\Delta^i \times \overline{D}^{r-i}}$ is a homeomorphism onto its image 
		\item $f(\Psi_{F^r}(s,d))= f(\Psi(s,d'))$ for $(s,d), (s,d') \in \overline{\Delta}^i \times \overline{D}^{r-i}$
		
		\item  For each $F^{r'} \in \mathcal{S}_{F^r}$, the following compatibility condition is satisfied: 
		\[
		\Psi_{F^{r'}}(q,d) = 	\Psi_{F^{r}}(q, \phi_{F^{r'},F^r}(d))
		\]
	\end{enumerate}
	

	For the base case, we know that $f$ maps the intersection of $f^{-1}(\Delta^i)$ with each $i$-simplex $F^i$  in $M$, if it is non-empty, homeomorphically onto $\Delta^i$.  We define $\Psi_{F^i}: \overline{\Delta}^i \times D^0 \rightarrow F^i$ by $\Psi_{F^i}(x) = f^{-1}(x) \cap F^i$.  We can then define $h_i^{\text{lift}}$ on $f^{-1}(U^{i-1}) \cap F^i$ by
 
	\[
	h_i^{\text{lift}} (x,t) = h_i( f(x),t).  
	\] 	
	
Next, assume we have made it to the $(r-i+1)$th stage of the induction. Let $F^{r}$  be an $r$-simplex in $M$ intersecting $f^{-1}(\Delta^i)$ nontrivially.   First we define $\Psi_{F^r}$  on $\overline{\Delta}^i \times \partial D^{r-i}$.  For each $(q,d) \in \overline{\Delta}^i \times \partial D^{r-i}$, we can find $F^{r'} \in \mathcal{S}_{F^r}$ so that for some $d' \in \overline {D}^{r'-i} $,  $ \phi_{F^{r'},F^r}(d')=d$.  We then define

 \[
\Psi_{F^r} ( q,d) = \Psi_{F^{r'}}(q,d').  
\]

It may be the case that there are multiple $F^{r'}$s we could have used to define $\Psi_{F^r}$ on $(q,d)$, but the compatibility condition (3) above ensures that each choice gives the same value.  This defines $\Psi_{F^r}$  on $\overline{\Delta}^i \times \partial D^{r-i}$.  

Next we define $\Psi_{F^r}$  on $\Delta^i \times \{0\}$. We do this by setting $\Phi^r(q,0)$ equal to the center of mass of $f^{-1}(q) \cap F^r$, where by using the local coordinate map $\Phi$ we can identify $F^{r}$ with a simplex $\Delta_{F^r}$ in $\mathbb{R}^r$.   Recall that the center of mass of a compact k-dimensional Hausdorff measurable subset $S$ of $\mathbb{R}^r$ is equal to the average of the coordinate vector $\Vec{x}$ with respect to Hausdorff measure $\mathcal{H}^k_S$ on $S$ normalized to have unit mass: 
\[
\int_{S} \vec{x} d \mathcal{H}^k_S(\Vec{x}).  
\]

As $q$ varies, the subsets $f^{-1}(q) \cap F^{r}$ will correspond to the preimages of the projection $\pi$ of $\Delta_{F^r}$ onto some affine plane $\Pi$ in $\mathbb{R}^r$.  The expression
\[
q \mapsto \text{center of mass}( \pi^{-1}(q) \cap \Delta_{F^r}) 
\]
defines a continuous map for $q$ contained in the closure of $\pi(\Delta_{F^r})$.  Therefore setting $\Phi^r(q,0)$ equal to the center of mass of $f^{-1}(q) \cap F^r$ defines a continuous map $\Delta^i \times \{0\} \rightarrow F^r$ that extends to a continuous map $\overline{\Delta}^i \times \{0\} \rightarrow F^r$.  

Finally, we identify $\overline{F}^r$ with the simplex $\Phi(\overline{F}^r)$ in $\mathbb{R}^r$, for $\Phi$ the map from Definition \ref{PLdefn}, and extend $\Psi_{F^r}$ radially: we map the line segment joining $(q,0)$ to $(q,v) \in \overline{\Delta}^i \times \partial D^{r-i}$ to the line segment joining their images in $\overline{F}^r$ (using again the map $\Phi$ to identify $\overline{F}^r$ with a simplex in $\mathbb{R}^d$ so that this makes sense.)  The map $\Psi_{F^r}$ is then a homeomorphism onto its image when restricted to $\Delta^i \times \overline {D}^{r-s}$, although it might not be on its whole domain  $\overline{\Delta}^i \times \overline {D}^{r-s}$. For example, for some $q\in \partial \Delta^i$  $\Psi_{F^r}$ might map $\{q\} \times \overline {D}^{r-s}$ to a single point, see e.g. figure \ref{figure:stratfibration}.  

We can then define $h_i^{\text{lift}}$ by setting 
\[
h_i^{\text{lift}}(x,t) = x 
\]

\noindent for $x \in f^{-1}(\partial \Delta^i \cap \overline{F}^r)$.  All other $x \in f^{-1}(U^{i-1})$ are equal to $\Psi_{F^r}(q,d)$ for a unique $q$ and $d$ because we constructed $\Psi_{F^r}$ so that it satisfies condition (1), so for $x \in f^{-1}(U^{i-1}) -f^{-1}(\partial \Delta^i \cap \overline{F}^r)$  we can set 
\[
h_i^{\text{lift}}(x,t) = \Psi_{F^r} (h_i(q,t),d).  
\]

That the compatibility condition (3) is met for $\Psi_{F^r}$ ensures that this definition of $h_i^{\text{lift}}$ agrees with the definition given on the $\overline{F}^{r'}$  that have been covered by a previous stage in the induction.  This completes the induction and constructs the $h_i^{\text{lift}}$.

\end{proof}

\section{Topological Waists}\label{sec-topwaist}




Recall the discussion in Section \ref{section:naturalvectorfields}, and the definition of $\Omega_d$.  The growth indicator function $\psi_{\Lambda}:\mathfrak{a}^+ \rightarrow \mathbb{R}$ of a discrete subgroup $\Lambda \subset \text{Isom}(X)$ measures the exponential growth rate of elements of $\Lambda$ with Cartan coordinate restricted to a cone in $\mathfrak a^{+}.$ In this section we do not require $X, G=\Isom(X)$ to be a rank one symmetric space, rank one semisimple Lie group respectively.


\begin{thm} \label{largefiber}
	 Let $\Gamma$ be a uniform lattice in $G$, and let $f:\Gamma \backslash X \rightarrow \mathbb{R}^d$ be a generic piecewise linear map. Then for some $p\in \mathbb{R}^d$ the image $\Lambda_p$ of the inclusion of the fundamental group of some connected component of the preimage  $f^{-1}( \{p\})$  in $\Gamma \cong \pi_1(\Gamma \backslash X)$ satisfies the following: for no $\eta \in \Omega_d$ is it the case that the growth indicator function $\psi_{\Lambda_p}$ of $\Lambda_p$ satisfies $\psi_{\Lambda_p} < \eta$ on $\mathfrak{a}^+.$ 
  
  

\end{thm}

Theorem \ref{intro:critical exponent} in the introduction is just the rank one case of Theorem \ref{largefiber} by the fact that the growth indicator function and critical exponent are related by (\ref{criticalexponentgrowthindicator}) and Definition \ref{defn-CXk}.  Theorem \ref{mt-topwaist} in the introduction regarding the $\SL_n(\mathbb R)$ case follows from Theorem \ref{largefiber}, the fact verified in the proof of Proposition \ref{prop-novereight} that $2\rho -\Theta$ is in $\Omega_d$ for $d$ in the range specified in the statement of the theorem, and the Quint-Lee-Oh \cite{quint2003propriete, LeeOh24, oh02} gap theorem \ref{thm-indfctiongap}.

\begin{proof} 
  Assume for contradiction that the fundamental groups of all connected components of all preimages of points under $f$ include to subgroups $\Lambda$ of $\pi_1(\Gamma \backslash X)$ with $\psi_{\Lambda}$ bounded above by some $\eta \in \Omega_d$.   Denote by $\mathcal{T}$ the triangulation of $\mathbb{R}^d$ given by the Lemma \ref{generic} and relative to which $f$ is a stratified fibration. The subsequent steps of the proof are split into Sections \ref{sec-opencover},\ref{sec-PsPHlifts},\ref{sec-piecing},\ref{sec-newsweepouts} and \ref{sec-conclusion}.
	
	\subsection{Construction of an open cover adapted to $\mathcal{T}$}\label{sec-opencover}
	We will write $\Delta^m$ for a simplex of dimension $m$. This will lead to slight abuses of notation but since we will never explicitly refer to more than one simplex of the same dimension simultaneously there should be no confusion. For each  simplex $\Delta^m$ in $\mathcal{T}$ denote by $\Delta_{\rm int}^m(\epsilon)$ the set of points in $\Delta^m$ at a distance \textit{greater than} $\epsilon$ from the boundary of $\Delta^m$ for $m>0$, and for $m=0$ we let $\Delta_{\rm int}^m(\epsilon)$ equal $\Delta^m$ (a single point.) We can choose $\epsilon(m)>0$ and neighborhoods $N(\Delta_{\rm int}^m(\epsilon(m)))$ of $\Delta_{\rm int}^m(\epsilon(m))$ in $\mathbb{R}^d$ so that the following properties hold.

 \begin{enumerate}
\item The $N(\Delta_{\rm int}^m(\epsilon(m)))$ are an open cover of $\mathbb{R}^d$
\item  For two simplices $\Delta^{m'}$ and $\Delta^m$ 
\[
N(\Delta_{\rm int}^m(\epsilon(m))) \cap N(\Delta_{\rm int}^{m'}(\epsilon(m')))= \emptyset
\]
unless $\Delta^{m'}$ is contained in the closure of $\Delta^m$ or vice versa.  
\item There is a deformation retraction from $N(\Delta_{\rm int}^m(\epsilon(m)))$  to a subset of $\Delta_{\rm int}^m(\epsilon(m)/2)$ and this deformation retraction lifts to a homotopy $h^m$ of the identity map of $f^{-1}(N(\Delta_{\rm int}^m(\epsilon(m))))$  to a map with image contained in $f^{-1}(\Delta_{\rm int}^m(\epsilon(m)/2))$ 
 \end{enumerate}

 The $N(\Delta_{\rm int}^m(\epsilon(m)))$ can be constructed by inducting on $m$, starting with $m=0$, in the following way.  After the $m$-th stage in the induction we can ensure that the $N(\Delta_{\rm int}^{m'}(\epsilon(m')))$ we have constructed for $m' \leq m-1$ cover the $(m-1)$-skeleton of the triangulation of $\mathbb{R}^d$. 
 We can also ensure that for any face $\Delta^m$ of the triangulation of $\mathbb R^d$ the set $\Delta^m\setminus\bigcup N(\Delta^{m'}_{\rm int}(\epsilon(m'))$, where the sum is taken over all the neighborhoods $N(\Delta^{m'}_{\rm int}(\epsilon(m')), m'\leq m$ constructed up to this point, is homeomorphic to the interior of $\Delta^m$ and is contained in $\Delta^m_{\rm int}(\epsilon(m))$ for $\epsilon(m)$ chosen sufficiently small. 
 Finally, we can show that condition (3) is met, provided $\epsilon(m)$ and the neighborhoods $N(\Delta_{\rm int}^m(\epsilon(m)))$ were chosen sufficiently small, by using the fact that $f$ is a stratified fibration and repeatedly applying condition (2) in Definition \ref{sf}. 

 \begin{figure}[h!]
		\includegraphics[width=40mm]{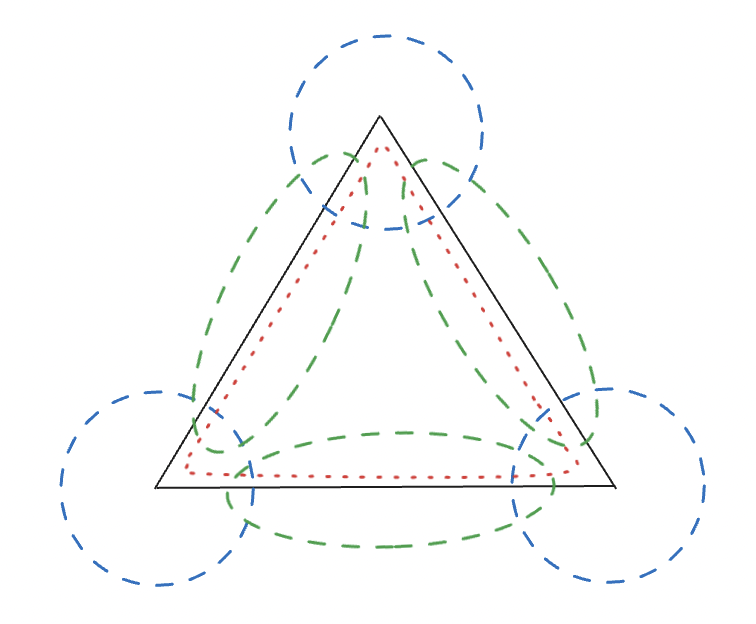}
		\caption{Cover of $\Delta^m$ by the $N(\Delta_{\rm int}^{m'}(\epsilon(m')))$ } 
	\end{figure}

	Choose a partition of unity for $\mathbb R^d$ subordinate to the open cover by the sets $N(\Delta_{\rm int}^m(\epsilon(m)))$.  Denote by $\phi_m$ the element of the partition of unity corresponding to $N(\Delta_{\rm int}^m(\epsilon(m)))$. 
    
    \subsection{PSpH functions for lifts of preimages of $f$}\label{sec-PsPHlifts}
    
    Because $f$ is a stratified fibration relative to $\mathcal{T}$, we know that for any $\Delta^m$, $f$ defines a trivial fiber bundle of $f^{-1}(\Delta^m_{\rm int})$ over the interior $\Delta^m_{\rm int}$ of $\Delta^m$. Therefore choosing some $p_m \in \Delta^m$, we have that $f^{-1}(\Delta^m_{\rm int}) \cong \Delta^m_{\rm int} \times f^{-1}(p_m)$.  Denote by $C_\ell(p_m)$, where $\ell$ ranges over $\ell=1,..,\overline{\ell}$, the connected components of $f^{-1}(p_m)$. 
	
	To the lifts $L_{\alpha}=L_{\alpha}(\ell,p_m)$ of  $C_{\ell}(p_m)$ to the universal cover $X$, where $\alpha$ ranges over the cosets of $\pi_1(C_{\ell}(p_m))$ in $\pi_1( \Gamma \backslash X)$, we associate $(n-d)$-PSpH functions $\psi_{p_m,\ell}^{i}:X \rightarrow \mathbb{R}$ as follows. Note that each $L_{\alpha}$ is invariant under the action of a subgroup $\Gamma_\alpha$ of the deck group $\Gamma$ conjugate to the image of the inclusion of the fundamental group of $C_{\ell}(p_m)$ in the fundamental group of $\Gamma \backslash X$ with some basepoint on $C_{\ell}(p_m)$.  Note also that for any given lift $L_{\alpha}$ preserved by $\Gamma_\alpha < \Gamma$, the left cosets of $\Gamma_\alpha$ are in bijection with the lifts of $C_{\ell}(p_m)$ to $X$.  This bijection is canonical given a choice of lift of $C_{\ell}(p_m)$ to $X$: for each $C_{\ell}(p_m)$, fix such a lift 
 $L_{\overline{\alpha}}$

We choose a  $\Gamma_{\overline{\alpha}}$-invariant PSpH function $\psi_{p_m,\ell}^{\overline{\alpha}}: X \rightarrow \mathbb{R}$.   We know there is such a function by our assumption that $\Gamma_{\overline{\alpha}}$ includes to a subgroup of $\pi_1(\Gamma \backslash X)$ whose growth indicator function $\psi_{\Gamma_{\overline{\alpha}}}$ is bounded above by some $\eta$ in $\Omega_d$. By scaling  $\eta$ by a factor less than or equal to $1$ so that it is tangent to $\psi_{\Gamma_{\overline{\alpha}}}$ but still contained in $\Omega_d$, which is possible by Lemma \ref{lem-omegakconvex}, the existence of the $\Gamma_{\overline{\alpha}}$-invariant PSpH function follows from Theorem \ref{thm-existenceofpattersonsullivan}.

For any coset $H_\alpha$ of $\Gamma_{\overline{\alpha}}$ corresponding to the lift $L_\alpha$, we choose a representative $h_\alpha$ of $H_\alpha$ and define $\psi_{p_m,\ell}^{\alpha}$ to be the composition of $ \psi_{p_m,\ell}^{\overline{\alpha}}$ with the isometry of $X$ given by $h_{\alpha}^{-1}$.  Because $ \psi_{p_m,\ell}^{\overline{\alpha}}$ is $\Gamma_{\overline{\alpha}}$-invariant, this is well-defined independent of the choice of coset representative $h_{\alpha}$.  Note that $ \psi_{p_m,\ell}^{\alpha}$ is $h_{\alpha} \Gamma_{\overline{\alpha}} h_{\alpha}^{-1}$-invariant.

	\subsection{Piecing together the PsPH functions }\label{sec-piecing}
	Since $f$ is a generic PL map, by Proposition \ref{prop:sweepout} the preimages of $f$ give a sweepout of $\Gamma \backslash X$, which we denote by $\mathcal{F}$.  This means that, adding a point at infinity to $\mathbb{R}^d$, the preimages of $f$ define a continuous map $\rho_f: \mathbb S^d \rightarrow cl_{n-d}(\Gamma \backslash X)$, where  $ cl_{n-d}(\Gamma \backslash X)$ is the space of $(n-d)$-cycles in the flat topology, so that the induced map $\rho_f^*: H^{d}(cl_{n-d}(\Gamma \backslash X);\mathbb{F}_2) \rightarrow H^d(\mathbb S^d;\mathbb{F}_2)$ is nontrivial.  Points outside of a compact set are mapped to the empty cycle under this map. 

	We will define a continuous family of sweepouts $\mathcal{F}_t$ presently, so that $\mathcal{F}_0=\mathcal{F}$.  Choose a point $p \in \mathbb R^d$.  The preimage $f^{-1}(p)$ may have several connected components which we label $C_\ell(p)$.

	Suppose that $p$ is contained in $N(\Delta_{\rm int}^{m_1}(\epsilon(m_1))),..., N(\Delta_{\rm int}^{m_i}(\epsilon(m_i)))$, for $m_1<..<m_i\leq d$.  Fix  a connected component $C_\ell(p)$ of $f^{-1}(p)$. Then for each $\Delta^{m_r}(\epsilon(m_r))$, the homotopy $h^{m_r}$ from (3) above takes $C_\ell(p)$ into some $\Delta_{\rm int}^{m_r}(\epsilon(m_r)/2) \times C_{\ell_r}(p_{m_r})$, where we have identified $f^{-1}(\Delta_{\rm int}^{m_r}(\epsilon(m_r)/2))$ with $\Delta_{\rm int}^{m_r}(\epsilon(m_r/2)) \times f^{-1}(p_{m_r})$.

Choose a connected component $\mathcal{L}(i)$ of the lift of $C_{\ell_i}(p_{m_i})$ to the universal cover $X$.  Then one can lift the homotopies $h^{m_r}$ to homotopies of $\mathcal{L}(i)$ that take $\mathcal{L}(i)$ into lifts $\mathcal{L}(r)$ of $C_{\ell_r}(p_{m_r})$. Let $\Gamma(C_{\ell_r}(p_{m_r}))$ be the copy of the inclusion of $\pi_1(C_{\ell_r}(p_{m_r}))$ in $\Gamma$ that preserves $\mathcal{L}(r)$.  Then if $r<r'$, it will be the case that $\Gamma(C_{\ell_r}(p_{m_r}))$ contains $\Gamma(C_{\ell_{r'}}(p_{m_{r'}}))$. This follows from the fact that $f$ is a stratified fibration and property (2) of Definition \ref{sf}.  Denote the  PSpH function for $\mathcal{L}(r)$ by $\psi_r$.  Then for each $p$ the sum
	\begin{equation} \label{psfoverp}
		\sum_{r=1}^{i} \phi_{m_r} (p) \psi_r
	\end{equation}
	defines a $\Gamma(C_{\ell_i}(p_{m_i}))$-invariant function on $X$.   Note that if $\phi_{m_i} (p)=0$ then we get  a $\Gamma(C_{\ell_{i-1}}(p_{m_{i-1}}))$-invariant function, and so on.  Note also that this sum defines a PSpH function, as a convex combination of PSpH functions, which varies continuously with $p$. We emphasize that $p$ is regarded as fixed in (\ref{psfoverp}), so that in particular no derivatives in the $p$ variable matter for whether that expression defines a PSpH function.  
	
	\subsection{Construction of new sweepouts $\mathcal{F}_t$}\label{sec-newsweepouts}
	
	We define the sweepout $\mathcal{F}_t: \mathbb S^d \rightarrow cl_{n-d}(\Gamma \backslash X)$ at the point $p$ we began with as follows.  Suppose $p\in \Delta^{\overline{i}}$, where necessarily $\overline{i} \geq i$. Each connected component $C_\ell(p)$ lifts to a collection of cycles in $X$ that are in bijection with the cosets of $\Gamma(C_{\ell_{\overline{i}}}(p_{m_{\overline{i}}}))$ in $\pi_1( \Gamma \backslash X)$.   Each such lift $\tilde{C}_\ell(p)$ determines a lift $\mathcal{L}(i)$ of $C_{\ell_i}(p_{m_i})$: the homotopy $h^{m_i}$ from (3) above determines a homotopy from $C_\ell(p)$ into $C_{\ell_i}(p_{m_i})$, which when lifted to the universal cover gives a homotopy from $C_{\ell}(p)$ into the lift of $C_{\ell_i}(p_{m_i})$ that we call $\mathcal{L}(i)$.  This in turn determines the $\mathcal{L}(r)$, $1 \leq r \leq i-1$, as above.  
 
 We can therefore flow $\tilde{C}_\ell(p)$ for time $t$ by the gradient of (\ref{psfoverp}) to obtain a new cycle which we denote by $\tilde{C}_\ell(p)(t)$.  Note that $\tilde{C}_\ell(p)(t)$ is  $\Gamma(C_{\ell_i}(p_{m_i}))$-invariant since this is the case for $\tilde{C}_\ell(p)$ and the function (\ref{psfoverp}).  We can then define $\mathcal{F}_t(p)$ to be the union over $\ell$ of the covering projections of the  $\tilde{C}_\ell(p)(t)$ down to $\Gamma \backslash X$.  This is well-defined independent of the choice of lift of $C_\ell(p)$ by the fact that the $\psi_r$ associated to the $\mathcal{L}(r)$ for two different lifts  $\tilde{C}_\ell(p)$ are related by precomposition with an isometry given by a $\gamma \in \Gamma $, independent of $r$, such that $\gamma H_r' = H_r$, where $H_r$  and $H_r'$ are the cosets of the image of the fundamental group of the $C_{\ell_r}(p_{m_r})$ in $\Gamma$ corresponding to the $\mathcal{L}(r)$.  
	
Although this observation is not necessary for the proof, since it might shed some light on our construction we point out that the sweepouts $\mathcal{F}_t$ we have constructed are each the image of $\mathcal{F}$ under a certain map $\boldsymbol{\Phi}_t$. This can be seen as follows.  Each $x\in \Gamma \backslash X$ is contained in some $C_{\ell}(f(x))$.  For each lift $\tilde{x}$ of $x$ to $X$, we then flow $\tilde{x}$ for time $t$ by the flow of the equivariant PSpH function (\ref{psfoverp}) corresponding to the lift $\mathcal L(i)$ of $C_{\ell}(f(x))$ containing $\tilde{x}$, and define $\boldsymbol{\Phi}_t(x)$ to be the projection of the result to $X$. The discussion of the previous paragraph shows that this does not depend on which lift $\tilde{x}$ of $x$ we chose, and so $\boldsymbol{\Phi}_t$ is well-defined, and by construction $\boldsymbol{\Phi}_t(\mathcal{F})= \mathcal{F}_t$.  

\subsection{Conclusion}\label{sec-conclusion}

	The cycles of $\mathcal{F}_t$  vary continuously in the flat topology.  The  $\mathcal{F}_t: \mathcal S^d \rightarrow cl_{n-d}(\Gamma \backslash X)$ are also nontrivial sweepouts because they are homotopic to $\mathcal{F}_0 = \mathcal{F}$, which is a nontrivial sweepout.  Moreover, since the gradients of the $\psi_r$ decrease $(n-d)$-dimensional area at a definite rate since the $\psi_r$ are $(n-d)$-PSpH, and this property is preserved under convex combinations, we know by Theorem \ref{contraction} that for any $\epsilon'$ we can choose $t$ large enough that every cycle in  $\mathcal{F}_t$ has area smaller than $\epsilon'$.   Here we are using the fact that the $\psi_r$ descend to a quotient of $X$ by a subgroup of $\Gamma$ isomorphic to $\Gamma(C_{\ell_i}(p_{m_i}))$ to which $C_\ell(p)$ lifts, and we can compute the rate of change of the area of $C_{\ell}(p)$ in this cover. But this contradicts the fact that the $(n-d)$-width of  $\Gamma \backslash X$ is positive (Proposition \ref{positivewaist}.)

\end{proof} 

\section{Property $\FA_r$ and the proof of Theorem \ref{mt-finiteorbit}} \label{section:thm2}
 The results of the previous section can be used to restrict the possible actions of $\Gamma$ on low dimensional simplicial complexes $Y$.  


\begin{thm} \label{thm-finiteorbit}
	Assume that $\Gamma$ is a lattice in a symmetric space $X$, and that it acts by  simplicial automorphisms on a contractible $d$-dimensional simplicial complex $Y$, for $d< r(X)$. Then there is a finite orbit for the action of $\Gamma$.   
	\end{thm} 

 It seems likely that the condition that $\Gamma$ act by simplicial automorphisms could be weakened (see \cite[II.6]{bh99}.)  In any case, the condition is satisfied by the examples that naturally occur (e.g., the action of a lattice on the corresponding symmetric space or the action of a group on its Bass-Serre tree, if the group splits as an amalgam or an HNN-extension)  

We restate Corollary \ref{corfanover8} from the introduction.
\begin{repcor}{corfanover8}
\begin{enumerate}\item Cocompact lattices in $\SL_n(\mathbb R)$ have property $\FA_{\lfloor \frac{n}{8}\rfloor-1}.$
\item Cocompact lattices in $\rm{F}_4^{-20}$ have property $\FA_{2}.$
\end{enumerate}
\end{repcor}

\begin{proof} 

Suppose we have an action of a lattice $\Gamma$ by isometric simplicial automorphisms on a $\text{CAT}(0)$ simplicial complex $Y$ of the appropriate dimension.  We know by the previous theorem that the action $\Gamma$ has a finite orbit:   there is some $p\in Y$ whose orbit under $\Gamma$ is a finite set of points. Then if the orbit of $p$ consists of the points $p_1,..,p_l$, the unique minimum of $g(y):= \sum_{i=1}^l d^2(y,p_i)$  is a fixed point of $\Gamma$ (see \cite[II.2]{bh99}.) 

 \end{proof}

The following lemma will be useful in the proof of Theorem \ref{thm-finiteorbit}.

\begin{lem} \label{lemma-epsilonapart}
Let $X$, $\Gamma$ be as in the statement of Theorem \ref{thm-finiteorbit}, and assume that $\Gamma$ acts (not necessarily properly discontinuously) by  simplicial automorphisms on a contractible $d$-dimensional simplicial complex $Y'$.  Suppose that there is a $\Gamma$-equivariant map $X \rightarrow Y'$ that descends to a map $f: \Gamma \backslash X \rightarrow \Gamma \backslash Y':=Z$ between the quotients by the actions, and assume that there is a point $p\in Z$ with the following property: some connected component of $f^{-1}(U_p)$ virtually $\pi_1$-surjects onto $\Gamma \backslash X $ for every open neighborhood $U_p$ of $p$.  Then some connected component of $f^{-1}(p)$ virtually $\pi_1$-surjects onto $\Gamma \backslash X $.  
\end{lem}

\begin{proof}

Equip $Y'$ with a metric so that all simplices of a given dimension are isometric to a fixed simplex in Euclidean space.  Then $\Gamma$ acts by isometries on $Y'$, and for any point $y\in Y$ there is some $\epsilon'$ that only depends on the orbit of $y$ so that if $\gamma \cdot y \neq y$, then the distance between $\gamma \cdot y$ and $y$ is greater than $\epsilon'$.  


The quotient $Z=\Gamma \backslash Y'$ inherits a metric from $Y'$.   Choose a ball $B$ in $Z$ centered at $p'$ and with radius less than $\epsilon'/2$.  Let $\rho: Y' \rightarrow \Gamma \backslash Y'$  be the quotient map.  By the previous paragraph $\rho^{-1}(B)$ is a disjoint union of connected regions $\tilde{B}_i$ with $\rho(\tilde{B}_i)=B$.


 
 We know that the fundamental group of some connected component $C_{f^{-1}(B)}$ of $ f^{-1}(B)$ includes to a finite index subgroup of $\pi_1(\Gamma \backslash X)$, and that the action of each element $\gamma \in \pi_1(C_{f^{-1}(B)})$ on $Y'$ preserves $\rho^{-1}(B).$  We claim that the action of $\gamma$ also preserves each $\tilde{B}_i$.  Fix $\tilde{B}_i$, and choose a lift $\tilde{c}_{\gamma}$ of a loop $c_{\gamma}$ representing $\gamma$ in $\pi_1(C_{f^{-1}(B)})$.  We can choose $\tilde{c}_{\gamma}$ so that its startpoint  maps to $\tilde{B}_i$. Then because $f(\tilde{c}_{\gamma})$ contains one point of $\tilde{B}_i$, and the connected component of $f(\tilde{c}_{\gamma}(0))$ in $\rho^{-1}(B)$ is $\tilde{B}_i$, it must be the case that all of $f(\tilde{c}_{\gamma})$ is contained in $\tilde{B}_i$.  It follows that $\gamma \cdot \tilde{B}_i \cap \tilde{B}_i \neq \emptyset$, and so because we chose $\epsilon'$ as above, it must be the case that $\gamma$ fixes the center of $\tilde{B}_i$. 

 This proves that the connected component of $f^{-1}(p')$ contained in $C_{f^{-1}(B)}$ virtually $\pi_1$-surjects onto $\Gamma \backslash X$. 
 
\end{proof}

We now prove Theorem \ref{thm-finiteorbit}.  

 \begin{proof}

To start off we give $Y$ a metric so that each simplex is isometric to the regular simplex of that dimension with edges of unit length.  Then since $\Gamma$ acts by simplicial automorphisms, it also acts by isometries.

There is a natural fibering $\Gamma \backslash (X \times Y) \rightarrow \Gamma \backslash X$, where the fibers are copies of $Y$.  Since the fibers are contractible, we can find a section $s: \Gamma \backslash X \rightarrow \Gamma \backslash (X \times Y)$. Choose a triangulation of $\Gamma \backslash X$. By lifting the triangulation to $X$ we obtain a $\Gamma$-invariant triangulation of $X$. Let $\overline{s}$ be the $\Gamma$-equivariant map $X \rightarrow Y$ obtained from $s$.   By simplicial approximation (which works in the same way for $\Gamma$-equivariant maps of the kind considered here) we can homotope $\overline{s}$ to a nearby simplicial map $F$, provided we took the initial triangulation of $\Gamma \backslash X$ fine enough that $s$ maps adjacent vertices to points at distance less than $1/2$.  Note that $F(X)$ is a subcomplex of $Y$.

By \cite[Section 6.5]{a79} (see the remark on \cite[pg. 146]{a79}, that states that the results of Section 6.5 are also true for locally finite complexes with infinitely many simplices, such as $F(X)$), the quotient $Z:=\Gamma \backslash F(X)$ has the structure of a simplicial complex,  where moreover the quotient map from $F(X)$ to $Z$ can be taken to be a simplicial map, possibly after subdividing $F(X)$. The statements of \cite[Section 6.5]{a79} are for simplicial complexes in Euclidean space, but  $F(X)$  can be linearly embedded in any Euclidean space of large enough dimension.  Note that $Z$ may have boundary: in the case of a group action on a tree, for example, we could get a closed interval. The quotient $Z$ is compact, since $\Gamma$ is cocompact, and $Z$ is the image of the closure of a fundamental domain of $\Gamma$.  Let $\underline{F}: \Gamma \backslash X \rightarrow Z$ be the map to which $F$ descends.

Choose some metric on $Z$ so that all simplices of a given dimension are isometric to a regular simplex in Euclidean space with edges of unit length. Fix a simplicial map $\pi: Z \rightarrow \mathbb{R}^d$,  which is generic in the sense that it maps each $d'$-dimensional simplex to a non-degenerate $d'$-dimensional simplex in $\mathbb{R}^d$ for $d' \leq d$.  Such a map $\pi$ is determined uniquely by where it sends the vertices of $Z$.  Note that the preimage of any point under $\pi$ is a finite set.  

Even though $\pi$ is generic, $\pi \circ \underline{F}$ might fail to be. For each $m\in \mathbb N$ we perturb the map $\pi \circ \underline{F}$ to a generic piecewise linear map $g_m: \Gamma \backslash X \rightarrow \mathbb{R}^d$ relative to a triangulation of $\Gamma \backslash X$ by simplices of diameter smaller than $1/m$, so that the $C^0$ distance between $\pi \circ \underline{F}$ and $g_m$ tends to zero as $m$ tends to infinity.  By Theorem \ref{largefiber}, the fundamental group of some connected component $C_{p_m}$ of the preimage of some point $p_m\in \mathbb R^d$ under $g_m$ includes to a finite index subgroup of $\pi_1(\Gamma \backslash X)$.  By passing to a subsequence, we can assume the $p_m$ converge to $p \in \mathbb{R}^d$.  Choose disjoint balls $N_1,..,N_\ell$ centered at the finitely many preimages of $p$ under $\pi$. For $m$ large enough,  $C_{p_m}$ is contained in $\underline{F}^{-1}(N_{\ell(m)})$ for some $\ell(m)\in \{1,\ldots, l\}$, and so the fundamental group of some connected component of $\underline{F}^{-1}(N_{\ell(m)})$ includes to a finite index subgroup of $\pi_1(\Gamma \backslash X)$.   By making the radii of the $N_i$ tend to zero and repeating this argument we can find a point $p'\in \pi^{-1}(p)$ so that the preimage under $\underline{F}$ of every open set containing $p'$ has a connected component that virtually $\pi_1$-surjects.  Therefore by Lemma \ref{lemma-epsilonapart}, the preimage of of $p'$ has a connected component that virtually $\pi_1$-surjects, and so the orbit of $p'$ under $\Gamma$ is finite, which completes the proof. 

\begin{figure}[h!]
		\includegraphics[width=80mm]{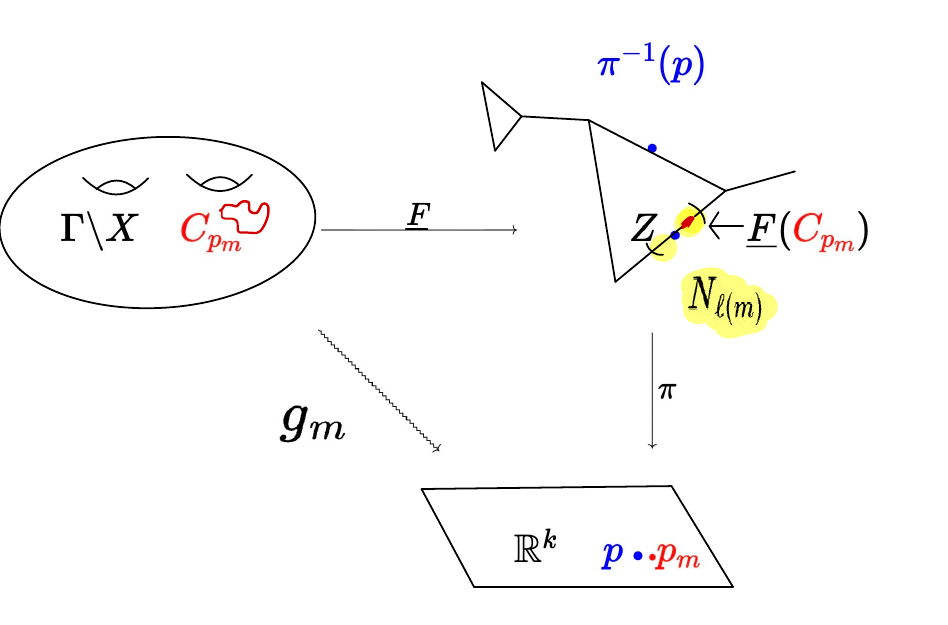}
		\caption{} \label{diagram:maps}
	\end{figure}

 \end{proof}

\section{Branched Cover Stability and non-abelian expansion}\label{sec-branched}
The main result of this section is Theorem \ref{mt-branched}. We recall the statement below. 
\begin{reptheorem}{mt-branched} Let $\varepsilon>0$. There exists $\alpha, c>0$ with the following properties. Let $M$ be a compact octonionic manifold and let $\iota_1\colon M_1\to M$ be a branched cover of $M$ with branching locus  a rectifiable $14$-varifold $N\subset M$. Suppose, $\mathcal H^{14}(N)\leq \alpha \vol(M).$ Then, there exists a rectifiable $15$-varifold $W\subset M$ with $\mathcal H^{15}(W)\leq c\mathcal H^{14}(N)$, a genuine cover $\iota_2\colon M_2\to M$ and an isometry $\mu\colon (M_1\setminus \iota_1^{-1}(W))\to (M_2\setminus \iota_2^{-1}(W))$ such that $\iota_2\circ \mu=\iota_1.$ 
\end{reptheorem}

One can view this as a stability result for covers of octonionic hyperbolic manifolds: any branched cover of an octonionic hyperbolic manifold with small branching locus $N$ can be cut and re-glued along a small codimension one subset $W$ to completely remove the branching and turn it into a genuine cover. We prove Theorem \ref{mt-branched} in Section \ref{subsection:riemanniancoverstability}. The geometric picture explained above is a Riemannian analogue to the type of cover stability studied in \cite{dinur2022near} and \cite[4.4]{chapman2023stability}. In fact, we use Theorem \ref{mt-branched} to produce families of simplicial complexes which are covering stable in the sense of \cite{dinur2022near}. Let $H$ be a group and let $h_1(Y,H)$ be the non-abelian Cheeger constant defined in \cite[1.3] {dinur2022near} (also in the formula (\ref{eq-nacheeger})). 
\begin{thm}\label{thm-CovStab}
Let $\Gamma\subset F_4^{(-20)}$ be a cocompact lattice. Let $\mathcal T$ be a triangulation of $\Gamma\bs X$ and let $Y$ be the resulting simplicial complex. There is a $\delta>0$ such that 
$h_1(Y_1,H)\geq \delta$, for any finite cover $Y_1$ of $Y$ and any group $H$. 
\end{thm}
In \cite[Thm 1.8]{dinur2022near}, it was proved that spherical complexes $A_3(\mathbb F_q)$ admit a uniform lower bound on the first non-abelian Cheeger constant. It was also known \cite[Prop 3.1]{meshulam2013bounded} for simplices $\Delta^n, n\in\mathbb N.$ The first bounded geometry families of simplicial complexes with this property were exhibited in the work of Dikstein-Dinur \cite{Dikstein2023CoboundaryAC}, where they were obtained from finite quotients of Bruhat-Tits buildings of higher rank $p$-adic groups, for $p$ high enough. As far as we are aware, Theorem \ref{thm-CovStab} provides the first Riemannian examples of a bounded geometry family of non-commutative expanders in the sense of \cite{dinur2022near}. In plain terms, Theorem \ref{thm-CovStab} asserts that an ``almost cover" of $Y_1$ can be made into a genuine cover, by modifying it over a small number of simplices. We define the non-commutative Cheeger constant and prove Theorem \ref{thm-CovStab} in Section \ref{sec-ccstability}. In Section \ref{subsection:moving} we develop several methods to ``adjust" branched covers by moving the branching locus, which might be of independent interest.

\subsection{Moving the branching locus} \label{subsection:moving} 
We collect several lemmas that will allow us to ``move" the branching locus of a branched cover of a manifold. This is made precise in Lemma \ref{lem-MovingLocus}. For the initial part of this section the Riemannian structure is not important since the results have a purely topological flavour. The metric will intervene when we consider the effect of Federer-Fleming deformation on the branching loci.  Let $M$ be a manifold and let $N_1,N_2\subset M$ be closed subsets. One can think of them as candidates for branching loci, although we do not require them to be of codimension two. Let $W$ be a closed set containing both $N_1$ and $N_2$. We consider the maps between the fundamental groupoids \cite[\S 2.5]{may1999concise} $\iota_i\colon \pi_1(M\setminus W)\to \pi_1(M\setminus N_i)$ induced by inclusions $(M\setminus W)\hookrightarrow (M\setminus N_i), i=1,2.$ One does not lose anything by fixing a family of  basepoints (one per each connected component) and considering only the fundamental groups. When we say that a map between groupoids is surjective we mean that the image is a dense, full subcategory. Similarly, when we will say that a map between groupoids is an isomorphism we really mean that is is an equivalence of categories. In the fundamental group language these coincide with the usual meaning of  surjectivity and isomorphism respectively.

\begin{defn}\label{def-Wdomination}
We say that $N_1$ $W$\textit{-dominates} $N_2$, and write $N_1\xRightarrow{W} N_2$, if there exists a groupoid homomorphism $\kappa$ making the following diagram commute
\[\begin{tikzcd}[ampersand replacement=\&]
	\& {\pi_1(M\setminus N_1)} \\
	{\pi_1(M\setminus W)} \& {\pi_1(M\setminus N_2)}
	\arrow["{\iota_1}", from=2-1, to=1-2]
	\arrow["{\iota_2}"', from=2-1, to=2-2]
	\arrow["\kappa"', dashed, from=2-2, to=1-2]
\end{tikzcd}\]
\end{defn}
\begin{rem}
Equivalently, we can require the above diagram for the fundamental groups with basepoints in every path connected component of $M\setminus W.$ Similarly, the statements that follow can be replaced by analogous statements for fundamental groups with basepoints in every path connected component.
\end{rem}
The Lemma below provides a more functorial characterization of domination, and gives a simple criterion. 
\begin{lem}\label{lem-WDomExt}

\begin{enumerate}
\item $N_1\xRightarrow{W} N_2$ if and only if for every groupoid $H$ and every homomorphism $\rho_1\colon\pi_1(M\setminus N_1)\to H$ there is a homomorphism $\rho_2\colon \pi_1(M\setminus N_2)\to H$ with $\rho_1\circ\iota_1=\rho_2\circ \iota_2.$
\item Suppose that $\ker \iota_2\subset \ker \iota_1$ and $\iota_2$ is surjective. Then, $N_1\xRightarrow{W} N_2$.
\end{enumerate}
\end{lem}
\begin{proof}
\begin{enumerate}
\item If we have $\kappa$ as in Definition \ref{def-Wdomination}, we put $\rho_2:=\rho_1\circ \kappa.$ For the other direction choose $H=\pi_1(M\setminus N_1)$ and let $\rho_1$ be the identity map. Then we can take $\kappa=\rho_2.$
\item 
For any $\gamma\in \pi_1(M\setminus W)$ put $\kappa(\iota_2(\gamma)):=\iota_1(\gamma).$ The condition $\ker\iota_2\subset \ker\iota_1$ guarantees that $\kappa$ is well defined on the image of $\iota_2$. Then $\kappa$ extends to the whole groupoid as the image of $\iota_2$ is a full dense subcategory. 

\end{enumerate}
\end{proof}
The following Lemma explains the relevance of $W$-domination to the problem of moving the branching locus. 

\begin{lem}\label{lem-MovingLocus}
Suppose $N_1\xRightarrow{W} N_2$. For any  cover $\rho_1\colon M_1\to (M\setminus N_1)$ there exists a cover $\rho_2\colon M_2\to (M\setminus N_2)$ and a map $\varphi\colon (M_1\setminus \rho_1^{-1}(W))\to (M_2\setminus \rho_2^{-1}(W))$ which makes the diagram commute. 
\[\begin{tikzcd}
	{M_1\setminus \rho_1^{-1}(W)} && {M_2\setminus \rho_2^{-1}(W)} \\
	& {M\setminus W}
	\arrow["\varphi", from=1-1, to=1-3]
	\arrow["{\rho_1}"', from=1-1, to=2-2]
	\arrow["{\rho_2}", from=1-3, to=2-2]
\end{tikzcd}\]
\end{lem}
\begin{proof}

Let $d\in \mathbb N\cup \{+\infty\}$ be the degree of the cover $\rho_1$ and let $\alpha_1\colon \pi_1(M\setminus N_1)\to {\rm Sym}[d]$ be the monodromy map.  Then applying (1) from Lemma \ref{lem-WDomExt} with $\rho_1=\alpha_1$  gives a homomorphism $\rho_2\colon  \pi_1(M\setminus N_2)\to {\rm Sym}[d]$ that is the monodromy map of the desired cover.

\end{proof}
\begin{lem}\label{lem-WdominationTrans}
Let $N_i,i=1,\ldots, k$ be closed subsets of $M$ and let $W_i,i=1,\ldots, k-1$ be closed such that $N_i,N_{i+1}\subset W_i, i=1,\ldots, k-1.$ Suppose $N_i\xRightarrow{W_i}N_{i+1}, i=1,\ldots k-1$. Then, 
$$N_1\xRightarrow{W} N_k, \quad  W:=\bigcup_{i=1}^{k-1}W_i.$$
\end{lem}
\begin{proof}
We first prove the case $k=3.$ 
The diagram below commutes. 
\[\begin{tikzcd}[ampersand replacement=\&]
	\&\& {\pi_1(M\setminus N_1)} \\
	\& {\pi_1(M\setminus W_1)} \\
	{\pi_1(M\setminus [W_1\cup W_2])} \&\& {\pi_1(M\setminus N_2)} \& \\
	\& {\pi_1(M\setminus W_2)} \\
	\&\& {\pi_1(M\setminus N_3)}
	\arrow["{\iota_{1,1}}"', from=2-2, to=1-3]
	\arrow["{\iota_{1,2}}", from=2-2, to=3-3]
	\arrow["{\eta_1}", from=3-1, to=2-2]
	\arrow["{\eta_2}"', from=3-1, to=4-2]
	\arrow["{\kappa_1}"', dashed, from=3-3, to=1-3]
	\arrow["{\iota_{2,1}}"', from=4-2, to=3-3]
	\arrow["{\iota_{2,2}}", from=4-2, to=5-3]
	\arrow["{\kappa_2}"', dashed, from=5-3, to=3-3]
\end{tikzcd}\]

Using Definition \ref{def-Wdomination} we construct the maps $\kappa_1,\kappa_2$ making the above diagram commute. Then, $\iota_{1,1}\circ \eta_1=(\kappa_1\circ \kappa_2)\circ (\iota_{2,2}\circ \eta_2)$.
This proves $N_1\xRightarrow{W_1\cup W_2} N_3$ and concludes the proof for $k=3$. To prove the general case, we use this step inductively to $N_1, N_j, N_{j+1}$ and $\bigcup_{i=1}^{j-1} W_i, W_j$, starting from $j=2$ and repeating until $j=k-1.$
\end{proof}

\begin{lem}\label{lem-FlowLocus}
Let $\Phi \colon M\times \mathbb R\to M$ be a continuous action of $\mathbb R$. Let  $N\subseteq M$ be a compact subset with an open neighborhood $U$ such that $\pi_1(M\setminus U)\to \pi_1(M\setminus N)$ is surjective. Let $N_t:=\Phi(N,t)$ and $W_{s,t}:=\Phi(N,[s,t]), s<t$. Then $N_s\xRightarrow{W_{s,t}}N_t$ for all $s<t\in\mathbb R.$ 

\end{lem}
\begin{proof}
We will first show that the statement holds for $s=0$ and small enough $t>0.$ Put $U_t:=\Phi(U,t).$ Let $\varepsilon>0$ be such that $N_{t}\in U$ for all $t\in [0,\varepsilon].$ Then the isomorphism $\pi_1(M\setminus U_{t})\to \pi_1(M\setminus N_t)$ factors through $\pi_1(M\setminus W_{0,t})\to \pi_1(M\setminus N_t)$ for all $0\leq t\leq \varepsilon$, which means the latter map must be surjective. Let $\alpha\colon [0,1]\to M\setminus W_{0,t}$ be a closed loop such that $[\alpha]\in \ker \left[\pi_1(M\setminus W_{0,t})\to \pi_1(M\setminus N_t)\right].$ Then $\alpha$ is freely homotopic to $\alpha_{-t}:=\Phi(\alpha,-t)$ in $M\setminus N_0$, because $\Phi(\alpha,-s)\subset M\setminus W_{-s,t-s}$ is disjoint from $N_0$ for all $s\in [0,t]$. The fact that $\alpha$ is contractible in $M\setminus N_t$ means that $\alpha_{-t}$ is contractible in $M\setminus N_0$. Being contractible is preserved by free homotopy, so $[\alpha]\in \ker \left[\pi_1(M\setminus W_{0,t})\to \pi_1(M\setminus N_0)\right].$ By Lemma \ref{lem-WDomExt} (2), $N_0\xRightarrow{W_{0,t}}N_t$ for all $0\leq t\leq \varepsilon.$ 
Since $\Phi(\bullet, s)$ is a diffeomorphism, we immediately deduce that $N_s\xRightarrow{W_{s,t}}N_{t}$ for all $s\leq t\leq s+\varepsilon.$ Now, we finish the proof using Lemma \ref{lem-WdominationTrans} and the fact that $W_{s,r}\cup W_{r,t}=W_{s,t}$ for all $s\leq r\leq t.$ 
\end{proof}

\begin{lem}\label{lem-HomeoCriterion}
Let $N_1,N_2\subset W$ be closed subsets of $M$. Suppose that there is a homeomorphism $(M\setminus W)\to (M\setminus N_2)$ which is homotopic to the inclusion map. Then $N_1\xRightarrow{W} N_2.$
\end{lem}
\begin{proof}
In this case $\iota_2$ coincides with map induced by the homeomorphism, so $\ker \iota_2=\{1\}\subset \ker \iota_1$, and $\iota_2$ is surjective, hence $N_1\xRightarrow{W} N_2$ by Lemma \ref{lem-WDomExt} (2). 
\end{proof}

\begin{lem}\label{lem-BiggerW}
Let $N_1,N_2\subset W\subset W'$ be closed subsets of $M$. Suppose 
that $N_1\xRightarrow{W} N_2.$ Then, $N_1\xRightarrow{W'}N_2.$
\end{lem}
\begin{proof}


By the diagram below we can use the map $\kappa$ provided by $N_1\xRightarrow{W} N_2$ to show $N_1\xRightarrow{W'}N_2.$ 

\[\begin{tikzcd}
	&&& {\pi_1(M-N_1)} \\
	{\pi_1(M - W') } && {\pi_1(M-W)} \\
	&&& {\pi_1(M-N_2)}
	\arrow[from=2-1, to=1-4]
	\arrow[from=2-1, to=2-3]
	\arrow[from=2-1, to=3-4]
	\arrow[from=2-3, to=1-4]
	\arrow[from=2-3, to=3-4]
	\arrow["\kappa", dashed, from=3-4, to=1-4]
\end{tikzcd}\]
\end{proof}
From this point onward let us assume that $M$ is a Riemannian $n$-manifold equipped with the structure of a simplicial complex (i.e. a triangulation). We recall that for $0\leq m\leq n$ we write $M^{(m)}$ for the $m$-skeleton, $M^{[m]}$ for the set of $m$-simplices of $M$ and $M^{[\cdot]}$ for the set of all simplices of $M$. Our simplices are closed by default. 
We start with a few simple observations.
\begin{lem}\label{lem-CodimDomination} Let $M$ be a Riemannian $n$-manifold equipped with a triangulation.
\begin{enumerate}
\item Let $N\subset M^{(n-2)}$ be a closed subcomplex and let $Q\subset M$ be a closed subset of $M^{(n-3)}$, Then $(N\cup Q)\xRightarrow{(N\cup Q)} N$. In particular, setting $N=\emptyset$ we get $Q\xRightarrow{Q}\emptyset.$
\item Let $N_1,N_2\subset W \subset M$, $Q\subset M^{(n-2)}$ be closed sets such that $N_1\xRightarrow{W\cup Q}N_2.$ Then,
$N_1\xRightarrow{W}N_2.$
\item Let $N\subset M^{(n-1)}$ be a closed sub-complex. Then $N\xRightarrow{N}N^{(n-2)}.$
\end{enumerate}
\end{lem}
\begin{proof}
\begin{enumerate}
\item The map $\pi_1(M\setminus (N\cup Q))\to \pi_1(M\setminus N)$ is an isomorphism\footnote{This is well known property of codimension $3$ sets, see \cite[Thm 2.3, p. 146]{Godbillion} for a proof in the case $N,M,Q$ are smooth submanifolds.} so the Definition \ref{def-Wdomination} is trivially satisfied.

\item  

If $N_1\xRightarrow{W}N_2$, then the same map $\kappa$ works for  $N_1\xRightarrow{W\cup Q}N_2$ and $N_1\xRightarrow{W}N_2$.  This follows from the diagram below and the fact that $\pi_1(M-(W\cup Q))\to \pi_1(M-W)$ is surjective\footnote{See \cite[Thm 2.3, p. 146]{Godbillion} for the proof of surjectivity in case $M,W,Q$ are smooth manifolds.}:

\[\begin{tikzcd}
	&&& {\pi_1(M \setminus N_1)} \\
	{\pi_1(M\setminus(W\cup Q))} && {\pi_1(M \setminus W)} \\
	&&& {\pi_1(M\setminus N_2)}
	\arrow[from=2-1, to=1-4]
	\arrow[two heads, from=2-1, to=2-3]
	\arrow[from=2-1, to=3-4]
	\arrow[from=2-3, to=1-4]
	\arrow[from=2-3, to=3-4]
	\arrow["\kappa"', from=3-4, to=1-4]
\end{tikzcd}\]

\item We need to show that the map $\pi_1(M\setminus N)\to \pi_1(M\setminus N^{(n-2)})$ induced by inclusion is injective. One can prove this by gluing the interiors of $(n-1)$-cells of $N$ back to $M\setminus N$, one by one, and using van Kampen's theorem to justify that these operations induce injective maps on fundamental groups at each step. We indicate how to do this for one such cell $\sigma\in N^{[n-1]}$. Let $U$ be an open neighborhood of $\sigma$ such that $\sigma$ is a deformation retract of $U$ and $U$ is contractible. The van Kampen theorem for the sets $M\setminus\sigma$ and $U\setminus \sigma^{(n-2)}$ yields an isomorphism \[\pi_1(M\setminus\sigma^{(n-2)})\cong \pi_1(M\setminus\sigma)\ast_{\pi_1(U\setminus \sigma)}\pi_1(U\setminus \sigma^{(n-2)}).\] The group $\pi_1(U\setminus \sigma)\cong \pi_1(\mathbb S^{n-1})$ is trivial if $n\geq 3$ so the map $\pi_1(M\setminus\sigma)\to \pi_1(M\setminus\sigma^{(n-2)})$ is injective. For $n=2$, once checks by hand that $\pi_1(U\setminus \sigma)\to \pi_1(U\setminus \sigma^{(n-2)})$ is injective so $\pi_1(M\setminus\sigma)\to \pi_1(M\setminus\sigma^{(n-2)})$ is still injective.
\end{enumerate}
\end{proof}

Our next step is to show that for any $k$-varifold $S\subset M, k\leq n$, the Federer-Fleming deformation yields a pair $S\xRightarrow{W}\FF(S)$, where $W$ is the image of the homotopy between $S$ and $\FF(S)$ provided by the Federer-Fleming deformation. 

The following lemma will serve as the ``building block" for the subsequent constructions.  Let $\Delta$ be the standard closed $n$-simplex in $\mathbb{R}^n$.  Let $F$ be a closed subset of $\Delta$, and let $p \notin F$ be a point in $\Delta$.  For each $p' \in F$, let $\ell_{p'}$ be the closed line segment joining $p'$ to its image in $\partial \Delta$ under radial projection at $p$.   Let
\[
K_F:= \cup_{p' \in F} \ell_{p'}.  
\]
 \begin{figure}[h!]
		\includegraphics[width=100mm]{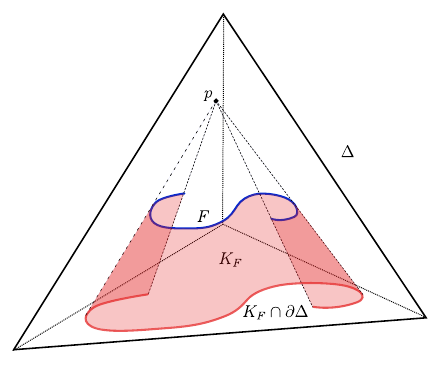}
		\caption{Sets $K_F$ (red) is the ``shadow" of $F$ (blue) seen from $p$.} \label{fig-KF}
	\end{figure}
\begin{lem}\label{lem-TTlemma}
For any closed proper subset $F$ of $\Delta$, let $K_F \subset \Delta$ be defined as above.  Then, there exists a map $\Phi\colon \Delta\times [0,1]\to \Delta$ such that $\Phi(x,0)=x, x\in \Delta$,  $\Phi(x,t)=x, x\in\partial \Delta, t\in [0,1]$, $\Phi(\bullet, 1)\colon \Delta\setminus K_F\to \Delta\setminus (K_F\cap \partial \Delta)$ is a homeomorphism and $\Phi(\bullet, t)\colon \Delta\to\Delta$ is a homeomorphism for all $t\in [0,1).$
\end{lem}
\begin{proof}

 \begin{figure}[h!]
		\includegraphics[width=130mm]{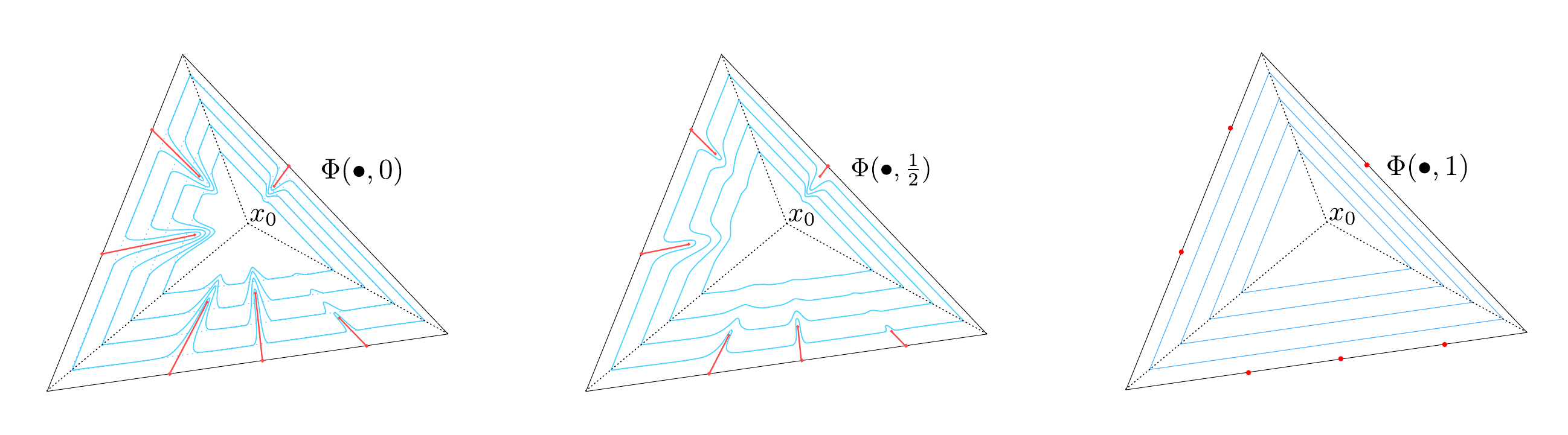}
		\caption{The maps $\Phi(\bullet,0), \Phi(\bullet,1/2)$ and $\Phi(\bullet,1)$. Blue lines are mapped to blue lines.} \label{fig-TT}
	\end{figure}
We identify $\Delta$ with $\partial \Delta\times [0,1]/\sim$ where $(x,0)\sim (y,0), x,y\in\partial\Delta$ via the map sending $x\mapsto (x',t)$, where $x'\in\partial\Delta$ is the unique point of the boundary such that $x$ lies on the line segment connecting $x'$ to $p$ and $t:=\frac{d(x,p)}{d(x',p)}.$ The point $p$ is sent to the class of $(p,t), t\in [0,1].$ Let $f\colon \partial \Delta\times [0,1]$ be the inverse map, that is $f(x',t)$ is the unique point of $\Delta$ which lies on the segment connecting $p$ to $x'$ and $d(p, f(x',t))=t d(p,x').$ Let $\ell_x=\min_{t\in [0,1]}f(x,t)\in K_F.$ We note that the definition of $K_F$ implies that $f(x,t)\in K_F$ if and only if $t\geq \ell_x.$

For any point $x\in \partial \Delta, t, s\in [0,1]$ and $s<1$ or $t<\ell_x$ let 
$$\phi(x,t,s):=\int_0^t \min \{d(f(x,t),K_F)^{-1},s(1-s)^{-1}\}dt.$$
This function is strictly increasing in $t$ for $s\in (0,1)$ and constant and equal to $0$ if $s=0$. For $s=1$ it is strictly increasing for $t\in [0,\ell_x)$ and tends to $+\infty$ as $t\to \ell_x.$
The function
\[g(x,t,s):=\begin{cases}
1-e^{-\phi(x,t,s)}+te^{-\phi(x,1,s)} & s<1 \text{ or } t<\ell_x\\
1 & \text{ otherwise,}\\
\end{cases}\]
is continuous on $\partial\Delta\times [0,1]^2$. We note $g(x,t,0)=t$ and $g(x,1,s)=1$. In addition, for $s\in [0,1),x\in\partial\Delta$ it defines a monotone increasing bijection $g(x,\bullet,s)\colon [0,1]\to [0,1]$, for $s=1$ it defines a monotone increasing bijection $g(x,\bullet,s)\colon [0,\ell_x]\to [0,1]$ and finally $g(x,t,1)=1$ for $t\in [\ell_x,1].$ 
We put $\Phi(f(x,t),s):= f(x, g(x,t,s)).$
The function $\Phi\colon \Delta\times [0,1]\to\Delta$ satisfies $\Phi(f(x,t),0)=f(x,t),t\in [0,1]$ and $\Phi(f(x,1),s)=f(x,1), s\in[0,1]$. For $s\in [0,1)$ it defines a continuous bijection on $\Delta,$ hence a homeomorphism. For $s=1$, it restricts to a homeomorphism $\Phi(\bullet,1)\colon \Delta\setminus K_F\to \Delta\setminus (K_F\cap \partial \Delta)$ (see Figure \ref{fig-TT})

\end{proof}

\begin{lem}\label{lem-SimplexExtension}
Let $\Delta$ be an affine simplex in $\mathbb R^n$. Let $n\geq l\geq m.$ Write $\Delta^{(m)}$ for the $m$-skeleton of $\Delta$. Let $W\subset \Delta^{(m)}$ be a closed set and suppose we have a continuous function 
 $\Phi^{(\ell)}\colon \Delta^{(l)}\times [0,1]\to \Delta^{(l)}$ satisfying 
\begin{enumerate}[label=(A\arabic*)]
\item $\Phi^{(l)}(x,0)=x$ for all $x\in \Delta$,
\item $\Phi^{(l)}(x,t)=x$ for all $x\in \Delta^{(m-1)}$ and $t\in [0,1],$
\item $\Phi^{(l)}(\bullet,t)\colon \Delta^{(l)}\to \Delta^{(l)}$ is a homeomorphism for all $t\in [0,1),$ 
\item $\Phi^{(l)}(\bullet,1)$ restricts to a homeomorphism $\Delta^{(l)}\setminus W\to \Delta^{(l)} \setminus (W\cap \Delta^{(m-1)}).$
\item $\Phi^{(l)}(\bullet,t)$ maps each $l$-cell to itself. 
\end{enumerate}
Then, for any $n\geq k\geq l$ there is an extension $\Phi^{(k)}$ satisfying the same conditions (with $l$ replaced by $k$).
\end{lem}
\begin{proof}
We prove the statement inductively on $k$, with the case $k=l$ serving as the base case. Suppose we have constructed the map $\Phi^{k}\colon \Delta^{(k)}\times [0,1]\to \Delta^{(k)}.$ Let $\sigma$ be a $(k+1)$-cell and choose a point $p$ in the interior of $\sigma$. 

Let $f\colon \partial \sigma \times [0,1]\to \sigma$ be the map that sends $(x,t)$ to the unique point $f(x,t)$ on the line segment joining $p$ and $x$ such that $d(p,f(x,t))=t d(p,x).$ Define \[\Phi^{(k+1)}(f(x,t),s):=f(\Phi^{(k)}(x,ts),t).\]
We note $\Phi^{(k+1)}(f(x,0),s)=f(\Phi^{(k)}(x,0),0)=p$ so it is a well defined map on $\partial \sigma\times [0,1]$. 
We verify that it indeed extends $\Phi^{(k)}$: for $x\in \partial\sigma$, 
\[ \Phi^{(k+1)}(x,s)=\Phi^{(k+1)}(f(x,1),s)=\Phi^{(k)}(x,s).\]
The maps $\Phi^{(k+1)}$, defined so far for each $(k+1)$-cell $\sigma$, agree on the boundaries of adjacent $(k+1)$-cells so they give a continuous map $\Phi^{(k+1)}\colon \Delta^{(k+1)}\times [0,1]\to \Delta^{(k+1)}.$ Now we check that it verifies the conditions (A1)-(A4), (A5) is guaranteed because we define it on each $k$-simplex $\sigma$ separately. 
\begin{align*}\Phi^{(k+1)}(f(x,t),0)=&f(\Phi^{(k)}(x,0),t)= f(x,t), \quad \text{(A1)}\\
\Phi^{(k+1)}(f(x,1),s)=&f(\Phi^{(k)}(x,1),1)=f(x,1), x\in \Delta^{(m-1)}\cap \partial \sigma, \quad \text{(A2)}\\
\end{align*}
For $s\in[0,1)$ the map $\Phi^{(k+1)}$ is continuous bijection from $f(\partial \sigma, t)$ to itself. These maps glue to a continuous bijection $\sigma\to\sigma,$ they agree on the boundaries of $(k+1)$-simplices so they glue to a continuous bijection  $\Delta^{(k+1)}\to\Delta^{(k+1)}$, which verifies (A3). 

It remains to verify (A4). We do it simplex by simplex, since we already know the condition is satisfied on the $(k-1)$-skeleton. First we check it is a bijection. For $s=1$ and $t<1$ the map $\Phi^{(k+1)}$ is a continuous bijection from $f(\partial \sigma, t)$ to itself. These maps glue to a continuous bijection from the set $\{f(\partial \sigma, t)\mid t\in [0,1-\varepsilon]\}$ to itself, for any $\varepsilon>0$. We deduce that it must be a homeomorphism on $\sigma\setminus \partial\sigma$. The map is continuous on $f(\partial \sigma, 1)=\partial \sigma,$ where it coincides with $\Phi^{(k)},$ so it restricts to a homeomorphism $\partial \sigma \setminus W\to \partial \sigma \setminus (W\cap \Delta^{(m-1)}).$ These observations imply that $\Phi^{(k+1)}(\bullet,1)\colon \sigma \setminus W\to \sigma \setminus (W\cap \Delta^{(m-1)})$ is a continuous bijection. 
To show that it is a homeomorphism we prove the continuity of the inverse map. Suppose $x_n,y_n\in \partial \sigma$ and $u_n,v_n\in  [0,1]$ are sequences such that $f(x_n,u_n),f(y_n,v_n)\in \sigma\setminus W,n\in\mathbb N,$ and the sequences  
\begin{align*}\Phi^{(k+1)}(f(x_n,u_n),1)&=f(\Phi^{(k)}(x_n,u_n),u_n), \\ \Phi^{(k+1)}(f(y_n,v_n),1)&=f(\Phi^{(k)}(y_n,v_n),v_n)\end{align*} converge to the same point $f(z,w)\in \sigma\setminus (W\cap \Delta^{(m-1)}), z\in \partial\sigma, w\in  [0,1].$ We already know that $\Phi^{(k+1)}$ restricts to a homeomorphism on $\sigma\setminus \partial\sigma$ so we can assume $w\neq 0.$ Clearly $u_n\to w$ and $v_n\to w$ as $n\to\infty$. Therefore, the sequences $\Phi^{(k)}(x_n,w),\Phi^{(k)}(y_n,w)$ also converge to $z$. We deduce that $x_n\to z$ and $y_n\to z$ either by hypothesis (A3) for $\Phi^{(k)}$ when $0<w<1$, or by (A4) for $\Phi^{(k)}$ when $w=1.$ In this way we showed that the inverse map 
$\Phi^{(k+1)}(\bullet,1)^{-1}\colon \sigma \setminus (W\cap \Delta^{(m-1)})\to \sigma\setminus W$ is continuous, which confirms (A4).
\end{proof}
\begin{lem}\label{lem-FFmovinglocus}
Let $S\subset M$ be a closed $k$-varifold, $k\leq n-1$. Let $\FF\colon S\to M$ be the Federer-Fleming deformation and let $h\colon S\times [0,1]\to M$ be the homotopy connecting the identity and the Federer-Fleming map constructed in Section \ref{sec-FF}. Then,
$$ S\xRightarrow{h(S,[0,1])} \FF(S).$$
\end{lem}
\begin{proof}
Let $\FF^m$ and $h^m$ be the maps defined in Section \ref{sec-FF}. We put $S_i:=\FF^{n-i}, W_i:=h^{n-i}(S_i,[0,1]), i=0,1,\ldots,n-k+1$ with $S_0:=S$. 
We'll write $M^{(m)}$ for the $m$-skeleton of the triangulation of $M$ used to define the Federer-Fleming deformation. 
Recall that $h^{n-i}$ is the homotopy between $S_i$ and $S_{i+1}$ (defined in Section \ref{sec-FF}) and $h$ is obtained by concatenating $h^n,h^{n-1},\ldots, h^k,$ so \[W:= h(S,[0,1])=W_0\cup W_1\cup \ldots\cup W_{n-k}.\] We shall prove that $S_i\xRightarrow{W_i}S_{i+1}.$ Then by Lemma \ref{lem-WdominationTrans}, $S=S_0\xRightarrow{W} S_{n-k+1}=\FF(S)$ which is what we are after. 

The proof of the statement  $S_i\xRightarrow{W_i}S_{i+1}$ will look different for $i=0, 1\leq i\leq n-k-1$ and $i=n-k$ but the rough outline is the same. Our goal is to construct a homotopy $\Phi\colon M\times [0,1]\to M$ between the identity map and a map which restricts to a homeomorphism $M\setminus W_i\to M\setminus S_{i+1}$. This will enable us to conclude $S_i\xRightarrow{W_i}S_{i+1}$ by Lemma \ref{lem-HomeoCriterion}. First we use Lemma \ref{lem-TTlemma} to construct the desired homotopy on the $(n-i)$-skeleton of $M$, basically by applying it on each $(n-i)$-simplex separately. Next we need to carefully apply Lemma \ref{lem-SimplexExtension} several times to extend it to a homotopy on the whole $M$. 

Let $i=0$. For each $n$-simplex $\sigma\in M^{[n]}$ the set $S\cap \sigma$ is closed. Let $x_0\in \sigma$ be the point used to define the $\FF^{n-1}$ on $\sigma$. Then $\FF^{n-1}(S)\cap \sigma=p_{x_0}(S\cap \sigma)$ and $W_0\cap\sigma =h_{x_0}(S\cap\sigma, [0,1]).$ Use the map $\iota_\sigma$ to identify $\sigma$ with an affine simplex $\Delta$ and apply Lemma \ref{lem-TTlemma} to $F=\iota_\sigma(S\cap \sigma)$. Then $K_F=\iota_\sigma(W_0\cap \sigma)$ and the Lemma \ref{lem-TTlemma} yields a homotopy $\Phi_0^{(n)}\colon \sigma \times [0,1]\to \sigma$ such that $\Phi_0^{(n)}(\bullet, 0)$ is the identity, $\Phi_0^{(n)}(x,t)=x$ for all $x\in \partial\sigma, t\in [0,1]$ and $\Phi_0^{(n)}(\bullet, 1)$ restricts to a homeomorphism $\Phi_0^{(n)}(\bullet, 1)\colon \sigma \setminus W_0\to \sigma \setminus S_1.$ These maps agree on the boundaries of $n$-simplices (where they are identity maps) so they glue to a global homotopy $\Phi_0^{(n)}\colon M\times [0,1]\to M$ between the identity and a map that restricts to a homeomorphism $M\setminus W_0\to M\setminus S_1.$ This proves the case $i=0.$

Let $0<i\leq n-k$. By construction, the image of $\FF^{n-i}$ is contained in $M^{(n-i)}$. As in the case $i=0$ for each $(n-i)$-simplex $\sigma\in M^{[n-i]}$ apply Lemma \ref{lem-TTlemma} to get a map $\Phi_i^{(n-i)}\colon \sigma \times [0,1]\to\sigma$ with the following properties. The map $\Phi_i^{(n-i)}(\bullet,0)$ is the identity on $\sigma$,  $\Phi_i^{(n-i)}(\bullet,t)$ is the identity on $\partial\sigma$ for all $t\in [0,1]$ and $\Phi_i^{(n-i)}(\bullet,1)$ restricts to a homeomorphism $\sigma \setminus W_i\to \sigma \setminus S_{i+1}.$ These maps all agree on the boundaries of $(n-i)$-simplices, where they are identity maps, so they glue to a map 
$\Phi_i^{(n-i)}\colon M^{(n-i)}\times [0,1]\to M,$ satisfying the axioms (A1)-(A5) of the Lemma \ref{lem-SimplexExtension} on each $(n-i+1)$-simplex of $M^{(n-i+1)}$. We successively apply Lemma \ref{lem-SimplexExtension} to construct a sequence of extensions $\Phi_i^{(m)}\colon M^{(m)}\times [0,1]\to M^{(m)}, m=n-i+1,\ldots, n.$ The final map $\Phi_i^{(n)}$ is a homotopy between the identity on $M$ and a map that restricts to a homeomorphism $M\setminus W_i\to M\setminus S_{i+1}$. By Lemma \ref{lem-HomeoCriterion}, $S_i\xRightarrow{W_i}S_{i+1}.$

For $i=n-k$ the proof is the same but we have to define $\Phi_i^{(n-i)}$ differently for the $(n-i)=k$-cells which are left unaffected by the last step of the Federer-Fleming deformation. For those cells $\sigma,$ we have $W_i\cap \sigma= S_{i+1}\cap \sigma \subseteq \partial\sigma$. We put $\Phi_i^{(n-i)}(x,t)=x, x\in\sigma, t\in [0,1].$ Otherwise, the contruction of $\Phi_i^{(m)}, m=n-i,\ldots,n$ proceeds as in the previous step. 
\end{proof}

\begin{lem}\label{lem-FFWdomination} Suppose $M$ is equipped with a triangulation. Let $N_1,N_2\subset M^{(n-2)}$ be closed subsets of $M$ and let $W$ be a rectifiable $(n-1)$-varifold, such that $N_1\xRightarrow{W} N_2$. Then $N_1\xRightarrow{\FF(W)}N_2.$
\end{lem}
\begin{proof} 
We prove the lemma in two steps, first that $N_1\xRightarrow{ \FF^{n-1}(W)}N_2$ and then $N_1\xRightarrow{\FF^{n-2}(W)}M_2.$
Let $Q_1:=h^n(W,[0,1])$ be the image of the homotopy connecting $W$ and $\FF^{n-1}(W)$. 
By Lemma \ref{lem-BiggerW} $N_1\xRightarrow{Q_1}N_2.$ 

Using Lemma \ref{lem-TTlemma} on each $n$-simplex, as in the proof of Lemma \ref{lem-FFmovinglocus}, we construct a homotopy $\Phi_1\colon M\times [0,1]\to M$ such that $\Phi_1(\bullet,0)$ is the identity map, $\Phi_1(\bullet,t)$ is the identity on $M^{(n-1)}$ and $\Phi_1(\bullet,1)$ restricts to a homeomorphism $M\setminus Q_1\to M\setminus \FF^{n-1}(W).$ The following diagram commutes.

\[\begin{tikzcd}
	{\pi_1(M\setminus N_1)} \\
	& {\pi_1(M\setminus Q_1)} && {\pi_1(M \setminus \FF^{n-1}(W))} \\
	{\pi_1(M \setminus N_2)}
	\arrow["{\iota_1}", from=2-2, to=1-1]
	\arrow["{\Phi_1(\cdot,1)}", from=2-2, to=2-4]
	\arrow["{\iota_2}"', from=2-2, to=3-1]
	\arrow["{\iota_1'}"', from=2-4, to=1-1]
	\arrow["{\iota_2'}", from=2-4, to=3-1]
	\arrow["\kappa", from=3-1, to=1-1]
\end{tikzcd}\]
The map $\Phi_1(\bullet,1)$ is homotopic to identity on $M$ and $\Phi_1(\bullet,1)$ is the identity on $N_1,N_2$. Therefore the diagram above implies that the same $\kappa$ works to show $N_1\xRightarrow{\FF^{n-1}(W)}N_2.$ 

The proof that $N_1\xRightarrow{\FF^{n-2}(W)}N_2$ follows the same outline. Put \[Q_2:=h^{n-1}(\FF^{n-1}(W),[0,1])\] and note that $\FF^{n-1}(W)\subset Q_2$ implies $N_1\xRightarrow{Q_2}N_2.$ Using Lemma \ref{lem-TTlemma} as in the proof of Lemma \ref{lem-FFmovinglocus}, we construct a new homotopy $\Phi_2\colon M\times [0,1]\to M$ such that $\Phi_2(x,0)=x, x\in M$, $\Phi_2(x,t)=x, x\in M^{(n-2)}$ and $\Phi_2(\bullet,1)$ restricts to a homeomorphism $M\setminus Q_2\to M\setminus \FF^{n-2}(W)=M\setminus \FF(W).$ Using a diagram similar to the one above we show $N_1\xRightarrow{\FF(W)}N_2$ and conclude the proof using Lemma \ref{lem-WdominationTrans}.

\end{proof}

\subsection{Riemannian cover stability} \label{subsection:riemanniancoverstability}
Here we will prove Theorem \ref{mt-branched}. Recall that $M=\Gamma\bs \mathbb H_\mathbb O^2$ is a compact octonionic hyperbolic manifold with injectivity radius at least $\varepsilon$ and $\alpha$ is the constant from Theorem \ref{mt-collapse}, implicitly dependent on $\varepsilon$ but not on $\Gamma$ or the varifold $S$. The manifold $M$ is equipped with an $(r,\delta)$-uniform triangulation $\mathcal T$ (see Definition \ref{def-UnifTrian} and Theorem \ref{thm-UTboundedgeo}) which allows us to perform Federer-Fleming deformations on subsets of $M$. 

\begin{proof}[Proof of Theorem \ref{mt-branched}]
Our condition on the $14$-dimensional volume of $N$ needs to be slightly stronger than in Theorem \ref{mt-collapse}, this is for technical reasons.
We suppose that $\mathcal H^{14}(N)\leq c_0^{-1}\alpha,$ where $c_0$ is the constant from Lemma \ref{lem-FFproperties}.  We say that a compact subset $E\subset M$ is \emph{usable} if there is an open neighborhood $U$ containing $E$ such that $\pi_1(M\setminus U)\to \pi_1(M\setminus E)$ is surjective.

Let $N_0:=N$. We will construct a sequence of closed subsets $N_i$ and $W_i$ such that 
\begin{enumerate}[label=(C\arabic*)]
\item $N_i\xRightarrow{W_i}N_{i+1}$, $\mathcal{H}^{15}(W_i)\leq c_2 (\mathcal H^{14}(N_i)-\mathcal H^{14}(N_{i+1})),$
\item $\mathcal H^{14}(N_{i+1})\leq \min\{\frac{1}{2}\mathcal{H}^{14}(N_i), \mathcal{H}^{14}(N_i)-\tau_0\}$ for $i\geq 1,$
\item $N_i$ is usable for $i\geq 1$. 
\end{enumerate}
Moreover, the sequence eventually ends at an empty set. 

The construction uses two types of steps. We will later describe how we choose which step to apply. 
\begin{enumerate}[label=(S\arabic*)]
\item The first type of step is a ``tidy Federer-Fleming deformation". Suppose $N_i$ is already constructed, then we put $N_{i+1}':=\FF(N_i)$ and let $W_i'=h(N_i, [0,1])$, where $h$ is the homotopy afforded by the Federer-Fleming deformation. We have $N_i\xRightarrow{W_i'}N_{i+1}'$, $\mathcal H^{15}(W_i')\leq c_1 \mathcal H^{14}(N_i)$ by Lemmas \ref{lem-FFmovinglocus} and \ref{lem-FFproperties} respectively. We define $N_{i+1}$ as the union of the $14$-simplices contained in $N_{i+1}'$ and $Q:=N_{i+1}'\setminus N_{i+1}$. By Lemma \ref{lem-FFcodim3remainder}, $Q\subset M^{(13)}$ and $N_{i+1}$ is a closed subcomplex of $M^{(14)}.$ By Lemma \ref{lem-CodimDomination}, $N_{i+1}'\xRightarrow{N_{i+1}'} N_{i+1}$. Put $W_i:=W_i'\cup N_{i+1}'$. By Lemma \ref{lem-WdominationTrans}, $N_i\xRightarrow{W_i}N_{i+1}$ and $\mathcal H^{15}(W_i)=\mathcal H^{15}(W_i')\leq c_1 \mathcal H^{14}(N_i).$ Moreover, by Lemma \ref{lem-FFproperties} $\mathcal H^{14}(N_{i+1})\leq c_0\mathcal H^{14}(N_i)$, unconditionally on $N_i$. 

\item The second type of step will be a flow. In this case we assume that $N_i$ is usable, and that $N_i$ contains a $t$-heavy point for $t\leq r_1=r_1(\varepsilon).$ Arguing in exactly the same way as in the proof of Theorem \ref{mt-collapse}, we produce a smooth flow $\Phi\colon M\times \mathbb R\to M$ and a time $\nu_0>0$ such that $N_{i+1}:=\mathcal H^{14}(\Phi(N_i,\nu_0))\leq \mathcal H^{14}(N_i)-\tau_0$, where $\tau_0=\tau_0(\varepsilon)>0$ is the constant from Lemma \ref{lem-flowsqueezing}. By the same lemma, $W_i=\Phi(N_i,[0,\nu_0])$ satisfies $\mathcal H^{15}(W_i)\leq c_1 (\mathcal H^{14}(N_i)-\mathcal H^{14}(N_{i+1})$. By Lemma \ref{lem-FlowLocus}, $N_i\xRightarrow{W_i}N_{i+1}$. Since $N_{i+1}$ is the image of $N_i$ by a diffeomorphism of $M$, we have that $N_{i+1}$ is usable.
\end{enumerate}

As the first step of the construction we always choose (S1). This has the important effect of rendering $N_1$ usable, at the cost of possibly increasing $14$-dimensional Hausdorff measure by a multiplicative constant $c_0$. Then $N_1$ is a usable subset of $M$, satisfying $\mathcal H^{14}(N_1)\leq \alpha \vol(M).$ All subsequent steps of the construction will decrease $\mathcal H^{14}(N_i)$. 

Using Lemma \ref{lem-heavypintdichotomy}, either $N_i$ has a $t$-heavy point, $t\leq r_1=r_1(\varepsilon)$ or $\mathcal H^{14}(\FF(N_i))\leq\frac{1}{2}(\mathcal H^{14}(N_i)).$ In the first case, apply step (S2). In the second we apply step (S1). We verify inductively that the properties (C1)-(C3) are satisfied, with $c_2:=\max \{c_1, 2c_0\}.$ In particular, after finitely many steps $\mathcal H^{14}(N_{i_0})\leq \eta$ where $\eta=\eta(\varepsilon)>0$ is the constant from Lemma \ref{lem-FFvanishing}. Then, the next step must be (S2) which produces $N_{i_0+1}=\emptyset.$

We obtain a sequence $N=N_0\xRightarrow{W_0}N_1\xRightarrow{W_1}\ldots N_{i_0}\xRightarrow{W_{i_0}}\emptyset.$ Put $W:=\bigcup_{i=0}^{i_0} W_i$. Then $N\xRightarrow{W}\emptyset$ and $$W\leq c_0\mathcal H^{14}(N)+\sum_{i=1}^{i_0} c_2(\mathcal H^{14}(N_{i+1})-\mathcal H^{14}(N_i))\leq (c_0+c_2)\mathcal H^{14}(N).$$
This concludes the proof by Lemma \ref{lem-MovingLocus}.

\end{proof}

\subsection{Combinatorial cover stability}\label{sec-ccstability} In this section we prove Theorem \ref{thm-CovStab}. We follow the exposition from \cite{dinur2022near}, with minor modifications in the notation. 

Let $Y$ be a finite $n$-dimensional simplicial complex. We shall write $Y^{(m)}$ for the $m$-skeleton of $Y$, $Y^{[m]}$ for the set of oriented $m$-simplices and $Y^{\{m\}}$ for the set of unoriented $m$-simplices of $Y$.  A simplicial $n$-complex $Y$ is \emph{pure} if every cell is contained in some $n$-simplex. In \cite{dinur2022near}, authors define a natural system of weights $c_Y\colon Y^{\{m\}}\to [0,1], m=0,\ldots,n$ such that $\sum_{\sigma\in Y^{\{m\}}}c_Y(\sigma)=1.$ For us, the only relevant property of these weights is that for simplicial complexes of bounded geometry (i.e. when $Y^{(1)}$ has uniformly bounded degree) the weights are at most a multiplicative constant away from uniform weights: there is a constant $A=A(\max \deg Y^{(1)})>0$ such that 
\[ A^{-1} |Y^{\{m\}}|^{-1}\leq c_Y(\sigma)\leq A|Y^{\{m\}}|^{-1}, \text{ for all } \sigma\in Y^{\{m\}}.\]
For any subset $\Sigma\subset Y^{\{m\}},$ we put $\|\Sigma\|:=\sum_{\sigma\in \Sigma}c_Y(\sigma).$

Let $H$ be a group. We write $C^0(Y,H)$ for the group of $H$ valued functions on $Y^{[0]}$ equipped with pointwise multiplication. We will call $C^0(Y,H)$ the $0$-cochains. Let 
\[C^1(Y,H):=\{\varphi\colon Y^{[1]}\to H\mid \varphi(u,v)=\varphi(v,u)^{-1}, (u,v)\in Y^{[1]}\}\] be the group of $1$-cochains with values in $H$. We define the coboundary operator $\delta\colon C^0(Y,H)\to C^1(Y,H)$
\[\delta\psi(u,v):=\psi(u)\psi(v)^{-1}, (u,v)\in Y^{[1]}.\]
Finally, for $\varphi\in C^1(Y,H)$ and $(u,v,w)\in Y^{[2]}$ put 
\[\delta\varphi(u,v,w):=[\varphi(u,v)\varphi(v,w)\varphi(w,u)]_H,\] where $[h]_H$ stands for the conjugacy class of $h\in H$. 
Then $\delta\varphi$ is a well defined function from $Y^{[2]}$ to the set of conjugacy classes of $H$, denoted $[H]$. The set of $1$-cocycles is defined as
\[Z^1(Y,H):=\{\varphi\in C^1(Y,H)\mid \delta\varphi=[1]_H\}.\]
We observe that $\delta\varphi\in Z^1(Y,H)$ for any $\varphi\in C^0(Y,H).$
The group of $0$-cochains acts on $1$-cochains by 
$$\psi\cdot\varphi(u,v):=\psi(u)\varphi(u,v)\psi(v)^{-1}, \quad (u,v)\in Y^{[1]},\psi\in C^0(Y,H), \varphi\in C^1(Y,h).$$
We have $\psi\cdot 1=\delta\psi,$ where $1$ stands for the trivial $1$-cochain. A simple computation shows that $\delta (\psi\cdot \varphi)= \psi \cdot \delta\varphi, \psi\in C^0(Y,H), \varphi\in C^1(Y,H).$ In particular, the set of $1$-cocycles is preserved by the action of $C^0(Y,H).$ Motivated by mathematical physics terminology we call two $1$-cochains $\varphi_1,\varphi_2\in C^1(Y,H)$ \emph{gauge equivalent} if there exists $\psi\in C^0(Y,H)$ such that $\varphi_2=\psi\cdot\varphi_1.$ The set of gauge equivalence classes of $1$-cocycles forms the first cohomology set 
\[ H^1(Y,H):=C^0(Y,H)\backslash Z^1(Y,H).\]
We are now in a position to define the non-abelian Cheeger constant $h_1(Y,H)$ \cite{dinur2022near}. For any function $\psi\colon Y^{[i]}\to H$ or $\psi\colon Y^{[i]}\to [H]$, where $[H]$ is the set of conjugacy classes of $H$,  let $\supp \psi=\{\sigma\in Y^{[i]}\mid \psi(\sigma)\neq 1\}.$ We put 
\[ \|\psi\|:=\frac{1}{2}\sum_{\sigma\in \supp \psi}c_Y(\sigma).\]
The \emph{cosystolic expansion} of $\varphi\in C^1(Y,H)\setminus Z^1(Y,H)$ is given by
\[h(\varphi):=\frac{\|\delta\varphi\|}{\min_{\psi\in Z^1(Y,H)}\|\varphi\psi^{-1}\|}.\]
Finally, the non-abelian Cheeger constant (named cosystolic expansion constant in \cite{dinur2022near}) of $Y$ is
\begin{equation}\label{eq-nacheeger} h_1(Y,H):=\min\{h(\varphi)\mid \varphi\in C^1(Y,H)\setminus Z^1(Y,H)\}.\end{equation}
Uniform lower bounds on the non-abelian Cheeger constants are known for $n$-simplices \cite{meshulam2013bounded}, for families of spherical buildings $A_3(\mathbb F_q)$ \cite{dinur2022near} and for quotients of certain Bruhat-Tits buildings \cite{Dikstein2023CoboundaryAC}. Theorem \ref{thm-CovStab} provides the first example of bounded geometry family of Riemannian manifolds with a positive non-abelian Cheeger constant.

There is a classical correspondence between finite covers of $Y$ of degree $d$ and cohomology classes in $H^1(Y, \Sym[d]).$ In \cite{dinur2022near}, it is shown that the inequality $h_1(Y,\Sym[d])>\varepsilon$ implies a cover stability statement similar to Theorem \ref{mt-branched}. We review the correspondence between cohomology classes and covers and indicate how one can derive cover stability theorems. Suppose $Y_1\rightarrow Y$ is degree $d$ cover. Fix an isomorphism between $Y_1^{[0]}$ and $Y^{[0]}\times [d]$, for $[d]=\{1,2,..,d\}$. We will write $[u,i], i=1,\ldots,d$ for the vertices of $Y_1$ lying above $u\in Y^{[0]}$. The cover $Y_1$ determines the monodromy cocycle $\varphi$ by the condition $\varphi(u,v)(i)=j$, where $j\in [d]$ is the unique index such that $([u,i],[v,j])\in Y_1^{[1]}.$ Conversely, every cocycle $\varphi\in Z^1(Y,\Sym[d])$ determines a unique degree $d$-cover. It is an easy computation to check that two cocycles determine isomorphic covers if and only if they are gauge equivalent, thus establishing the bijection between $d$-covers and $H^1(Y,\Sym[d]).$ 
Now, if $\varphi\in C^1(Y,\Sym[d])$ is a cochain and $N\subset Y$ is a pure $n$-subcomplex such that $\supp \delta\varphi\subset N^{[2]},$ then $\varphi$ determines a cocycle on $Y\setminus N, $ and hence a degree-$d$ cover of $Y\setminus N.$ We can therefore think of cochains as representing branched covers of $Y$, with branching locus containing $\supp \delta\varphi$. By definition of $h_1(Y,\Sym[d])$, we can find a cocycle $\psi\in Z^1(Y,\Sym[d])$ such that $\|\varphi\psi^{-1}\|\leq h_1(Y,\Sym[d])^{-1}\|\delta\varphi\|.$ Let $W$ be any pure subcomplex of $Y$ such that $\varphi=\psi$ outside of $W$. Then, the branched covers determined by $\varphi$ and the genuine cover determined by $\psi$ coincide outside of $W$. If $Y$ is bounded geometry, then the size of $W$ can be bounded linearly by the size of $N$ times $h_1(Y,\Sym[d])^{-1}.$ For details we refer to \cite{dinur2022near}. This concludes the introductory part of this section. 


\begin{lem}\label{lem-cocyclefilling}
Let $Y$ be a finite $n$-dimensional pure simplicial complex. Suppose $\varphi\in C^1(Y,H)$ and $N_1,N_2$ are closed subcomplexes of $Y$ such that $\supp \delta\varphi\cap (Y\setminus N_1)^{[2]}=\emptyset $ (i.e. $N_1$ contains $\supp \delta\varphi$ in its interior). Let $W$ be a closed subcomplex such that $N_1\xRightarrow{W}N_2.$ Then, there is a cochain $\psi\in C^1(Y,H)$ with $\supp \delta\psi\subset N_2^{[2]}$ and $\varphi=\psi$ outside of $W$.
\end{lem}
\begin{proof}
Since $\supp \delta\varphi\cap (Y\setminus N_1)^{[2]}=\emptyset$, the restriction $\varphi|_{Y\setminus N_1}$ is a cocycle. Let $\varphi^*\colon \pi_1(Y\setminus N_1)\to H$ be the monodromy representation. Since $N_1\xRightarrow{W}N_2,$ there is a representation $\psi^*\colon \pi_1(Y\setminus N_2)\to H,$ such that $\varphi^*=\psi^*$ on $\pi_1(Y\setminus W).$ We can therefore select a cochain $\psi\in C^1(Y,H)$ which restricts to a cocycle on $Y\setminus N_2$ and induces $\psi^*$ on $\pi_1(Y\setminus N_2)$. Both $\varphi,\psi$ determine the same monodromy representation on $\pi_1(Y\setminus W)$, so by replacing $\psi$ by gauge equivalent cochain we can arrange for $\varphi=\psi$ outside of $W$. Finally, observe that since $\psi$ is a cocycle on $Y\setminus N_2$, we have $\supp \delta \psi\subset N_2^{[2]}.$
\end{proof}

\begin{proof}[Proof of Theorem \ref{thm-CovStab}]
We will show that any $(r,\delta)$-uniform triangulation of a compact octonionic hyperbolic manifold $M$ with $\injrad M\geq \varepsilon$ and any group $H$ satisfies 
\[h_1(M,H)\geq h=h(r,\delta,\varepsilon)>0.\]
All the constants appearing in the argument below depend only on $r,\delta,\varepsilon.$ Let $\varphi\in C^1(M,H).$ We have 
\[|\supp \delta\varphi|\geq c_0\|\delta\varphi\|\cdot |M^{[2]}|.\] Let $N$ be the union of all closed $16$-simplices intersecting $\supp \delta\varphi$ non-trivially. Then $|N^{[\cdot]}|\leq c_1 \|\delta \varphi\| |M^{[\cdot]}|$ and $\varphi|_{M\setminus N}$ is a cocycle. 

There is a constant $\alpha_1>0$ such that for $\|\delta\varphi\|\leq \alpha_1$ the inequality $\mathcal H^{14}(N^{(14)})\leq \alpha \vol(M)$ holds,  where $\alpha$ is the constant from Theorem \ref{mt-branched}. Suppose that this is the case, i.e. $\|\delta\varphi\|\leq \alpha_1$. By Theorem \ref{mt-branched}, there is a rectifiable $15$-varifold $W\subset M$ such that $N^{(14)}\xRightarrow{W}\emptyset$ and \[\mathcal H^{15}(W)\leq c_2 \mathcal{H}^{14}(N)\leq c_3 \|\delta\varphi\| |M^{[\cdot]}|.\]
By Lemma \ref{lem-FFWdomination}, $N^{(14)}\xRightarrow{\FF(W)}\emptyset.$ Let $W_1$ be the union of all $(15)$-simplices $\sigma$ such that $\FF(W)\cap {\rm int}\,\sigma \neq \emptyset.$ Then $W_1\setminus \FF(W)\subset M^{(14)}$, so the map $\pi_1(M-(W_1 \cup FF(W)))\rightarrow \pi_1(M)$ is still surjective, and $N^{(14)} \xRightarrow{W_1\cup \FF(W)}\emptyset.$ By Lemma \ref{lem-CodimDomination} (2), $N\xRightarrow{W_1}\emptyset.$ Since all the $(15)$-cells $\sigma$ with $\FF(W)\cap {\rm int}\,\sigma \neq \emptyset$ satisfy $\mathcal H^{15}(\FF(W)\cap {\rm int}\,\sigma)=\mathcal H^{15}(\sigma)\geq \eta=\eta(r,\delta)>0$ (see Section \ref{sec-FF}), we get 
\[ |W_1^{[2]}|\leq c_4 \mathcal{H}^{15}(\FF(W))\leq c_5 \|\delta\varphi\| |M^{[\cdot]}|.\] For the second inequality we used $\mathcal H^{15}(\FF(W))\leq c_6 \mathcal H^{15}(W)$ from Lemma \ref{lem-FFproperties}. We have $N\xRightarrow{N}N^{(14)}$ by Lemma \ref{lem-CodimDomination} (3) and $N^{(14)}\xRightarrow{W_1}\emptyset$ so $N\xRightarrow{N\cup W_1}\emptyset$ by Lemma \ref{lem-WdominationTrans}. Finally, by Lemma \ref{lem-cocyclefilling}, we find $\psi\in Z^1(M,H)$ such that $\psi=\varphi|_{M\setminus(N\cup W_1)}.$ In particular, \[\|\varphi\psi^{-1}\|\leq \|N^{[2]}\cup W_1^{[2]}\|\leq c_7 \|\delta\varphi\|.\]
This proves the Theorem under the assumption $\|\delta\varphi\|\leq \alpha_1.$ If this is not satisfied we just take $\psi=1$. In every case we have found $\psi\in Z^1(M,H)$ such that $\|\varphi\psi^{-1}\|\leq c_8 \|\delta\varphi\|, \text{ with } c_8:=\max\{c_7,\alpha_1^{-1}\}.$
\end{proof}

\appendix{} \label{sectiononsweepouts}
\section{Sweepouts}\label{sec-sweepouts}

In this section we collect some facts about sweepouts that are used in the main body of the paper.  First we define the space of mod-2 Lipschitz cycles in the flat norm.  Our presentation closely follows \cite{g09}. For a closed Riemannian manifold $M$, a mod-2 Lipschitz k-chain in $M$ is a sum $\sum_{i=1}^k a_i f_i$ where $a_i \in \mathbb{F}_2$ and each $f_i$ is a Lipschitz map from the standard k-simplex to $M$.  The space of all mod-2 Lipschitz k-chains is naturally a vector space over $\mathbb{F}_2$.  The boundary of a mod-2 Lipschitz chain is defined as in singular homology.  We denote the set of  mod-2 Lipschitz k-chains that lie in the kernel of the boundary map-- the space of mod-2 Lipschitz \textit{k-cycles}-- by $Z^*(M,k)$.  Since Lipschitz functions are differentiable almost everywhere, one can define the volume of an element of $Z(M,k)$ in the natural way. 

Let $Z_0^*(M,k)$ be the subset of $Z^*(M,k)$ consisting of cycles that are trivial as elements of singular homology.  We define the flat metric $d_{\mathbb{F}}$ on $Z_0^*(M,k)$ as follows.  For $a,b \in Z_0^*(M,k)$, let $\mathcal{C}$ be the collection of all mod-2 Lipschitz $k+1$-chains $C$ so that the boundary of $C$ is $a-b$, and set 
\[
d_{\mathbb{F}}(a,b):= \inf_{C\in \mathcal{C}} \vol_{k+1}(C).  
\]
\noindent Note that since $a$ and $b$ are both null-homologous, the set $C$ is non-empty.  The flat distance between two different cycles $a$ and $b$ can be zero, if for example one of them contains $k$-cycles of zero volume.  We quotient $Z_0^*(M,k)$ by the equivalence relation of being at flat distance zero to obtain a metric space whose completion we denote by $Z_0(M,k)$.  We define the volume of an element $c$ of  $Z_0(M,k)$ to be 
\[
\inf \lim_{m \to \infty } \inf \vol_k(c_m),
\]
\noindent where the infimum is taken over all sequences $\{c_m\}$ in $Z_0^*(M,k)$ converging to $c$.

A k-sweepout $F:X \rightarrow Z_0(M,k)$ is a map from a polyhedral complex $X$ that is continuous in the flat topology.  For each sweepout $F$ we define the volume $V(F)$ of $F$ to be $\sup_{x\in X} \vol(F(x))$.  Almgren showed that $H^{n-k}(Z_0(M,k);\mathbb{F}_2)$ contains a non-zero element, which we call $a(M,k)$, that corresponds to the fundamental class of $M$ \cite{almgren1962homotopy}.  We say that a sweepout $F: X \rightarrow Z_0(M,k)$ is nontrivial if it detects $a(M,k)$ in the sense that $F^*(a(M,k))$ is nonzero in $H^{n-k}(X;\mathbb{F}_2)$.  We define the \textit{k-waist} of $M$ to be the infimal volume over all nontrivial k-sweepouts.   The following theorem due to Almgren will be important for us .  Later Gromov gave a shorter proof of a similar statement \cite[pg. 134]{gromov1983filling} (see also \cite[Proposition 2.1]{g09}.)   

\begin{thm} \label{positivewaist}
	Let $M$ be a closed Riemannian manifold of dimension $n$.  Then the k-waist of $M$ for $0<k<n$ is positive.  
	\end{thm}


\subsection{} 
We now describe a useful criterion for a sweepout to be non-trivial.  To do this we first need to introduce the notion of a complex of cycles, which will serve as a discrete approximation to a sweepout.  A complex of cycles $\mathcal{C}$ is parametrized by a polyhedral complex $X$.  It associates to each $i$-face $\Delta$ of $X$ an $i+k$ chain $\mathcal{C}(\Delta)$ in $M$.  These chains must satisfy the following compatiblity condition: if $A$ is an $i$-face in $X$ with $\partial A= \sum B_\ell$, then $\partial \mathcal{C}(A)= \sum \mathcal{C}(B_{\ell})$.  Note that this implies that vertices of $X$ are sent to cycles.    This also implies that a sum of $i$-dimensional faces in $X$ that adds up to a cycle gives an $i+k$ cycle in $X$, and so we get a map $H_i(X;\mathbb{F}_2) \rightarrow H_{i+k} (M;\mathbb{F}_2)$.  The case $i=n-k$ is the one that will matter for us.  

We describe how to build a complex of cycles $\mathcal{C}_F$ from a sweepout $F$.  Let $F:X \rightarrow Z_0(M,k)$  be a sweepout, and for some $\delta>0$ to be specified later choose a subdivision of the polyhedral complex structure on $X$ so that adjacent vertices $v_1$ and $v_2$ satisfy $d_{\mathbb{F}}(F(v_1)',F(v_2)')<3\delta$, for $F(v_1)'$ and $F(v_2)'$ Lipschitz k-chains in $Z_0(M,k)$ satisfying $d_{\mathbb{F}}(F(v_i),F(v_i)')<\delta$ and $\vol_k(F(v_i)') < \vol_k(F(v_i))+ \delta$.  (Recall that $Z_0(M,k)$ is the completion of the space of null-homologous Lipschitz k-chains relative to the flat metric.)  

For a vertex $v$ we define $\mathcal{C}_F(v)$ to be $F(v)'$.  For an edge $E$ joining vertices $v_1$ and $v_2$, we define $\mathcal{C}_F(E)$ to be a chain $C_E$ with $\partial C_E = F(v_1)' - F(v_2)'$ and $\vol_{k+1}(C_E) <5 \delta$, where to find $C_E$ we have used the Federer-Fleming isoperimetric inequality (see \cite[Section 3.4]{gromov1983filling})  Continuing in this way and assuming that we have defined $\mathcal{C}_F$ on all $i-1$-faces, for each i-face $A$ of $X$ we define $\mathcal{C}_F(A)$ to be a chain $C_A$ with boundary equal to 
\[
\sum_{B \in \partial A} \mathcal{C}_F(B)
\]
\noindent and $\vol_{i+k}(C_A) < C(n) \delta$, where we have again used the Federer-Fleming isoperimetric inequality.  

Having defined $\mathcal{C}_F$, we obtain a map $H_{n-k}(X;\mathbb{F}_2) \rightarrow H_{n} (M;\mathbb{F}_2)$.  Provided we choose $\delta$ small enough, this map will be well defined independent of the complex of cycles we chose \cite{g09}.  The following criterion for $F$ to be a nontrivial sweepout makes precise the sense in which $a(M,k)\in H^{n-k}(Z(M,k);\mathbb{F}_2)$ corresponds to the fundamental class of $M$.  

\begin{prop} 
The sweepout $F:X \rightarrow Z(M,k)$ is nontrivial exactly if the image of the map $H_{n-k}(X;\mathbb{F}_2) \rightarrow H_{n} (M;\mathbb{F}_2) \cong \mathbb{F}_2$ obtained from $\mathcal{C}_F$ is nontrivial.  
\end{prop}

\subsection{}

We now explain why generic piecewise-linear maps to lower-dimensional Euclidean spaces give sweepouts.  Let $M$ be a closed smooth  manifold of dimension $n$, and let $\Phi :M \rightarrow \mathbb{R}^d$ be a piecewise linear map for $d<n$ for some triangulation of $M$.  We also require $\Phi$ to be generic, which recall means that it maps any $d+1$ vertices of a $d$-dimensional face of the triangulation of $M$ to a set of points in $\mathbb{R}^d$ that are not all contained in an affine hyperplane.  This can be arranged by an arbitrarily small perturbation of $\Phi$, by perturbing the images of the vertices and extending piecewise linearly.  Because $\Phi$ is generic, the preimage $\Phi^{-1}(p)$ of any point $p \in \mathbb{R}^d$ is a mod-2 Lipschitz $(n-d)$-cycle.  We see that it is a null-homologous cycle, and thus an element of $Z_0(M, n-d)$, by noting that the image of $\Phi$ is contained in a compact set, choosing a path joining $p$ to any point outside of this compact set, and taking the preimage under $\Phi$.  In more detail, one can choose the path $\gamma$ so that it intersects the positive-codimension faces of the triangulation of $\mathbb{R}^d$ in a finite set of points $p_1,..,p_l$. Since $\Phi$ defines a fiber bundle over the sub-interval of $\gamma$ between $p_i$ and $p_{i+1}$, it is clear that the preimages of points on this sub-interval are homologous.  One can also check that the preimage of a sub-interval $I$ of $\gamma$ containing one of the $p_i$ defines a  chain that bounds the preimages of the two endpoints of $I$.  Therefore all fibers are null-homologous.

The map $\Phi$ thus defines a continuous map $F_\Phi: \mathbb S^d \cong \mathbb{R}^d \cup \infty \rightarrow Z_0(M, n-d)$, where $F_{\Phi}$ sends the point at infinity to the empty cycle (as it does every point outside of a compact set containing the image of $F$.)  The following proposition seems to be well-known, at least if we were to replace ``piecewise-linear" with ``smooth," but we were unable to find a reference.  

\begin{prop} \label{prop:sweepout}
	For a generic piecewise-linear map $\Psi$, the map $F_{\Psi}$ is a nontrivial sweepout.  
	\end{prop} 
\begin{proof} 
	We prove this by induction on $d$.  For $d=1$, $\Psi$ is a generic PL map $M \rightarrow \mathbb{R}$.  For any $\delta>0$, we can choose a set of points $x_0<x_1<..<x_n$ sufficiently close together so that $d_{\mathbb{F}}(\Psi^{-1}(x_i), \Psi^{-1}(x_{i+1}))< \delta$.  We also require that the $\Psi$-preimages of both $x_0$ and $x_n$ are empty.  In building the complex of cycles $\mathcal{C}$ for $F_{\Psi}$ for the simplicial complex structure on $\mathbb{R}$ given by the $x_i$, we can set $\mathcal{C}(x_i)$ equal to $\Psi^{-1}(x_i)$.  The preimage $\Psi^{-1}((x_i,x_{i+1}))$ will be an open submanifold whose boundary is equal to the $n-1$-chain
	\[
	\Psi^{-1}(x_{i-1}) \cup \Psi^{-1}(x_i), 
	\]
	and so we can set $\mathcal{C}([x_{i-1},x_i])$ equal to $\Psi^{-1}([x_i,x_{i+1}])$.  The $\mathcal{C}([x_{i-1},x_i])$ sum to a cycle representing the mod-2 fundamental class, which finishes the $d=1$ case. 

	Assume the proposition is true for $d-1$, and let $\Psi: M \rightarrow \mathbb{R}^d$ be a generic PL map.  If $\pi: \mathbb{R}^d \rightarrow \mathbb{R}^{d-1}$ is the projection onto a generic hyperplane, then $\pi \circ \Psi$ is also generic.  Let $\mathcal{C}_\pi$ be a complex of cycles for $\pi \circ \Psi$ and some simplicial complex structure $\Delta$  on $\mathbb{R}^{d-1} \cup \infty$ comprised of simplices with edge lengths less than $\delta$ for $\delta$ sufficiently small, where $\mathcal{C}_\pi$ assigns to each face in $\Delta$ the preimage under $\pi \circ F$.

     \begin{figure}[h!] \label{figure:picomplexofcycles}
		\includegraphics[width=40mm, height=70mm]{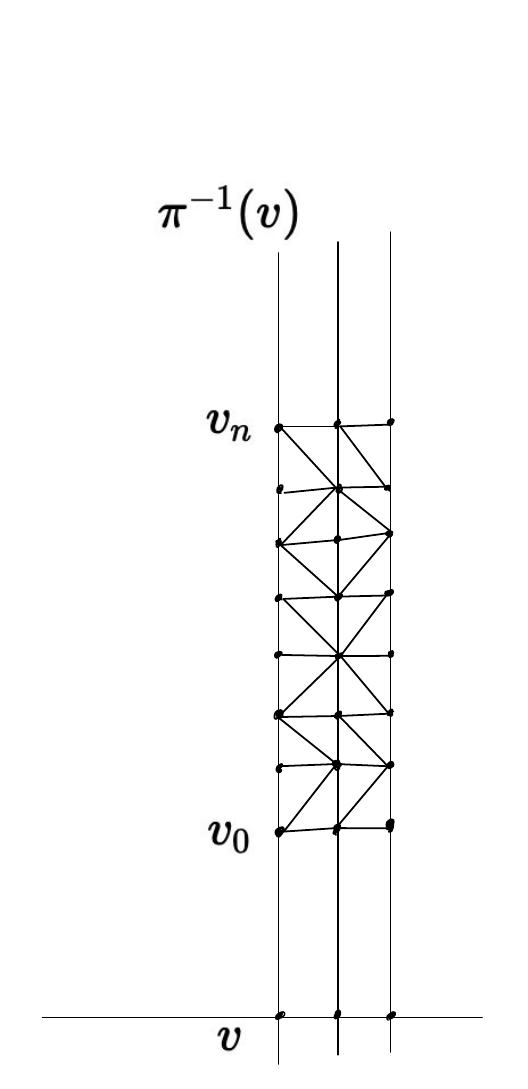} 
		\caption{The subdivided polyhedral complex structure for the map $\pi$ \iffalse from $\mathbb{R}^2$\fi  in the case $d=2$ } 
	\end{figure}	
	

	By taking the preimage of a large ball in $\mathbb{R}^{d-1}$ under $\pi$, we obtain a polyhedral complex structure on $\mathbb{R}^d \cup \infty $ a sub-complex of which containing the image of $\Psi$ is isomorphic to the product of a sub-complex $\Delta'$ of $\Delta$ with $\mathbb{R}$.  We take a subdivision of $\Delta' \times \mathbb{R}$ so that for each vertex $v\in \Delta$, the preimage $\pi^{-1}(v)$ contains vertices $v_0,v_1,..,v_{m} $ of the subdivision, listed in order of increasing $\mathbb{R}$-component, so that $d(v_i,v_{i-1})< \delta$.

 We claim that we can build a complex of cycles $\mathcal{C}$ for $F_\Psi$ such that the sum of the mod-2 $n$-chains of $\mathcal{C}$ is equal to the sum of the mod-2 $n$-chains of  $\mathcal{C}_{\pi}$.  We do this inductively in the following way.  To the vertices of the polyhedral structure on $\mathbb{R}^d$ we assign the preimages under $F$.  For edges $[v,w]$   we set $\mathcal{C}([v,w])$ equal to $F^{-1}([v,w])$, for $[v,w]$ the straight-line segment joining the vertices $v$ and $w$.  Note that the sum of the d+1-chains in $\mathcal{C}$ is equal to the sum of the d+1-chains in $\mathcal{C}_\pi$ plus the sum of the d+1-chains in $\mathcal{C}$ corresponding to edges that are not contained in the preimage of a point under $\pi$.   
	
	 Continuing in this way, we have that the sum of the $i$-chains of $\mathcal{C}$ is equal to the sum of the $i$ chains of $\mathcal{C}_{\pi}$ plus the sum of the $\Psi$-preimages of $(i-d)$-faces of the simplicial complex structure on $\mathbb{R}^d$ that are not contained in the preimage of any $(i-d-1)$-face  of $\Delta$ under $\pi$.  Since every $(n-d)$-face in $\mathbb{R}^d$ is contained in the preimage of an $(n-d-1)$-face of $\Delta$, we have that the sum of the n-chains  in $\mathcal{C}$ is equal to the sum of the n-chains in $\mathcal{C}_\pi$.  Therefore if $\pi \circ \Psi$ defines a non-trivial sweepout, then $\Psi$ defines a nontrivial sweepout, which completes the induction.  

  \end{proof}

The following theorem is a consequence of \cite[Theorem 4.10]{pittsbook} and the regularity theory for almost minimizing varifolds contained in \cite{pittsbook}.  

\begin{thm}
    There is a stationary integral rectifiable k-varifold in $M$ whose area is equal to the k-waist of $M$.   
\end{thm}

\bibliography{bibliography}
\bibliographystyle{alpha}

\end{document}